\documentclass[a4paper,reqno,12pt]{amsart}

\usepackage{amsthm,amsfonts,amsxtra,amssymb,amscd}

\usepackage[all]{xy}

\usepackage{enumerate}

\usepackage{fullpage}

\usepackage{xspace}
\usepackage{color}
\definecolor{mycolor}{rgb}{0,0,0}
\definecolor{mycolorG}{rgb}{0,0,0}
\definecolor{mycolorP}{rgb}{0,0,0}

\newcommand{\tPhi}     {\tilde \Phi}
\newcommand{\tgra}     {\tilde \gra}
\newcommand{\tgrD}     {\tilde \grD}
\newcommand{\bH}       {\bar H}

\newcommand{\NuG}      { {^e\!{G}}}
\newcommand{\Nug}      { {^e\!{\gog}}}

\newcommand{\Nut}      { {^e\!{\got}}}
\newcommand{\Nugos}    { {^e\!{\gos}}}
\newcommand{\Nub}      { {^e\!{\gob}}}
\newcommand{\NuB}      { {^e\!\!{B}}}

\newcommand{\NuP}      { {^e\!\!{P}}}

\newcommand{\Nun}      { {^e\!{\gon}}}
\newcommand{\NuU}      { {^e\!{\calU}}}
\newcommand{\Nusun}    { {^e\!{\hat{\gon}}}}
\newcommand{\NusuU}    { {^e\!{\hat{\calU}}}}
\newcommand{\NuPhi}    { {^e\!{\Phi}}}
\newcommand{\NutPhi}   { {^e\!{\tilde{\Phi}}}}
\newcommand{\NugrD}    { {^e\!\!{\grD}}}
\newcommand{\NutgrD}   { {^e\!\!{\tilde{\grD}}}}
\newcommand{\NuW}      { {^e\!{W}}}
\newcommand{\NutW}     { {^e\!{\wt{W}}}}

\newcommand{\NuE}      { {^e\!\!{E}}}
\newcommand{\NutE}     { {^e\!\!{\tilde{E}}}}
\newcommand{\NuO}      { {^e\!\Omega}}
\newcommand{\nugre}    { {^e\!{\gre}}}
\newcommand{\nuom}     { {^e\!{\omega}}}
\newcommand{\nupi}     { {^e\!{\pi}}}
\newcommand{\nueta}     { {^e\!{\eta}}}
\newcommand{\nutom}    { {^e\!{\wt{\omega}}}}

\newcommand{\Nugr}     { {^e\!{gr}}}

\newcommand{\NuCartan} {\tilde{\mathbf A}}
\newcommand{\Nur}      {\boldsymbol{r}}
\newcommand{\Nugrs}    {\boldsymbol{\grs}}
\newcommand{\Nukappa}  {\boldsymbol{\kappa}}

\newcommand{\sueta}    {\hat{\eta}}
\newcommand{\sutau}    {\hat{\tau}}
\newcommand{\extended} {extended }

\newcommand{\alt}      {hgt}
\newcommand{\suA}      {\hat{A}}
\newcommand{\Gr}       {{\calG\!r}}
\newcommand{\Rich}     {\calR}
\newcommand{\Schub}    {\calS}
\newcommand{\Xbar}{{\bar{X}}}
\newcommand{\tw}{\tilde w}
\newcommand{\meno}{\text{-}}

\theoremstyle{plain}
\newtheorem{lem}{Lemma}
\newtheorem{teo}[lem]{Theorem}
\newtheorem*{teo*}{Theorem}
\newtheorem{prp}[lem]{Proposition}
\newtheorem{cor}[lem]{Corollary}
\newtheorem{con}[lem]{Conjecture}

\theoremstyle{definition}
\newtheorem{dfn}[lem]{Definition}

\newtheorem{dfnteo}[lem]{Definition-Theorem}
\newtheorem{oss}[lem]{Remark}

\newtheorem{conv}[lem]{Convention}

\theoremstyle{remark}

\newcommand{\mA}  {\mathbb A} 
\newcommand{\mB}{{\mathbb B}} \newcommand{\mC}{\mathbb C} 
 \newcommand{\mF}{\mathbb F} \newcommand{\mG}{\mathbb G}
  
  \newcommand{\mM}{\mathbb M}
\newcommand{\mN}{\mathbb N}  \newcommand{\mP}{\mathbb P}
 \newcommand{\mR}{\mathbb R} \newcommand{\mS}{\mathbb S}
  \newcommand{\mV}{\mathbb V}
  
\newcommand{\mZ}{\mathbb Z}

\newcommand{\calB}{\mathcal B}  
 \newcommand{\calF}{\mathcal F} \newcommand{\calG}{\mathcal G}
  
 \newcommand{\calL}{\mathcal L} \newcommand{\calM}{\mathcal M}
 \newcommand{\calO}{\mathcal O} 
 \newcommand{\calR}{\mathcal R} \newcommand{\calS}{\mathcal S}
\newcommand{\calT}{\mathcal T} \newcommand{\calU}{\mathcal U}

\newcommand{\gob}{\mathfrak b}  
  \newcommand{\gog}{\mathfrak g}
\newcommand{\goh}{\mathfrak h}  
 \newcommand{\gol}{\mathfrak l} 
\newcommand{\gon}{\mathfrak n} \newcommand{\goo}{\mathfrak o} \newcommand{\gop}{\mathfrak p}
  \newcommand{\gos}{\mathfrak s}
\newcommand{\got}{\mathfrak t}

\newcommand{\sfA}  {\mathsf A}
\newcommand{\sfB}{\mathsf B} \newcommand{\sfC}{\mathsf C} \newcommand{\sfD}{\mathsf D}
\newcommand{\sfE}{\mathsf E} \newcommand{\sfF}{\mathsf F} \newcommand{\sfG}{\mathsf G}

  \newcommand{\sfS}{\mathsf S}

\newcommand{\sfe}{\mathsf e} \newcommand{\sff}{\mathsf f}

\newcommand{\gra}{\alpha}
\newcommand{\grb}{\beta}
\newcommand{\grg}{\gamma}
\newcommand{\grd}{\delta}
\newcommand{\gre}{\varepsilon}

\newcommand{\grl}{\lambda}
\newcommand{\grs}{\sigma}
\newcommand{\grf}{\varphi}

\newcommand{\om} {\omega}

\newcommand{\grG}{\Gamma}
\newcommand{\grD}{\Delta}
\newcommand{\grL}{\Lambda}

\newcommand{\mi}  {\imath}
\newcommand{\mj}  {\jmath}
\newcommand{\mk}  {\Bbbk}

\newcommand{\Aut}        {\operatorname{Aut}}

\renewcommand{\Im}       {\operatorname{Im}}

\newcommand{\card}       {\mathop{\mathnormal card}}
\newcommand{\Pic}        {\operatorname{Pic}}

\newcommand{\incluso} {\hookrightarrow}

\newcommand{\lra}     {\longrightarrow}
\newcommand{\isocan}  {\simeq}

\newcommand{\infinito}{\infty}
\newcommand{\cech}    {\spcheck}
\newcommand{\comp}    {\!\circ\!}

\newcommand{\coinc}   {\equiv}

\newcommand{\defi}    {:=}

\renewcommand{\geq}   {\geqslant}
\renewcommand{\leq}   {\leqslant}
\newcommand{\senza}   {\smallsetminus}
\newcommand{\ristretto}{\bigr|}

            \newcommand{\st}       {\, : \,}
       \newcommand{\mfor}     {\text{ for }}
         \newcommand{\mand}     {\text{ and }}
        \newcommand{\mif}      {\text{ if }}
  
  \newcommand{\mforall}  {\text{ for all }}

\newcommand{\wt}[1]    {\widetilde {#1} }

\newcommand{\iacc}   {\`{\i} }

\newcommand{\id}         {\mathrm{id}} 
\newcommand{\blU}        {\mathbf U}

\newcommand{\Q}{\mathbb{Q}}  
\newcommand{\B}{\mathbb{B}}

\newcommand{\al}{\alpha}   
   
\newcommand{\ep}{\epsilon} 
     \newcommand{\la}{\lambda}
     \newcommand{\De}{\Delta}
\newcommand{\La}{\Lambda}  \newcommand{\Ga}{\Gamma}

\newcommand{\Si}{\Sigma}

 \newcommand{\tal}{\wt\gra}
\newcommand{\tom}{\wt\om}              

 \newcommand{\osi}{\overline{\grs}}

\newcommand{\tDe}{\wt{\De}}

\newcommand{\tW}{\wt{W}}

\newcommand{\ts}{\tilde{s}}
\newcommand{\tbe}{\wt{\grb}}

\newcommand{\tlb}{\mathcal{O}(1)}
\newcommand{\lb}{\mathcal{L}}

\newcommand{\Pj}{\mathbb{P}}

\newcommand{\lih}{\mathfrak{h}}

\newcommand{\Res} {\mathop{\rm Res}\nolimits}
\newcommand{\sfBC}{\mathsf{BC}}

\renewcommand{\iacc}{\`\i\xspace}

\newcommand{\provf}{\sff}
\newcommand{\prove}{\sfe}

\newcommand{\myRank}{n}

\begin{document}

\title{Equations defining symmetric varieties and affine Grassmannians}

\author{Rocco Chiriv\iacc, Peter Littelmann and Andrea Maffei}



\maketitle

\begin{abstract}
Let $\grs$ be a simple involution of an algebraic semisimple group $G$ and let $H$ be the subgroup
of $G$ of points fixed by $\grs$. If the restricted root system
is of type $\sfA, \sfC$ or $\sfB\sfC$ and
$G$ is simply connected
or if the restricted root system is of type $\sfB$ and $G$ is adjoint, then we
describe a standard monomial theory and the equations for the coordinate ring $\mk[G/H]$ using the
standard monomial theory and the Pl\"ucker relations of an appropriate (maybe
infinite dimensional) Grassmann variety.
\end{abstract}
The aim of this paper is the description of the coordinate ring of the symmetric varieties and of
certain rings related to their wonderful compactification. The main tool to achieve this goal is a
(possibly infinite dimensional) Grassmann variety associated to a pair consisting of a symmetric
space and a spherical representation.

More precisely, let $G$ be a semisimple algebraic group over an algebraically closed field $\mk$ of
characteristic $0$ and let $\grs$ be a simple involution of $G$ (i.e. $G\rtimes \{\id,\grs\}$ acts
irreducibly on the Lie algebra of $G$). Let $H=G^\sigma$ be the fixed point subgroup. The quotient
$G/H$ is an affine variety, called a \emph{symmetric variety}.

A simple finite dimensional $G$-module $V$ is called spherical (for $H$) if $V^H\not=0$.
By results of Helgason \cite{Helg} and Vust \cite{Vust}, these modules are
parametrized by a submonoid $\Omega^+$ of the dominant weights of a suitable root system,
called the \emph{restricted root system}. As a $G$-module, $\mk[G/H]$ is well understood:
it is the direct sum $\bigoplus_{V\text{spherical}} V^*$.

Fix a spherical dominant weight $\gre$ in $\Omega^+$. We add a node $n_0$ to the Dynkin diagram of
$G$ and, for all simple roots $\alpha$, we join $n_0$ with the node $n_\alpha$ of the simple root
$\alpha$ by $\gre(\alpha^\vee)$ lines, and we put an arrow in direction of $n_\alpha$ if
$\gre(\alpha^\vee)\ge 2$. In the cases relevant for us, the Kac-Moody group $\NuG$ associated to
the extended diagram will be of finite or affine type. Let $\calL$ be the ample generator of
$\Pic(\Gr)$ for the generalized Grassmann variety $\Gr=\NuG/\NuP$.
The homogeneous coordinate ring $\grG_\Gr=\bigoplus_{j\ge0}\grG(\Gr,\calL^j)$ is the quotient of
the symmetric algebra $\sfS(\grG(\Gr,\calL))$ by an ideal generated by quadratic relations, the
{\it generalized Pl\"ucker relations}.

Since our aim is to relate these Pl\"ucker relations to $\mk[G/H]$, we say that the monoid
$\Omega^+$ is \emph{quadratic} if (it is free and) its basis has the following property with
respect to the dominant order of the restricted root system: any element of $\Omega^+$ that is less
than the sum of two elements of the basis is the sum of at most two elements of the basis. In
\ref{quadraticlattice} we show that this condition is equavalent to: either the group $G$ is simply
connected and the restricted root system is of type $\sfA$, $\sfBC$ or $\sfC$, or the group $G$ is
of adjoint type and the restricted root system is of type $\sfB$.

To analyse the structure of $\mk[G/H]$, we construct a $G$-equivariant ring homomorphism
$\grf:\grG_{\Gr}\lra \mk[G/H]$. If $\NuG$ is of finite type, then the morphism is just the pull
back of a canonical $G$-equivariant map $G/H\rightarrow \Gr$. In the general case, the underlying
idea is the same, but the construction is more involved.

Roughly speaking, the main result of this paper can be formulated as follows:
\begin{quote}
if $\Omega^+$ is quadratic, then the defining relations for $\mk[G/H]$ can be obtained from the
defining relations of $\Gr$ and a standard monomial theory for $\mk[G/H]$ can be obtained from the
standard monomial theory of a suitable $G$-stable Richardson variety $\Rich$ of $\Gr$.
\end{quote}

Let us formulate the result more precisely. If $\Omega^+$ is free (for example if $G$ is simply
connected), then let $V_1,\ldots,V_\ell$ be the simple spherical modules corresponding to the basis
of $\Omega^+$; in this case a basis of $\mV^*=\bigoplus_{j=1}^\ell V_j^*$ is a canonical set of
generators for $\mk[G/H]$. We have a canonical surjective map from the symmetric algebra
$\Psi:\sfS(\mV^*)\rightarrow \mk[G/H]$, and our aim is to describe of the kernel $Rel$ of $\Psi$.

We need to recall a few facts about the generalized Pl\"ucker relations. In \cite{L:SMT}, a basis
$\mF\subset\grG(\Gr,\calL)$ has been constructed together with a partial order "$\ge$``, such that
the monomials $\mF^2=\{ff'\mid f,f'\in \mF,f\le f'\}\subset \grG(\Gr,\calL^{\otimes 2})$ form a
basis. For a pair $f,f'\in \mF$ of not comparable elements let $R_{f,f'}\in
\sfS^2(\grG(\Gr,\calL))$ be the relation expressing the product $ff'$ as a linear combination of
elements in $\mF^2$. It was shown in \cite{LLM} that the $R_{f,f'}$ generate the defining ideal of
$\Gr\hookrightarrow \mP(\grG(\Gr,\calL)^*)$.

If $\Omega^+$ is quadratic, then we can construct a $G$-equivariant injection
$i:\mV^*\hookrightarrow \grG(\Gr,\calL)$ such that $\grf\circ i:\mV^*\rightarrow \mk[G/H]$ is an
isomorphism onto the image and $i(\mV^*)$ is compatible with  $\mF$. I.e., there exists a (finite)
subset $\mF_0\subset \mF$ forming a basis for the image of $\mV^*$. For $f\in \mF_0$ set
$g_f=\grf(f)$, then
$\mG=\{g_f\mid f\in\mF_0\}$ 
is a generating set for $\mk[G/H]$.

The relations $R_{f,f'}$ for $f,f'\in \mF_0$ involve also elements in $\mF-\mF_0$. Let
$\mF_1\sqcup\mF_0$ be the (finite) set of functions appearing in some polynomial $R_{f,f'}$ for
$f,f'\in\mF_0$. Denote by $\hat R_{f,f'}\in \sfS^2(\mV^*)$ the relation obtained from  $R_{f,f'}$
by replacing a generator $h\in \mF_0$ by $g_h\in \mG$ and a generator $h\in \mF_1$ by the function
$F_h=\grf(h)$ of $\mG$.
\begin{teo*}
The relations $\{\hat R_{f,f'} \st f,f'\in \mF_0$ not comparable$\}$ generate the ideal $Rel$
of the relations among the generators $\mG$ of $\mk[G/H]$.
\end{teo*}

Now we want to give a more detailed description of the construction of the $G$-equivariant ring
homomorphism $\grf:\grG_{\Gr}\lra \mk[G/H]$.

Consider the generalized Grassmann variety $\Gr=\NuG/\NuP$, let $\calL$ be as above and continue to
assume that the monoid $\Omega^+$ is quadratic.
We show that there exist a $G$-stable Richardson variety $\Rich$ in $\Gr$ which has a homogenous
coordinate ring that looks, as a $G$-module, like $\mk[G/H]$, i.e.,
$\mk[G/H]\simeq_G\Gamma_\Rich=\bigoplus_{j\geq0}\grG(\Rich,\calL^j)$. Moreover
(see Corollary \ref{cor:sezioniR}):
$$
\grG(\Rich,\calL)= V_1^* \oplus \dots \oplus V_\ell^*.
$$
The basis $\mF$ of $\Gamma(\Gr,\calL)$ we have introduced above is compatible with
a certain Richardson
subvariety; furthermore the standard monomials of elements in the set $\mF_0$ are a basis of
$\Gamma_\Rich$.


To relate this standard monomial theory to the symmetric space, note that the Lie algebra of $\NuG$
has by construction a natural grading such that in degree $0$ there is a maximal torus and the Lie
algebra of $G$, and in degree $-1$ there is the $G$ spherical module $V$ of highest weight $\gre$.
In particular, there exists a vector $h_{-1}$ fixed by $H$ in degree $-1$. If $\NuG$ is of finite
type, then we may consider the exponential $e^{h_{\meno 1}}$ and the $H$ stable point
$x=e^{h_{\meno 1}}\NuP$ in the Grassmannian $\NuG/\NuP$. So in this case we are able to define a
$G$-equivariant map from $G/H$ to $\Gr$ using the map $gH\longmapsto g x$. The pull back of such
map gives the ring homomorphism $\grf:\grG_{\Gr}\lra \mk[G/H]$. In fact, the morphism
$\grf:\grG_{\Gr}\lra \mk[G/H]$ can also be defined when $\NuG$ is not of finite type (see Section
\ref{sez:equazioni}).

Moreover we are able to show that the previous theorem may be strengthened to

\begin{teo*}
Consider $\mG = \{g_f\ |\ f \in \mF_0\}$ as a partially ordered
set with the same partial order as on $\mF$.
Then $\mG$ is a basis of $\mV^*\subset \mk[G/H]$, the set $\mS\mM_0$ of ordered
monomials in $\mG$ realizes a standard monomial theory for $\mk[G/H]$ and the
relations $\hat R_{f,f'}$ for the non standard $ff'$ are a set of straightening relations.
\end{teo*}

A key point in the proof of the theorems above is Theorem \ref{teo:hnonsvanisce}, whose proof in
turn uses some results about the product in $\mk[G/H]$ from \cite{CM2} which hold only in
characteristic zero. We need this hypothesis of course also for the definition of $e^{h_{\meno
1}}$. However, we want to point out that in most of the cases where the restricted root system of
type $\sfA$, it is possible to directly define the point $x$. If one is able to check the
conclusions of Theorem~\ref{teo:hnonsvanisce} in these cases, then the corresponding result is
valid in arbitrary characteristic since the remaining arguments are characteristic free.

The standard monomial theory is compatible with the decomposition in $G$-modules in
the following sense:  there exists a filtration of $\mk[G/H]$ by $G$-modules $F_i$ with
simple quotients such that for all $i$ the set $
\mS\mM_0\cap F_i$ is a $\mk$-basis of $F_i$ (Remark \ref{oss:smt}).

%
%

We want to stress that the relations $\hat R$ describing the ideal $Rel$ cannot be considered as
completely explicit. The actual computation of the functions $F_f$ depends only on the exponential
$e^{h_{-1}}$ and on the representation theory of $G$ (see remark \ref{oss:polF}). Such computations
may be considered as algorithmic, but it seems very difficult to obtain more explicit formulas.
Clearly it should be interesting to have more information on such formulas.

If $\NuG$ is of finite type (or, equivalently, the restricted root system is of type $\sfA$) we can
show that $\mF_1$ is given by just two elements $f_0,f_1$ and that
$$ F_{f_0}=F_{f_1}=1. $$
In particular, in these cases the explicit relations may be summarized in the following
description of the coordinate ring of the symmetric variety:
$$ \mk[G/H]\isocan \frac{\grG_{\Gr}} { (f_{0}=f_{1}=1) }. $$

The study of the coordinate rings $\mk[G/H]$ is strongly related to the study of the multicone
associated to the wonderful compactification of the symmetric varieties of adjoint type. De Concini
and Procesi \cite{DP1} defined the wonderful compactification $\Xbar$ of $G/\bH$ where $\bH$ is the
normalizer of $H$. In \cite{CM1} the total ring of sections
$\grG=\oplus_{\calM\in\Pic(X)}\grG(\Xbar,\calM)$ and a canonical set of generators for these rings
had been introduced. The computation of the relations among these generators is equivalent to the
computation of the relations in the ring $\mk[G/H]$ above.

In some special cases a standard monomial theory for $\mk[G/H]$ had been developed before
\begin{enumerate}
\item[-] for $G/H=SL(n)$, corresponding to the involution $(x,y)\mapsto(y,x)$ of the group
$SL(n)\times SL(n)$ and whose restricted root system is of type $\sfA$, here our construction
gives the same as the construction of De Concini, Eisenbud and Procesi \cite{DEP};

\item[-] for $G/H=$`symmetric quadrics', corresponding to the involution $x\mapsto(x^{\meno 1})^t$
of the group $SL(n)$ and whose restricted root system is of type $\sfA$, a theory of standard
monomials has been introduced by Strickland \cite{Str} and Musili \cite{Musili,Musili2}; however,
we do not know whether their SMT is equivalent to ours;

\item[-]for $G/H=Sp(2\,n)$, corresponding to the involution $(x,y)\mapsto(y,x)$ of the group
$Sp(2\,n)\times Sp(2\,n)$ and whose restricted root system is of type $\sfC$, a theory of standard
monomials has been introduced by De Concini in \cite{DC:Sp}. Also in this case we do not know
whether this SMT is equivalent to ours.
\end{enumerate}
The results above cover almost all cases with restricted root system of type $\sfA$; there are only
two families missing whose restricted root system is of type $\sfA_1$ (and hence they are very
simple), the `symplectic quadrics' and an involution of $\sfE_6$ which we discuss briefly at the
end of the paper.

Finally we want to stress that the condition on the restricted root system to be of type $\sfA$,
$\sfB$, $\sfC$ or $\sfB\sfC$, while looking strong, is actually fulfilled for many involutions. In
the Tables in \cite{VO} it holds for $12$ families of involutions out of a total of $13$ families
and in $4$ exceptional cases out of a total of $12$. Moreover one should add to such list of
families the involutions such that $G=H\times H$, $H$ is simple and the involution is given by
$(x,y)\mapsto(y,x)$; for these cases $\mk[G/H]$ is the coordinate ring of $H$ and our condition is
equivalent to $H$ equals to $SL(n)$ or $Sp(2n)$ or $SO(2n+1)$.

Now we want to describe the structure of the article. In the first section below we introduce
notation and gives some preliminary result on the comninatorics of the set of spherical weights.

In Section \ref{sez:sezionicomplete} we review the main properties of the De Concini Procesi
wonderful compactification of a symmetric variety. We relate the multiplication of sections of line
bundles on such compactification and the multiplication of functions on the symmetric variety.


In Section \ref{sez:extended} we study some simple properties of the group $\NuG$. In the cases
related to our problem stated above, the group $\NuG$ is of finite type if and only if the
restricted root system is of type $\sfA$, and it is of affine type if and only if the restricted
root system is of type $\sfB$, $\sfB\sfC$, $\sfC$ or $\sfD$ (see Proposition
\ref{prp:finitoaffine}).

In Section \ref{sez:schubert} we introduce and study a certain module of the extended Lie algebra
corresponding to the new node of the extended Dynkin dagram. In the same section we study also the
Richardson variety $\Rich$.

In Section \ref{sez:equazioni} all results of the previous sections are used to relate the
symmetric variety and the Grassmannian $\Gr$. And in Section \ref{sez:finito} we study the simpler
situation where the Grassmannian $\Gr$ is finite dimensional.

In the Appendix we prove that two standard monomial bases related to the symmetric
variety coincide. One of the two bases is the one
considered above, the other is the standard monomial basis one may
construct via lifting and pull back from the standard monomial theory of the multicone over the
closed orbit in the wonderful compactification.

Finally, for the convenience of the reader, we have reported in the Appendix B the Satake diagrams
of the involutions together with the additional node relevant for the constructions and other
informations.

\section{The coordinate ring of $G/H$ and quadratic lattices}\label{sez:1}
In this section we introduce some {notation} and we make some remarks
on the combinatorics of spherical weights.

Let $G$ be a semisimple simply connected algebraic group over an
algebraically closed field of characteristic zero.  Let $\grs$ be an
involution of $G$ and $H_{sc}$ its fixed point subgroup.  Since $G$ is
simply connected $H_{sc}$ is known to be connected (see for example
\cite{VO}).

Let now $q:G\lra G_q$ be an isogeny and let $K_q$ be the kernel of $q$.
If $\grs (K_q) = K_q$, then we can consider an involution $\grs_q$ of $G_q$
induced by $\grs$ and its fixed points $G_q^{\grs_q}$. We define also
$H_q$ as the inverse image of $G_q^{\grs_q}$ in $G$.
The groups $H_q$ are reductive so the quotients
$X_q = X_q(\grs) = G_q / G_q^{\grs_q} = G / H_q$ are
affine varieties. These varieties are
called \emph{symmetric varieties}.
When $q$ is the identity, then we use the subscript $sc$ instead of $q=id$. We also
denote by $q=ad$ the adjoint quotient; in this case $H_{ad}$ is known to be equal to the normalizer
of $H_{sc}$ in $G$ (see \cite{DP1} \S 1).

\subsection{Spherical representations}
If $V$ is an irreducible representation of $G$, then we say that it is
$q$-\emph{spherical} (resp. spherical) if there exists a non zero
vector fixed by $H_q$ (resp. $H_{sc}$). The subspace $V^{H_q}$
of $H_q$-fixed vectors is then one dimensional, and hence
$$
\mk[X_q]=\mk[G]^{H_q} = \bigoplus_{V \: irr. \: rep.} V^*\otimes
V^{H_q} = \bigoplus_{V \:q\text{-}spherical} V^*.
$$
We want now to give a more precise description of the set of
$q$-spherical representations.

Let $T$ be a maximally split $\grs$ stable maximal torus of $G$, that is a maximal torus of $G$
stable under $\grs$ such that the dimension of $\{t\in T \st \grs(t)=t^{\meno 1}\}$ is maximal, and
let $S$ be the identity component of this subgroup. The dimension of $S$ is called the \emph{rank}
of the symmetric variety $G/H$ and we denote it by $\ell$.  Let $\grL$ be the weight lattice of $T$
and let $\grL_q$ be the sublattice of weights trivial on $K_q$. The Killing form $\kappa$ defines a
positive definite bilinear form on $\grL$ and on $\grL_q$. A weight $\grl$ is said to be special if
$\grs(\grl) = - \grl$ and we denote by $\grL^s$ (resp. $\grL_q^s$) the sublattice of $\grL$
(resp. $\grL_q$) of special weights.

Denote by $\Phi \subset \grL$ the set of roots. We choose the set
of positive roots $\Phi^+$ in such a way that if $\gra$ is positive,
then $\grs(\gra)$ is either equal to $\gra$ or is a negative root
(see \cite{DP1} \S 1).
We denote by $\grD$ the set of simple roots of $\Phi$ defined by the
choice of $\Phi^+$. In exactly the same way let
$\grL^+\subset \grL$ be the monoid of dominant weights.
If $\gra \in \Phi$ is not fixed by $\grs$, then we define the \emph{restricted
  root} $\tgra$ as $\gra - \grs(\gra)$ and the \emph{restricted root
  system} $\tPhi \subset \grL^s$ as the set of all restricted roots.
This is a (not necessarily reduced) root system (see \cite{Ric}) of
rank $\ell$ and the subset $\tPhi^+$ (resp. $\tDe$) of restricted
roots $\tgra$ with $\gra$ positive (resp. $\gra$ simple) is a choice
of positive roots (resp. a basis of simple roots) for $\tPhi$.
For $\tgra \in \tPhi$
we define $\tgra \cech \in \got$ such that
$\langle \tgra\cech , \grl \rangle =
2\kappa(\grl,\tgra)/\kappa(\tgra,\tgra)$ for all
$\lambda \in \got^*$.
A special weight $\grl\in \grL^s$ is said to be \emph{spherical} if
$\langle \tgra\cech , \grl \rangle \in \mZ$ for all $\tgra \in
\tPhi$. The subset $\Omega_{sc}\subset\grL^s$ of spherical weights
is a weight lattice for $\tPhi$ (w.r.t. $\kappa$) and we observe
that if $\grl \in \Omega_{sc}$ then it is dominant with respect to $\Phi^+$
if and only if it is dominant with respect to $\tPhi^+$.  On $\Omega_{sc}$
one has two different dominant orders:
one with respect to $\Phi^+$ that we indicate {by}
$\leq$, and one {with} respect to $\tPhi^+$ that we indicate by
$\leq_\grs$: if $\grl, \mu \in \Omega_{sc}$, then $\mu
\leq_\grs \grl$ iff $\grl -\mu \in \mN\,[\tPhi^+]$.

We can now describe the set of spherical representations.  For $\grl
\in \Lambda^+$ let $V_\grl$ be the irreducible representation of
$G$ of highest weight $\grl$. Define the set
$$
\Omega_q^+ := \{\grl \in \grL^+ \st V_\grl \text{ is
  $q$-spherical}\}
$$
and let $\Omega_q$ be the lattice generated by $\Omega_q^+$. If $\grl
\in \Omega_{sc}^+$, then we denote by $h_\grl \in V_\grl$ a non zero vector fixed
by $H_{sc}$. Given $\grl, \mu \in \Omega_q$, we can think of $V_\grl, V_{\mu}$
as sections of a line bundles over the flag variety of $G$, and hence the product
of the sections $h_\grl \cdot h_\mu$ is a nonzero vector fixed by $H_q$ in
$V_{\grl +\mu}$. In particular, we see that $\Omega_q^+$ is a monoid.
In the simply connected case this definition of $\Omega_{sc}$ coincides
with the one given above using the restricted root system ({see
Helgason \cite{Helg}}).  In general the set of $q$-spherical weights
has been characterized by Vust \cite{Vust} who proved the following Theorem.

\begin{teo}[Vust \cite{Vust} {Th\'eor\`eme} 3]
  Let $S_q = q(S)$ and let $\grL(S_q)$ be the weight lattice of
  $S_q$ and let $\grl \in \grL^+_q$. Then $\grl \in \Omega^+_q$ if and
  only if $\grs(\grl) = - \grl$ and $\grl \ristretto _{S_q} \in
  2\grL(S_q)$.
\end{teo}

In the following corollary we collect some consequences of the characterization
by Helgason and Vust.
\begin{cor}\hfill
\begin{enumerate}[\indent i)]
\item For every $q$ we have $\Omega_q^+ = \Omega_q \cap \grL^+$;
\item For every $q$ we have $\Omega_q = \{\grl-\grs(\grl) \st \grl \in \grL_q\}$;
\item For every $q$ we have $\grL_q\cap\Omega \supset \Omega_q \supset
  \mZ[\tPhi]$;
\item In the adjoint case we have $\Omega_{ad} = \mZ[\tPhi]$;
\item If $K_q \subset K_{q'}$ then the natural map $G/H_q \lra
  G/H_{q'}$ is an isomorphism if and only if $\Omega_q = \Omega_{q'}$.
\end{enumerate}
\end{cor}

\begin{proof}
In the simply connected case, the statements are part of the results of Helgason.

The condition given by Vust's criterion is linear, so $i)$ follows by Vust's characterization.
In particular, if $\grl \in \grL_q$, then $\grl \in \Omega_q$ if and only if $\grs(\grl)= - \grl$
and $\grl \ristretto _{S_q} \in 2\grL(S_q)$.

  Let $L = \{\grl-\grs(\grl) \st \grl \in \grL_q\}$. The inclusion $L
  \subset \Omega_q$ is evident.  To prove the converse notice that the
  restriction map $\rho:\grL_q \lra \grL(S_q)$ is surjective, {
  $\rho \ristretto _{\grL^s}$ is injective and $\rho \comp \grs =
  - \rho$}. Let now $\mu{ \in}\Omega_q$ and consider $\rho(\mu)$. By
  Vust's criterion there exists $\grl \in \grL_q$ such that
  $2\rho(\grl)=  \rho(\mu)$. Now $2\rho(\grl)= \rho(\grl - \grs(\grl))$
  and $\mu, \grl - \grs(\grl) \in \grL^s$ so $\mu = \grl - \grs(\grl)$.
  This proves $ii)$.

Point $iv)$ follows directly from $ii)$, and $iii)$ follows from $iv)$
and $ii)$. Finally, $v)$ is an obvious consequence of the description of the
coordinate ring of $X_q$ given above.
\end{proof}

\medskip

\subsection{Quadratic lattices}\label{quadraticlattice}
As explained in the introduction, our construction of a
standard monomial theory for $G/H_q$ starts with a choice of canonical
generators of the coordinate ring. For this reason we want $\Omega^+ _q$
to be freely generated. The following combinatorial conditions will ensure
in addition that the relations between these generators are going to be quadratic.

\begin{dfn}\label{dfn:reticoli}
Let $R$ be a root system with a choice of positive roots $R^+$, let $P$ be
the weight lattice with $P^+$ as the monoid of dominant weights, and let
$Q\subseteq P$ be the root lattice.  For a sublattice $L\subseteq P$ set
$L^+ = L \cap P^+$. The sublattice $L$ is called \emph{admissible} if
\begin{enumerate}[\indent i)]
\item $L\supset Q$;
\item $L^+$ is a finitely generated free commutative monoid.
\end{enumerate}
The (unique) basis $\calB$ of the free monoid $L^+$ (note that $\calB$
is also a basis of $L$) is called the
\emph{admissible basis} of $L$.  If $\grl \in L$ and
$\grl= \sum_{\gre \in \calB} a_\gre \, \gre$, then we define $\alt_\calB (\grl)=
\sum_{\gre \in \calB} a_\gre$. An admissible lattice $L$ is called
\emph{quadratic} if the following additional property holds:
\begin{enumerate}[iii)]
\item If $\grl \in L^+$ is such that $\grl\le \gre+\eta$ for some $\gre,\eta \in \calB$
(with respect to the dominant order), then $\alt_\calB (\grl) \leq 2$.
\end{enumerate}
\end{dfn}

This definition is strongly related to the description of the coordinate
ring of $G/H_q$: take $R=\tPhi$ and suppose that $\Omega_q$ is
admissible. Let $\calB = \{\gre_1,\dots,\gre_\ell\}$,
then fixing a basis of $V^*_{\gre_1}\oplus \dots \oplus V^*_{\gre_\ell}$ is a canonical
choice for fixing a generating set of $\mC[X_q]$.  A rough description of the
relations between the generators is given in the next section.
From that description it will be clear that if
$\Omega_q$ is quadratic, then also the relations in these generators
are quadratic (see Corollary \ref{cor:quadratic}).

\begin{conv}\label{conv:fundBC}
Before proving the next proposition we fix a convention for the fundamental weights
of a root system $\Phi$ of type $\sfB\sfC_\ell$. Let $\gra_1,\dots,\gra_\ell$ be simple roots
of $\Phi$ such that $2\gra_\ell \in \Phi$. Notice that $\gra_\ell\cech = 2(2\gra_\ell)\cech$.
We define
the fundamental weights $\om_1,\dots,\om_\ell$ as the weights such that
$\langle \om_i , \gra_j\cech \rangle =
\grd_{ij}$ if $i \neq \ell$ or $j\neq \ell$ and $\langle \om_\ell , \gra_\ell\cech \rangle=2$.
With this definition $\{\om_1,\dots,\om_\ell\}$ is a basis of the weight lattice.
\end{conv}

Now we classify the quadratic lattices for an abstract root system.

\begin{prp}\label{prp:quadratici}
  Let $R, R^+, P, Q$ be as in Definition~\ref{dfn:reticoli} above with $R$ simple.
Then a lattice $L\subset P$ is quadratic only in the following cases
\begin{enumerate}[\indent i)]
\item $R$ is of type $\sfA_1$ and $L=P$ or $L=Q$;
\item $R$ is of type $\sfB\sfC_1$ and $L=P$;
\item $R$ is of type $\sfA_\ell, \sfC_\ell$ or $\sfB\sfC_\ell$ with
  $\ell \geq 2$ and $L=P$;
\item $R$ is of type $\sfB_\ell$ with $\ell \geq 2$ and $L=Q$.
\end{enumerate}
\end{prp}

\begin{proof}
  Let $\om_1,\dots,\om_\ell$ be the fundamental weights, $\gra_1,
  \dots, \gra_\ell$ the simple roots and let $n = \card{P/Q}$.
  For a quadratic lattice $L$ let $\calB$ be as in definition
  \ref{dfn:reticoli}. Notice that $L$ is a lattice of rank $\ell$
  so let $\calB=\{\gre_1, \dots , \gre_\ell \}$. By condition i)
  we know $n\,\omega_i \in L^+$ for all $i$, and hence by condition ii) the
  $\gre_i$ have to be multiples of the fundamental weights. So, up to
  renumbering them, we have $\gre_i = c_i \om_i$ for some $c_i \in \mN$.

 So given a simple root $\gra_i = \sum_j  c_{ij} \omega_j$, then
 $c_{ij} \omega_j \in L$.  In particular, if $R$ is of type $\sfA_\ell$
 ($\ell \geq 2$), $\sfC_\ell$ ($\ell \geq 3$), $D_\ell$ or $E_\ell$,
 for every $i$ there exists $j$ such that $c_{ij} =-1$, and hence
 $L=P$ in these cases.  In the cases $\sfB\sfC_\ell$, $\sfF_4$
 and $\sfG_2$ we have $n=1$, so $L=P=Q$.

The condition for $L$ to be quadratic is obviously equivalent to
$\alt_\calB (\gra) \geq 0$ for all simple roots $\gra$. So
if $L=P$ and the root system is of type
$\sfA_\ell$,  $\sfC_\ell$, $\sfB\sfC_\ell$, then the condition is satisfied; for
the root systems of type $\sfD_\ell$ or $\sfE_\ell$, if $\gra$ is the simple root
corresponding to the ramification node in the diagram, then we have
$\alt_\calB(\gra) <0$; and for  the root systems of type $\sfB_\ell$ with $\ell \geq 3$, $\sfG_2$
and $\sfF_4$, if $\gra$ is the simple long root ``near'' a short root, then we
have $\alt_\calB(\gra) <0$.

It remains to consider the cases $A_1$ and $\sfB_\ell$ with $L=Q$.
(Note in both cases $n=2$, so the only possibilities are $L=P$ or $L=Q$).
 For $A_1$ the proposition is trivially true, and for $B_\ell$ one has $Q=
\langle \om _1, \dots , \om_{l-1}, 2 \om_\ell \rangle$, again the fact that
the lattice is quadratic is easily verified.
\end{proof}

The proof shows also that the only admissible lattices $L$ for
which $L\neq P$ are the ones with $R=\sfB_\ell$ or $\sfA_1$ and $L=Q$.

Let $X_q$ be a symmetric variety such that $\Omega_q$ is quadratic.
In section 3 we will construct a group $\NuG$ with the properties briefly
explained in the introduction. In these cases the
restricted root system is always of type $\sfA$, $\sfB$, $\sfC$ or
$\sfB\sfC$. For convenience we introduce the following
convention that will be used in the next sections.

\begin{conv}\label{conv:basee}
Let $R$ be a simple root system of type $\sfA_\ell$, $\sfB_\ell$,
$\sfC_\ell$ or $\sfB\sfC_\ell$. Notice that a simple basis of $R$
is linearly ordered and we number it according to Bourbaki \cite{Bourbaki}.
In particular notice that we number in a different way the simple basis
of $\sfB_2$ and $\sfC_2$. Let $\om_1,\dots,\om_\ell$ be the fundamental weights and define
$$
\gre_i =
\begin{cases}
\om_i & \mif i \neq \ell \text{ or $R$ is not of type $\sfB_\ell$};\\
2\om_\ell  & \mif i = \ell \text{ and $R$ is of type $\sfB_\ell$}.
\end{cases}
$$
We refer to $\gre_1,\dots,\gre_\ell$ as the \emph{quadratic basis} since the lattice spanned by
$\gre_1,\dots,\gre_\ell$ is quadratic and all quadratic lattices (with the exception of $L=Q$ and
$R$ of type $\sfA_1$) of simple root systems are of this form. In order to have a uniform notation
we consider root systems of type $\sfA_1$ and $\sfB_1$ as different and we choose in the first case
$L=P$ and $\gre_1=\om_1$ and in the second case $L=Q$ and $\gre_1=2\om_1$.
\end{conv}

We will need later the following combinatorial lemma about basis of quadratic lattices.

\begin{lem} \label{lem:ord}
  Let $R$ be a simple root system of type $\sfA_\ell$, $\sfB_\ell$,
  $\sfC_\ell$ or $\sfB\sfC_\ell$ and let $\gre_1,\gre_2,\dots,\gre_\ell$
  be the quadratic basis according to Convention \ref{conv:basee}. Then for
  all $i= 1, \dots, \ell$ we have
\begin{enumerate}[\indent i)]
\item $\gre_i \leq \gre_1 + \gre_{i-1}$;
\item if $\grl, \mu \in P^+$ and $\mu \leq \gre_1$, $\grl \leq
  \gre_{i-1}$, $\gre_i \leq \grl+\mu$ then $\mu = \gre_1$, $\grl = \gre_{i-1}$.
\end{enumerate}
\end{lem}

\begin{proof}
$i)$ follows by $\gra_1+\dots+\gra_{i-1}=\gre_1+\gre_{i-1}-\gre_i$ for all $i$ and all types.

$ii)$ is trivial for $R$ of type $\sfA$. Assuming $R$ of type $\sfB$ we have that
$$0 < \gre_1 < \gre_2 < \dots < \gre_i$$
is a complete list of elements less or equal to $\gre_i$ for all $i$. So the statement follows from
$\gre_i \not \leq \gre_1+\gre_{i-2}$. For $R$ of type $\sfC$ we have that for all $i$
$$ \cdots \gre_{i-4} < \gre_{i-2} < \gre_i$$
is a complete list of elements less or equal to $\gre_i$. So the claim follows from $\gre_i \not
\leq \gre_1+\gre_{i-3}$.

Finally notice for this problem the arguments for $R$ of type $\sfB\sfC$
are the same as in the case $R$ of type $\sfB$.
\end{proof}

\section{The ring of sections of a complete symmetric variety}\label{sez:sezionicomplete}
In this section we recall some facts about the wonderful compactification of a symmetric variety of
adjoint type defined by De Concini and Procesi in \cite{DP1}. We describe the relation between the
multiplication of sections of line bundles on this compactification and the multiplication of
functions on the symmetric variety.

\subsection{The wonderful compactification of a symmetric variety}
We keep the notation introduced in the previous section, in particular,
$\ell$ is the rank of the lattice $\Omega_{sc}$. A spherical weight $\grl
\in \Omega_{sc}^+$ such that $\langle \tal \cech , \grl\rangle \neq 0$ for all
$\tal\in \tPhi$ is called regular. If $\grl$ is regular, then we have an
embedding of $X_{ad}=G/H_{ad} \incluso \mP(V_\grl)$ given by $[g]
\mapsto g[h_\grl]$. The wonderful compactification of De Concini and
Procesi of $X_{ad}$ is defined as the closure of this image and its
main properties are listed in the following theorem:

\begin{dfnteo}[Theorem 3.1 and Proposition 8.1 in \cite{DP1}]
  Up to isomorphism, the closure of $X_{ad}$ in $ \mP(V_\grl)$ does not depend on the
  choice of the spherical regular weight $\grl$. We call it the
  \emph{wonderful compactification} of $X_{ad}$ and we denote it by
  $\Xbar=\Xbar(\grs)$. This variety has the following properties:
\begin{enumerate}[\indent i)]
\item $\Xbar$ is a smooth projective $G$ variety;
\item $\Xbar\senza X_{ad}$ is a divisor with normal crossings and
  smooth irreducible components $S_1,\ldots,S_\ell$;
\item $\Xbar$ has a unique closed orbit $Y = Y(\grs)$ and the
  restriction of line bundles $\Pic(\Xbar) \lra \Pic(Y)$ is injective.
  In particular $\Pic(\Xbar)$ is identified with a sublattice of
  $\grL$ and we denote by $\calL_\grl$ the line bundle corresponding
  to a weight $\grl \in \Pic(\Xbar)$;
\item for every $\grl \in \Omega_{sc}^+$ (not necessarily regular) the map
  $[g] \mapsto g[h_\grl]$ from $X_{ad}$ to $\mP(V_\grl)$ extends to a
  morphism $\psi_\la:\Xbar\rightarrow\Pj(V_\la)$ and $\calL_\grl =
  \psi_\la^*\tlb$.
\end{enumerate}
\end{dfnteo}

By the properties iii) and iv) we know that $\Omega_{sc} \subset \Pic(\Xbar)$.
Moreover, the weights $\tal_1, \dots,\tal_\ell \in \tDe$ are the
weights corresponding to the line bundles
$\calO(S_1),\dots,\calO(S_\ell)$.  In particular, there exists a $G$
invariant section $s_i \in \grG(X,\calL_{\tal_i})$ such that $div
(s_i) = S_i$.

For an element $\nu= \sum_{i=1}^{\ell} n_i {\wt \gra}_i \in
\mN[\tPhi]$ the multiplication by $s^\nu \defi \Pi_i s_i^{n_i}$ gives
a $G$ equivariant map from $\Gamma (\Xbar,\lb_{\la-\nu})$ to $\Gamma
(\Xbar,\lb_\la)$.

Now we can describe the sections of a line bundle as a $G$-module.
Observe that every line bundle $\calL_\grl$ with $\grl \in \Omega_{sc}$ has
a natural $G$ linearization and{,} since the variety has a dense
orbit under the Borel subgroup, any irreducible $G$-module appears in
$\Gamma (\Xbar,\calL_\grl)$ with multiplicity at most one (see Lemma~8.2
in \cite{DP1}).

If $\mu \in \Omega_{sc}^+$, then by the construction of
$\calL_\mu$ we have a submodule of $\Gamma
(\Xbar,\calL_\mu)$ isomorphic to $V_\mu^*$ obtained by the pull back
of the homogeneous coordinates of $\mP(V_\mu)$ to $\Xbar$. Since the
multiplicity of any irreducible submodule is at most one, we can speak
of the submodule $V_\mu^*$ of $\Gamma (\Xbar,\calL_\mu)$ without
ambiguity. If now $\grl \in \Omega$ is such that $\mu \leq_\grs \grl$, then
we can consider the image of $V_\mu^* \subset \grG(\Xbar,\calL_\mu)$
under the multiplication by $s^{\grl-\mu}$.  We denote this image by
$s^{\grl -\mu}V_\mu ^*$.  We have the following Theorem:
\begin{teo}[Theorem 5.10 in \cite{DP1}]\label{prp:sezioni} If $\grl \in \Omega_{sc}$ then
  $$
  \Gamma (\Xbar,\calL_\grl) = \bigoplus_{\mu \in \Omega_{sc}^+ \st \mu
    \leq_\grs \grl} s^{\grl -\mu}V_\mu ^*. $$
\end{teo}

\subsection{Standard monomial theories}
\label{ssec:SMT}

We recall the definition of standard monomial theory.

Let $A$ be a commutative $\mk$-algebra. Let $\mA$ be a finite subset of $A$ and
 ``$<$" a partial order on $\mA$. If $a_1 \leq a_2 \leq \cdots \leq a_n$, then
we say that the monomial $a_1 \cdot a_2 \cdots a_n$ is a \emph{standard monomial}.
We denote by $\mS\mM(\mA)$ the set of all standard monomials. We
say that $(\mA,<)$ is a standard monomial theory (for short SMT) for $A$ if
$\mS\mM(\mA)$ is a $\mk$-basis of $A$.

The construction of a standard
monomial theory comes often together with the description of the
straightening relations,  i.e. a set of relations in the elements
of $\mA$ which provide an inductive procedure to rewrite a non-standard
monomial as a linear combination of standard monomials.

Let $(\mA, <)$ be a SMT for the ring $A$. In particular, $\mA$ generates $A$ and we denote by
$Rel_A$ the kernel of the natural morphism form the symmetric algebra $\sfS(\mA)$ to $A$. Let
$\mM(\mA)$ be the set of all monomials in the generators $\mA$ and let $<_t$ be a monomial order
which refines the order $<$ on $\mA$. (We recall that a monomial
order is a total order on the set
of monomials such that (i) if $m,m',m''$ are monomials and $m'<_t m''$ then $mm'<_t mm''$ and (ii)
$1<_t m$ for all monomials $m\neq 1$ (see \cite{Eisenbud}, section 15.2).) For any $a, a' \in \mA$
which are not comparable assume now that there exists $R_{a,a'} \in Rel_A$ such that
$$R_{a,a'}=a\,a' - P_{a,a'}$$
and $P_{a,a'}$ is a sum of monomials which are strictly smaller to
$a\,a'$ with respect to the order $<_t$. A set of relations
satisfying these properties is called a set of \emph{straightening
relations}. In this case we have the following simple lemma.

\begin{lem}\label{lem:SMT}
Let $(\mA,<)$ be a SMT for the ring $A$ and let $\calR = \{R_{a,a'}\st a,a'
\in \mA$ are not comparable$\}$ be a set of straightening
relations. Then $\calR$ generates $Rel_A$.
\end{lem}
\begin{proof}
Let $I$ be the ideal generated by $\calR$. We have a natural surjective morphism
$\grf : B=\sfS(\mA)/I \lra \sfS(\mA)/Rel_A =A$ induced by $I\subset Rel_A$. Now we prove that the
set of standard monomials generates the ring $B$ as a vector
space. This implies that $\grf$ is an isomorphism and hence $I=Rel_A$.

Let $m$ be any monomial and assume that it is not standard.  
The monomial can be written in the form $a\,a'\,m'$ where $a,a' \in \mA$ are not
comparable and $m'$ is a smaller monomial. So $m \coinc P_{a,a'}m' \;(mod \;I)$
and each monomial in $m'P_{a,a'}$ is strictly smaller with respect to $<_t$
than $m$ and we can conclude by induction.
\end{proof}

\subsection{Standard monomial theory for flag and Schubert varieties}
\label{ssec:SMTflag}
Let $A$ be the coordinate ring of the cone over a generalized
flag variety $\calF$ of a symmetrizable Kac-Moody group $\calG$.
For this type of algebras a standard monomial theory has been
constructed in \cite{L:SMT}. We are going to recall the main properties
of this SMT.

Fix a maximal torus $\calT$ and a Borel subgroup $\calB$ in $\calG$ such that $\calT \subset \calB$.
Let $\calL$ be a line bundle generated by global section over $\calF$ and consider the ring
$\grG_\calL(\calF)=\bigoplus _{n\geq 0} \grG(\calF,\calL^{n})$.
A basis $\mF_\calL$ of $\grG(\calF,\calL)$ has been constructed
in \cite{L:SMT} together with an order $<$ on this set such that $(\mF_\calL,<)$ is a SMT for $\grG_\calL(\calF)$.
We denote by $\mS\mM_\calL(\calL^n)$ the set of standard monomials
of degree $n$, by $\mS\mM_\calL$ the set of all standard
monomials and by $\mM(\mF_\calL)$ the set of all monomials in the set of generators $\mF_\calL$.

For $f,f' \in \mF_\calL$ that are not comparable,
the product $f\,f'$ can be expressed as a sum
$P_{f,f'}$ of standard monomials of degree two. In \cite{L:SMT} a total order
$<_t$ has been introduced on $\mM(\mF_\calL)$ with the properties required in the
previous discussion of a general SMT, so the relations $R_{f,f'}=f\,f' -P_{f,f'}$ are
a set of straightening relations. These relations are called Pl\"ucker relations since they
generalize the usual Pl\"ucker relations for the Grassmannian.

Furthermore, this theory is adapted to Schubert varieties.
Let $\calS\subset \calF$ be a closed $\calB$ stable subvariety
and set $\grG_\calL(\calS) =\bigoplus_{n\geq 0}
\grG(\calS,\calL^n\ristretto_\calS)$. Denote by $r :\grG_\calL(\calF) \lra
\grG_\calL(\calS)$ the restriction map, let $I_\calS$ be its kernel
and define $\mF_\calL(\calS) = \{a \in \mF_\calL \st r(a)\neq 0\}$. Then the set
$\{r(a) \st a \in \mF_\calL(\calS)\}$ with the order induced by the order $<$ on $\mF_\calL$
realizes a SMT for $\grG_\calL(\calS)$ and the monomials $m \in \mS\mM_\calL$ which contain
elements not in $\mF_\calL(\calS)$ form a $\mk$ basis of $I_\calS$. Finally, the {restriction
$r(R_{f,f'})$} of relations $R_{f,f'}$ {to $\calS$} for $f,f' \in \mF_\calL(\calS)$ not comparable
form a set of straightening relations. Summarizing we have:
\begin{teo}[\cite{L:SMT}]
\begin{itemize}
\item[i)] $(\mF_\calL,<)$ is a SMT for $\grG_\calL(\calF)$, and the relations $R_{f,f'}$ for $f,f'
\in \mF_\calL$ not comparable, are a set of straightening relations. \item[ii)] $(\{r(a)\mid a\in
\mF_\calL(\calS)\},<)$  is a SMT for $\grG_\calL(\calS)$, and the relations ${r}(R_{f,f'})$ for
$f,f' \in \mF_\calL(\calS)$ not comparable, are a set of straightening relations. Moreover, the
kernel $I_\calS$ of the restriction map has as basis the set of all standard monomials which
contain elements not in $\mF_\calL(\calS)$.
\end{itemize}
\end{teo}
The elements of $\mF_\calL$ are eigenvectors for the action of
$\calT$ and we denote by $weight(f)$ the weight of $f\in \mF_\calL$ w.r.t. the
action of $\calT$. The order $<$ is compatible with the dominant order
in the following way: if $f<f'$ then $weight(f)< weight(f')$ w.r.t. the
dominant order. Moreover $\mF_\calL$ has a minimum $f_0$ which is a lowest weight vector $f_0$.

The SMT described for a Schubert variety $\calS$ immediately generalizes to the Richardson
variety $\calS_0 =\{y\in \calS\st f_0(y)=0\}$ by choosing as set of
generators $\mF_0(\calS_0)=\mF(\calS)\senza\{f_0\}$. In this paper we will only need
the SMT for these particular types of Richardson varieties, a SMT for
general Richardson varieties has been constructed by Lakshmibai and
Littelmann \cite{LL}.

\medskip

In the case of the multicone over a flag variety some changes to this
general setting is needed (see \cite{C:SMT2}). We will not need these results in this paper
but we will briefly explain these changes to recall some results for the
total ring of $\Xbar$ proved in \cite{CM1}. Let $L^+ \subset
\Pic(\calF)$ be a free monoid contained in the set of elements of $\Pic(X)$
generated by global sections and let $\calL_1,\dots,\calL_r$ be the
generators of $L^+$. Define $\grG_{L^+}(\calF)=\bigoplus_{\calL \in
  L^+} \grG(\calF,\calL)$ and $\mF_{L^+}=\mF_{\calL_1}\cup \dots \cup
\mF_{\calL_r}$. It is still possible to define an order $<$ on $\mF_{L^+}$ and a total order $<_t$
on the set of all monomials with the same properties of the total order in section \ref{ssec:SMT}
such that the standard monomials of degree two are a basis of $\bigoplus_{i\leq j}
\grG(\calF,\calL_i\otimes \calL_j)$ and such that the products $f\,f'$ with $f,f'\in \mF_{L^+}$ not
comparable are a sum $P_{f,f'}$ of standard monomials that are strictly smaller to $f,f'$ with
respect to the total order $<_t$. Moreover the straightening relations $R_{f,f'}=f\,f'-P_{f,f'}$
generate the ideal of relations in the generators $\mF_{L^+}$. However, if one defines in this case
the standard monomials as above as monomials of ordered elements, then they are not anymore
linearly independent. One has to give a more restrictive definition of a standard monomial (see
\cite{C:SMT2}). For $\calL\in L^+$ we denote by $\mS\mM_{L^+}(\calL)$ the set of all standard
monomials w.r.t. to this new definition belonging to $\grG(\calF,\calL)$, and we denote by
$\mS\mM_{L^+}$ the set of all standard monomials.

\subsection{Standard monomial theory for the total ring of $\Xbar$}
\label{ssec:SMTX}
We  describe now the connection between $\mk[G/H]$ and the ring of
total sections on $\Xbar$, and we recall some properties of this ring.

For an admissible sublattice $L\subset \Omega_{sc}$ we introduce an analogue
of the total ring studied in \cite{CM1} (in that paper the ring was
called the ring of all sections):
$$
\grG_L(\Xbar) = \bigoplus_{\grl\in L} \grG(\Xbar,\calL_\grl).
$$
Set $L^+$ the subset of elements $\grl$ such that $\calL_\grl$ restricted to the closed orbit $Y$
is generated by global sections. To construct a SMT for $\grG_L(\Xbar)$ we use the SMT for the
multicone $\grG_{L^+}(Y)$ briefly explained above. Let $\gre_1, \dots, \gre_\ell \in L$ be the
basis $\calB$ of the admissible lattice $L$ as in Definition \ref{dfn:reticoli}. For all
$i=1,\dots,\ell$ and for all $f \in \mF_{L^+}(\calL_{\gre_i}\ristretto_Y)$ fix a section $f^X \in
\grG(\Xbar,\calL_{\gre_i})$ such that $f^X\vert_Y=f$. For a monomial $m=f_1\cdots f_r$ in the
elements of $\mF_{L^+}$ we denote by $m^X=f_1^X\cdots f_r^X$ the corresponding product of the
elements $f^X$.

We define $\mF^X_L = \{s_1,\dots,s_\ell\} \cup \{f^X \st f \in \mF_{L^+}\}$. We
order this set by setting $s_1<\cdots <s_\ell <f^X$ for all $f\in
\mF_{L^+}$ and we order the elements $f^X$ as in the set $\mF_{L^+}$.

We define the standard monomials $\mS\mM^X_L$ as the set $s^\nu m^X$,
where $\nu$ is a positive sums of the roots $\tal_i$ and $m \in \mS\mM_{L^+}$.
This set is a $\mk$-basis of $\grG_L(X)$ and more precisely:

\begin{teo}[\cite{CM1}]\label{SMTonWonderful}
The set  $\{s^{\grl-\mu} m^X \st \mu \in L^+,\quad \mu \leq_\grs \grl \mand
m \in \mS\mM_{L^+}(\calL_\mu\ristretto_Y)\}$
is a $\mk$ basis of $\grG(\Xbar,\calL_\grl)$.
\end{teo}
We can also give a rough description of a set of  straightening relations
in terms of the elements of $\mF^X_L$. We define a total order on the set of
monomials: let $\mu,\nu$ be positive sums of the restricted roots
$\tal_i$ and let $m,n$ be monomials in the elements in $\mF_{L^+}$.
{We set} $s^\mu m^X <_t s^\nu n^X$ if $\mu$ is less to $\nu$ with respect
to the lexicographic order, or if $\mu=\nu$ and $m <_t n$ with
respect to the total order of the monomials in the elements in $\mF_{L^+}$.

Let now $f \in \mF_{L^+}(\calL_{\gre_i}\ristretto_Y)$ and $h\in
\mF_{L^+}(\calL_{\gre_j}\ristretto_Y)$ be such that they are not comparable as elements of
$\mF_{L^+}$. Since the standard monomials form a basis of $\grG_L(\Xbar)$, we can express the
product $f^X h^X$ as a sum of standard monomials
$$
f^X h^X = P^X_{f, h} =
\sum_{\mu \in L^+ \mand \mu\leq_\grs \,\gre_i+\gre_j} s^{\gre_i+\gre_j-\mu}(P_{f,h}^\mu)^X
$$
where $P_{f,h}^\mu \in \mS\mM_{L^+}(\calL_\mu\ristretto_Y)$.

In the symmetric algebra $\sfS(\mF^X_L)$ set
$$
R_{f,h}= f\, h - P^X_{f,h}.
$$
This is a straightening relation. In fact:
\begin{teo}[\cite{CM1}]
The set of straightening relations $R_{f,h}$ for $f,h \in \mF_{L^+}$ not comparable,
generates the ideal of relations in the generators $\mF^X_L$ of  $\grG_L(\Xbar)$.
\end{teo}
The proof in \cite{CM1} of the theorem above has been only given
for $L=\Pic(\Xbar)$, but extends to the general case without changes.

The part of highest degree of the relation $R_{f,h}$ is easy to describe.
Indeed by restricting this equation to $Y$ we see
that $f\,h-P_{f,h}^{\gre_i+\gre_j}$ is the usual straightening
relation for the multicone $L^+$ over $Y$.
In a certain sense the aim of this paper is to give a description of the polynomials
$P_{\pi,\pi'}^\mu$ for $\mu \neq \, \gre_i+\gre_j$.

\subsection{A first description of the coordinate ring of $X_q$}\label{afirstdescription}
In the construction above we choose now $L=\Omega_q$ and we describe
the relation between $\grG_{\Omega_q}$ and $\mk[X_q]$.

We consider the map $\mj: X_q \lra \Xbar$ given by the composition
$X_q \lra X_{ad} \incluso \Xbar$ and we observe that for all $\grl \in
\Omega_q$ the pull back $\mj^*(\calL_\grl)$ is the trivial line
bundle. Indeed, as a representation of $H_q$,
the fiber of $\mj^*(\calL_\grl)$ over the point $H_q \in X_q$ is the
line $\mk \, h_\grl$, so the bundle is trivial.

In particular, if $\grl \in \Omega_q$ and we choose an isomorphism
$\grf_\grl: \mj^*(\calL_\grl)\lra \calO$, then we get an inclusion
$\grf_\grl: \grG(\Xbar, \calL_\grl) \incluso \mk[X_q]$. If $\Omega_q$
is admissible and $\calB$ is its {admissible basis,} then
we can choose isomorphisms $\grf_\gre$ for $\gre \in \calB$
and define $\grf_\grl = \bigotimes _{\gre\in\calB} \grf_\gre^{\otimes
  a_\gre} : \calO \lra \mj^*(\calL_\grl)=\bigotimes _{\gre\in\calB}
\mj^*(\calL_\gre)^{\otimes a_\gre}$  for $\grl = \sum_{\gre \in \calB}
a_\gre \gre$. With this choice of isomorphisms we get for all
$\grl , \mu \in \Omega_q$ the following commutative diagram:
$$
\xymatrix{ \Gamma (\Xbar,\calL_\grl)\otimes \Gamma(\Xbar,\calL_\mu)
  \ar[rr]^{\; \; multipl.} \ar[d]^{\grf_\grl\otimes \grf_\mu} & &
  \Gamma (\Xbar,\calL_{\grl+\mu}) \ar[d]^{\grf_{\grl+\mu}}  \\
  \mk[X_q]\otimes \mk[X_q] \ar[rr]^{multipl.} & & \mk[X_q]. }
$$
Hence we can define a morphism of rings $\mj ^* \defi \bigoplus_{\grl
  \in \Omega_q}
\grf_\grl : \grG_{\Omega_q} \lra \mk[X_q]$.

Observe also that $\mj^*_{\gra_i}(s_i)$ is a nonzero $G$ invariant
function on $X_q$, so we can normalize this function so that
$\mj^*_{\gra_i}(s_i) = 1$. The relation between the ring
$\grG_{\Omega_q}$ and the coordinate ring of $X_q$ is given by the
following proposition whose proof is easy.
\begin{prp}\label{liftsmt} The map $\mj^*$ gives an isomorphism
  $$\frac{\grG_{\Omega_q}}{(s_i-1 \st i = 1,\dots,\ell)} \isocan
  \mk[X_q] $$
\end{prp}

In particular we have the following corollary.

\begin{cor} \label{cor:quadratic}
  If $\Omega_q$ is quadratic{,} then the ring $\mk[X_q]$ has quadratic
  relations in the generators $\bigcup_{\gre\in \calB} V_\gre^*$.
\end{cor}

We believe that also the opposite is true:
\begin{con}
  Suppose that $\Omega_q$ is admissible. If
  $\Omega_q$ is not quadratic, then also the relations are not quadratic.
\end{con}

\subsection{Surjectivity of multiplication and applications}
We now discuss some consequences of the surjectivity of
multiplication of sections of line bundles generated by global sections.

If $\grl, \mu \in \Omega_{sc}^+$, then the line bundles $\calL_\grl, \calL_\mu$
are generated by global sections. By \cite{CM2}, the multiplication
map $m_{\grl,\mu}:\grG(\Xbar, \calL_\grl) \otimes \grG(\Xbar,
\calL_\mu) \lra \grG(\Xbar, \calL_{\grl+\mu})$ is surjective.  We
consider now the restriction $n_{\grl,\mu}$ of the multiplication map
to the submodules
$$
V_\grl^*\otimes V_\mu^* \subset \grG(\Xbar,
\calL_\grl) \otimes \grG(\Xbar, \calL_\mu)
$$
and we define $N(\grl,\mu)= \{\nu \in \grL^+ \st \nu \leq_\grs \grl+\mu \mand s^{\grl+\mu-\nu}
V_\nu^* \subset \Im n_{\grl,\mu}\}$.

We provide now a different construction of the set $N(\grl,\mu)$.
For $\grl, \mu \in \Omega_{sc}^+$ consider the element $h_\grl \otimes h_\mu
\in V_\grl\otimes V_\mu$. Let $W_\nu^{\grl,\mu}$ be the isotypic component of
type $V_\nu$ of $V_\grl\otimes V_\mu$. Denote by $\pi_\nu^{\grl,\mu}$ the
$G$ equivariant projection of $V_\grl\otimes V_\mu$ onto its isotypic
component of type $V_\nu$:
$$
\pi_\nu^{\grl,\mu} : V_\grl\otimes V_\mu \lra W_\nu^{\grl,\mu}.
$$
We define $N'(\grl,\mu) \defi \{\nu \in \grL^+ \st \pi^{\grl,\mu}_\nu(h_\grl \otimes h_\mu)
\neq 0\}$.

\begin{lem} \label{lem:decomphxh}
 With the same notation as above:
 for all $\grl, \mu \in \Omega^+$, we have $N(\grl,\mu) = N'(\grl,\mu)$.
\end{lem}
\begin{proof}
  Consider the Segre embedding $S: \mP(V_\grl)\times \mP(V_\mu) \lra
  \mP(V_\grl \otimes V_\mu)$ and define the morphism $\grD _{\Xbar} : \Xbar \lra
  \mP(V_\grl \otimes V_\mu)$ by $\grD_{\Xbar} (x) = S(\psi_\grl(x) ,
  \psi_\mu(x))$. The image of $\grD_{\Xbar}$ is the closure of the $G$ orbit
  of the vector $h_\grl \otimes h_\mu$ and $\grD_{\Xbar} ^*: V_\grl^*\otimes
  V_\mu^* \lra \grG(\Xbar,\calL_{\grl+\mu})$ is the multiplication map
  $n_{\grl,\mu}$. So $\Im n_{\grl,\mu} \supset s^{\grl+\mu-\nu}V_\nu^*$ if and only if
  $\langle G\cdot (h_\grl\otimes h_\mu) ; (W^{\grl,\mu}_\nu)^*
  \rangle \not \coinc 0$, where $(W^{\grl,\mu}_\nu)^*$ is the annihilator
  of a $G$ stable complement of $W^{\grl,\mu}_\nu$.
  Hence  $\Im n_{\grl,\mu} \supset V_\nu^*$ if and only if
  $\pi^{\grl,\mu}_\nu(G\cdot(h_\grl\otimes h_\mu)) \not \coinc 0$ and this happen if and only if
  $\pi^{\grl,\mu}_\nu(h_\grl\otimes h_\mu) \neq 0$.
\end{proof}

 The following corollary will be needed in section 4.

\begin{cor}\label{cor:hallan}
  Suppose $\tPhi$ is a simple root system of type $\sfA_\ell$, $\sfB_\ell$, $\sfC_\ell$
  or $\sfB\sfC_\ell$ and let $\gre_1,\dots,\gre_\ell$ be the quadratic basis
as in Convention \ref{conv:basee}. Then for $i=2,\dots,\ell$ we have
  $$\pi^{\gre_1,\gre_{i-1}}_{\gre_i}(h_{\gre_1}\otimes h_{\gre_{i-1}})\neq 0$$
\end{cor}
\begin{proof}
 The corollary follows by Lemma \ref{lem:ord}, Lemma \ref{lem:decomphxh},
 the description of the sections of a line bundle in Proposition \ref{prp:sezioni}
 and the surjectivity of the multiplication
 map $m_{\grl,\mu}$ proved in \cite{CM2}.
\end{proof}

\section{Construction and properties of the group  $\NuG$ and it's Lie
algebra}\label{sez:extended}

In this section we describe the Lie algebra $\Nug$ and some of
its properties. $\Nug$ is a Kac-Moody algebra endowed with a grading and an
involution. The involution contains the Lie algebra $\gog$ of $G$ as a Levi
factor in the part of degree $0$ and a spherical representation in degree
$1$. This construction depends on the choice of a spherical weight $\gre$
that we consider to be fixed.

We assume from now on the involution $\grs$ to be simple (i.e.,
$\gog$ is an irreducible $G\rtimes\{\id,\grs\}$-module) or equivalently,
$\tPhi$ is an irreducible root system.  We keep the
notation introduced in the previous
sections. In particular, the enumeration of the basis $\tal_1, \dots, \tal_\ell$
of the irreducible root system $\tPhi$ is as in \cite{Bourbaki}. Let
$\tom_1, \dots, \tom_\ell \in \Omega_{sc}$ be
the fundamental weights corresponding to this basis.

\subsection{The \extended Lie algebra}

To define the \extended Lie algebra we define its Dynkin diagram.  The new Dynkin diagram is
constructed by adding a node that we index with $0$ to the old Dynkin diagram. We join the new node
$0$ with the node corresponding to the simple root $\gra \in \grD$ with $\langle \gre, \gra\cech
\rangle$ lines and we put an arrow towards the node corresponding to $\gra$ if this number is
bigger or equal to $2$. In general this is not a Dynkin diagram of finite type. The new matrix
coefficients of the extended Cartan matrix are given by the following rules: if $\gra \in \grD$,
then we have
$$
\langle \gra_0, \gra \cech \rangle \defi - \langle \gre,
\gra \cech \rangle \quad \mand \quad \langle \gra, \gra_0\cech
\rangle \defi
\begin{cases}
  0  & \mif \langle \gre, \gra \cech \rangle = 0; \\
  -1 & \mif \langle \gre, \gra \cech \rangle \neq 0.
\end{cases}
$$

Let $\NugrD \defi \grD \cup \{\gra_0\}$ and choose a realization
$(\Nut, \NugrD, \NugrD\cech)$ of this Cartan matrix.  We define $\Nug$
as the Lie algebra constructed using this realization and we denote
by $\Nut$ its standard maximal toral subalgebra.  Denote by
$\NuPhi$ the set of roots of $\Nug$ with respect to $\Nut$
and $\NuPhi^+$ (resp. $\NuPhi^-$) the positive (resp. negative)
roots with respect to the basis
$\NugrD$.  For all $\gra \in \NugrD$ let
$\prove_\gra$ and $\provf_\gra$
be the Chevalley generators of
$\Nug$ and set $\gra \cech = [\prove_{\gra},\provf_{\gra}]$.
We can naturally identify the Lie algebra $\got$ of the maximal torus
$T$ with the subspace of $\Nut$ spanned by the $\gra\cech$ for $\gra \in \Delta$.
Moreover, we identify the Lie algebra $\gog$ of $G$ with the semisimple part
of the Levi subalgebra of $\Nug$ associated to the simple roots $\not=\al_0$.

We have also an
inclusion of $\got^* \subset \Nut^*$ induced by $\grD
\subset\NugrD$. Note that the restriction of the pairing
between $\Nut$ and $\Nut^*$ to the subspaces $\got$ and $\got^*$ induces the
usual pairing between $\got$ and $\got^*$. In
particular, if $\got^*_\perp$ is the annihilator
of $\got^*$ in $\Nut$,
then we have natural decompositions $\Nut = \got \oplus \got^*_\perp$ and
$\Nut^* = \got^* \oplus \got_\perp$. Here $\got_\perp$ denotes the annihilator
of $\got$ in $\Nut^*$.

We denote by $\Nug'$ the derived subalgebra of $\Nug$ and let
$\Nut' = \Nut \cap \Nug'$ be the subspace of $\Nut'$ spanned by the
elements in $\NugrD\cech$.  Choose an element $C$ generating
the intersection $\Nut' \cap \got_\perp^*$ and an element $D \in \got_\perp^*$ such that
$\langle D, \gra_0 \rangle = 1$. We normalize $C$ in such
a way that $\gra_0\cech \in C + \got$. Observe that
$C,D$ generate $\Nut_\perp$ and that they are linearly
independent if and only if the new Dynkin diagram is of affine type.

We grade $\Nug$ according to the action of $D$:
$$
\Nug_i \defi \{ x \in \Nug \st [D,x] = i x \}.
$$

We define now an involution $\Nugrs$ of $\Nug$ in the following
way:
$$
\Nugrs(x) = \grs (x) \mif x \in \gog; \quad \Nugrs(C) = -C; \quad
\Nugrs(D) = -D; \quad \Nugrs(\prove_0) = \provf_0 \,\mand \Nugrs(\provf_0) = \prove_0.
$$
We denote by $\Nugrs$ also the induced involution of $\Nut^*$ and
we observe that, since $\Nugrs(\got_\perp^*) =\got_\perp^* $,
we have $\Nugrs\ristretto
_{\got^*}=\grs$.  To verify that $\Nugrs$ is well defined note
that by definition for all $\gra \in \Delta$
we have $\langle \grs(\gra) , \gra_0\cech \rangle = -\langle \gra, \gra_0\cech
\rangle $ and $\langle \gra_0, \grs(\gra\cech)\rangle = -\langle
\gra_0, \gra\cech \rangle $, and hence $\Nugrs(\gra_0\cech) = -
\gra_0\cech$ and $\Nugrs(\gra_0) = - \gra_0$.

\subsection{Some remarks and conventions {concerning}
the weights of $\gog$ and $\Nug$}\label{oss:esprfund}
For  $\gra \in \De$ we denote by
  $\om_\gra\in \got^*$ the corresponding fundamental
  weight with respect to the basis $\De$.
  Let $\grD_0$ be the set of simple roots fixed by $\grs$
  and let $\grD_1$ be the complement of $\grD_0$ in $\grD$.
  Recall (\cite{DP1}) that $\grs$ induces an involution
  $\bar \grs$ of $\grD_1$ characterized by
  $\grs(\gra) + \bar\grs(\gra)$ is in the vector space spanned by
  $\grD_0$. Furthermore, $\osi$ is the restriction to $\Delta_1$ 
  of an automorphism of the Dynkin diagram of $\Phi$.

 The following connection between fundamental weights with respect to $\grD$ and fundamental
 weights with respect to $\tgrD$, as explained in \cite{CM1}, is a direct consequence of the Helgason
  criterion. For a weight $\tom_i$ we have three possibilities:
  $$
  \tom_i =
\begin{cases}
  \om_\gra            &\mif \tal=\tal_i \mand  \bar\grs (\gra) = \gra \mand \grs(\gra)\neq -\gra; \\
  2\om_\gra           &\mif \tal=\tal_i \mand  \bar\grs (\gra) = \gra \mand \grs(\gra) = -\gra; \\
  \om_\gra + \om_\grb &\mif \tal=\tal_i \mand \bar\grs (\gra) = \grb \neq \gra.
\end{cases}
$$

We fix some notation for the fundamental weights of $\Nug$.
Choose $\grg ,\grd \in \got_\perp$ univocally determined by
$\langle\grg,C\rangle=\langle\grd,D\rangle = 1$ if the new Dynkin diagram is not affine
and by $\langle\grg,C\rangle=\langle\grd,D\rangle = 1$ and
$\langle\grg,D\rangle=\langle\grd,C\rangle = 0$ if it is affine.
Notice that we have
\begin{equation}\label{eq:alpha0}
\gra_0\cech = C -\sum_{\gra\st \langle\tal; \gre\rangle\neq 0} \om_\gra\cech
\qquad\mand\qquad
\gra_0 = \grd - \gre
\end{equation}
where $\omega_\gra\cech \in \got$ are the fundamental weights w.r.t.\ $\grD\cech$. Notice also that
for $\gra \in \grD$ the weight $\om_\gra\in \got^*$ is not anymore the fundamental weight of $\gra$
w.r.t.\ to $\NugrD$ since we do not have $\langle \om_\gra , \gra_0\cech\rangle = 0$ in general. We
denote by $\nuom_\gra$ the fundamental weight of $\gra$ w.r.t.\ the extended root system. In the
affine case we normalize it in such a way that $\langle\nuom_\gra,D\rangle=0$ for all $\gra \in
\NugrD$. So we have
$$
\nuom_\gra = \om_\gra - \langle \om_\gra , \gra_0\cech \rangle\,\grg
\qquad\mand\qquad
\nuom_0 \defi \nuom_{\gra_0}= \grg.
$$
Beware that in the affine case, with these choices, we do not have
$\gra = \sum_{\grb\in \NugrD}\langle\gra,\grb\cech\rangle\nuom_\grb$
for all $\gra \in \NugrD$. Indeed
this formula holds for $\gra \neq \gra_0$ while for $\gra_0$ we
have $\gra_0 = \sum_{\grb \in \NugrD}\langle\gra_0,\grb\cech\rangle\nuom_\grb +\grd$.
In particular $\gra \ristretto _{\Nut'} = \sum_{\grb \in \NugrD}\langle\gra,\grb\cech\rangle
\nuom_\grb \ristretto_{\Nut'}$ still holds for every $\gra\in \NugrD$.

\subsection{The restricted root system of the \extended Lie algebra}

We want to study now some properties of the involution $\Nugrs$.

As in the case of the root system $\Phi$, if $\gra \in \NuPhi$ and $\Nugrs(\gra)
\neq \gra$, then we define $\tal \defi \gra - \Nugrs(\gra)$. In particular,
we have $\tal_0 \defi \gra_0 - \Nugrs(\gra_0)=2\,\gra_0$. For $i=1,\dots,\ell$
we consider the elements $\tal_i\cech \in \got$ defined in section 1
as elements of $\Nut\supset\got$ and
we define $\tal_0\cech = \tfrac{1}{2}\gra_0\cech$.

As in the classical case, we define $\nutom_0,\nutom_1,\dots,\nutom_\ell$
and we notice that we have
$$
\nutom_0 = 2\,\nuom_0 = 2 \,\gamma \qquad\mand\qquad
\nutom_i = \tom_i - \langle \tom_i , \gra_0\cech \rangle \grg \ \mfor i=1,\dots,\ell.
$$

In general we do not know
if the set of the $\tal$ with $\al \in \NuPhi$ is a root system (see Conjecture
\ref{con:tildeallargato} below for some comments).
But we can define always the Cartan matrix of this hypothetical
root system as the $(\ell+1)\times(\ell+1)$ matrix
$$
\NuCartan \defi \Big(\langle \tal_i ; \tal_j\cech \rangle \Big)_{i,j=0,\dots,\ell}.
$$
The next proposition implies that the Cartan matrix $\NuCartan$ is determined only by the restricted
root system and the weight $\gre$. In particular, it is very easy to compute.

\begin{prp}\label{prp:nucartan}
The Cartan matrix $\NuCartan$ is given by the coefficients
  of the Cartan matrix of $\tPhi$ and by the following numbers, where
  $i = 1,\dots,\ell$:
$$
\langle \tal_0 ; \tal_0 \cech \rangle = 2; \qquad
\langle \tal_0 ; \tal_i \cech \rangle = - 2\langle \gre ; \tal_i \cech
\rangle \quad \mand \quad
\langle \tal_i ; \tal_0 \cech \rangle =
\begin{cases}
 0  &\mif \langle \gre ; \tal_i \cech \rangle = 0; \\
 -1 &\mif \langle \gre ; \tal_i \cech \rangle \neq 0.
\end{cases}
$$
\end{prp}

\begin{proof}
Let $\gra \in \Delta$ be such that $\tal = \tal_i$.
We have
$$
\langle \tal_i;\tal_0\cech \rangle =
\tfrac{1}{2}\langle \gra -\grs(\gra);\gra_0\cech \rangle =
\langle \gra;\gra_0\cech \rangle
$$
which proves the third equality, while using \eqref{eq:alpha0} we obtain
$\langle \tal_0;\tal_i\cech \rangle =
\langle 2\,\grd - 2\,\gre;\tal_i\cech \rangle = - 2\langle \gre;\tal_i\cech \rangle$.
\end{proof}

In the same way the restricted root system controls many properties of $\gog$
related to the involution $\sigma$, the Cartan matrix $\NuCartan$ controls some of the properties
of $\Nug$.

\begin{prp}\label{prp:symmetrizable}
The Cartan matrix of $\Nug$ is symmetrizable if and only if $\NuCartan$ is symmetrizable.
Moreover in this case the standard bilinear form on $\Nug$  defined in \cite{Kac} is
$\Nugrs$ invariant.
\end{prp}
\begin{proof}
Recall that we assume that $\gog$ is
simple for the action of $G\rtimes\{\id,\grs\}$ (the proof in the general case is
similar). In this case there are two possibilities: either $\gog$ is simple
or $\gog=\goh\oplus\goh$, with $\goh$ a simple Lie algebra and
$\grs(x,y)=(y,x)$.

Assume first that $\gog$ is simple.
Let $d_\gra = \kappa(\gra,\gra)$ for $\gra \in \Delta$
and $\tilde d_\gra = \kappa (\tal,\tal)$ for $\tal \in \tDe$.
Then $\Nug$ is symmetrizable if and only if there exists $d_0$ such that
$d_0 \langle \gra, \gra_0\cech \rangle = d_\gra \langle \gra_0,
\gra\cech \rangle$
for all $\gra \in \Delta$.
Similarly $\NuCartan$ is symmetrizable if and only if there exists $\tilde d_0$
such that
$\tilde d_0 \langle \tal, \tal_0\cech \rangle = \tilde d_\gra \langle \tal _0, \tal\cech \rangle$
for all $\tal \in \tDe$.
Now notice that $\langle \tal, \tal_0\cech \rangle =
\tfrac{1}{2}\langle \gra-\grs(\gra), \gra_0\cech \rangle =
\langle \gra, \gra_0\cech \rangle$ and that
$\langle \tal _0, \tal\cech \rangle = 2 \langle \gra_0, \tal\cech
\rangle$
and since $\tal\cech \in \got^*$
which is spanned by the coroots $\gra\cech \in \Delta\cech$  we have
$2\langle \gra_0, \tal\cech \rangle = -2\langle \gre , \tal \cech
\rangle
= -\tfrac{4}{\tilde d _{\tal} } \kappa(\gre,\gra - \grs(\gra))
= \tfrac{4\, d _\gra}{\tilde d _{\tal} } \langle \gra_0, \gra\cech \rangle$.
So the two conditions are equivalent and $\tilde d_0 = 4 d_0$.

Assume now that $\gog=\goh\oplus\goh$ and $\grs(x,y)=(y,x)$, let
$\got_\goh$ be a maximal toral subalgebra of $\goh$ and $\grD_\goh$ a
choice of a simple basis for the roots of $\goh$. Then $\got =
\got_\goh\oplus \got_\goh$, $\grD = \grD_1 \cup \grD_2$ where $\grD_1=\{(\grb,0)\st \grb \in
\grD_\goh)\}$ and $\grD_2=\{(0,-\grb)\st \grb \in \grD_\goh)\}$
and $\gre = (\gre_\goh,-\gre_\goh)$. Since $\tPhi$ is simple the
condition for $\NuCartan$ to be symmetrizable is the same given in the
discussion of $\gog$ simple while the condition for
$\Nug$ to be symmetrizable
becomes equivalent to the existence of $d_0$ and two non zero
scalars $\grl_1$ and $\grl_2$ such that
$d_0 \langle \gra, \gra_0\cech \rangle = \grl_1 d_\gra \langle \gra_0
\gra\cech \rangle$
for all $\gra \in \Delta_1$ and
$d_0 \langle \gra, \gra_0\cech \rangle = \grl_2 d_\gra \langle \gra_0
\gra\cech \rangle$
for all $\gra \in \Delta_2$. Now if $\gre_1 =0$ the statement is
trivial while if $\langle \gre_1,\grb\cech \rangle \neq 0$ and $\gra_1
=(\grb,0)$ and $\gra_2=(0,-\grb)$ we deduce that we must have
$\grl_1=\grl_2$. The proofs can now be completed as above.

From this description it is also clear that the standard symmetric
bilinear form is $\grs$ invariant.
\end{proof}

If $\Nug$ is symmetrizable, then we denote by $\Nukappa$ the standard
symmetric bilinear form on $\Nug$ as defined in \cite{Kac}.

\begin{prp}\label{prp:finitoaffine}\hfill

\begin{enumerate}[\indent i)]
\item  The Lie algebra $\Nug$ is finite dimensional if and only if
  $\NuCartan$ is of finite type.
\item  The Lie algebra $\Nug$ is of affine type if and only if
  $\NuCartan$ is of affine type.
\end{enumerate}
\end{prp}

\begin{proof}
  In both cases we can assume that $\Nug$ is symmetrizable.

  Consider the bilinear form $\kappa'$ obtained by the restriction of
  $\Nukappa$ to the real span $\NuE$ of $\NugrD$ and the bilinear form
  $\tilde\kappa'$ obtained by the restriction of $\Nukappa$ to the real
  span $\NutE$ of $\NutgrD$.

  Then we have that $\Nug$ is of finite type if and only if $\kappa'$
  is positive definite and $\Nug$ is of affine type if and only if
  $\kappa'$ is degenerate and positive semidefinite. In the same way
  the Cartan matrix $\tilde{\mathbf A}$ is of finite (resp. affine)
  type if and only if $\tilde \kappa'$ is positive definite (resp.
  semidefinite).

  Let $E$ be the real span of $\grD$ and let $\tilde E$
  be the real span of $\tDe$, then $\NuE = E \oplus \mR C$ and $\NutE
  = \tilde E \oplus \mR C$. The restriction of the bilinear form $\kappa$
 to $E$ respectively $\tilde E$ is positive definite, and recall
  that $C$ is orthogonal to $E$.

  Hence $\kappa'$ is positive definite (resp. semidefinite) if and
  only if $\Nukappa(C,C) > 0$ (resp. $\Nukappa(C,C) = 0$). The same
  condition holds for $\tilde\kappa'$ so $\Nug$ is of finite (resp.
  affine) type if and only if the Cartan matrix $\tilde{\mathbf A}$ is
  of finite (resp.  affine) type.
\end{proof}

\begin{oss}\label{oss:levi}
If $\Nug$ is affine, then $\gra_0$ does not always correspond to the
``affine'' root in the new Dynkin diagram. In particular, the grading
$\bigoplus \Nug_i$ is not always the ``loop graduation''. It is clear that $\Nug$ is the
(non twisted) affinization of $\gog$
if and only if $\gog$ is a spherical representation and the {highest} root $\theta$
is equal to $\gre$.

If $\gog$ is spherical and $z \in \gog$ is a spherical vector,
then it is easy to prove that $H=Z_G(z)$
and hence it is a Levi subgroup. On the other hand, if $H$ is a Levi
subgroup, then $\gog$ is obviously spherical.
So $\sigma(\theta)=-\theta$ and $\tilde \theta = \theta - \sigma(\theta)$,
which is clearly the highest root of $\tPhi$, is equal to $2\theta$.
So $\tilde \theta$ divided by $2$ must be in the weight lattice of $\tPhi$.
This happens if and only if $\tilde\theta = 2 \tom_1$ and the restricted root
system is of type $\sfC$ or $\sfB\sfC$ or $\sfA_1$ (but it is
not always true that if the restricted root system is of type $\sfC$
or $\sfB\sfC$, then $\gog$ is a spherical representation).

In particular, $\theta =\gre$ if and only if $H$ is a Levi factor and
$\gre=\tom_1$
(in the reduced case of rank two the numbering is given by the fact
that here we consider the root system to be of type $\sfC_2$ and not $\sfB_2$).
\end{oss}

\begin{oss}\label{oss:nucartan}
For us the cases $\gre=\tom_1$ or $\gre=2\tom_1$ and $\tPhi$ of type $\sfB_1$ (see Convention
\ref{conv:basee}) will be of particular interest. We make explicit the results of Proposition
\ref{prp:nucartan} in these cases.

Notice that with our convention, for $\ell \geq 2$ there is no difference between the Cartan matrix
of a root system of type $\sfB_\ell$ and of a root system of type $\sfB\sfC_\ell$.
By the special choice of $\gre$ there is no difference between the Cartan matrix $\NuCartan$
obtained starting from $\tPhi$ of type $\sfB_\ell$
and that obtained starting from $\tPhi$ of type $\sfB\sfC_\ell$ (see
Conventions \ref{conv:fundBC} and \ref{conv:basee}), and in both cases
we obtain the Cartan matrix of the affine Dynkin diagram of type $\sfA_{2\ell}^{(2)}$.

In all other cases $\NuCartan$ is the Cartan matrix associated to the Dynkin diagram obtained by
adding a ``longer'' node and a double arrow from this node to the node associated to $\tal_1$, so it
is very easy to compute.

In particular, note that $\NuCartan$ is a Cartan matrix of finite type
if and only if $\tPhi$ is of type $\sfA$ and it is of affine type if and only if
$\tPhi$ is of type $\sfB, \sfB\sfC, \sfC$ or $\sfD$. Here is the list of what
we obtain in these cases.

\smallskip

\begin{center}
\begin{tabular}{|c|c|c|c|c|c|}
\hline type of $\tPhi$      &
$\sfA_\ell$, $\ell \geq 1$ &
$\sfB_\ell$, $\ell \geq 1$ &
$\sfC_\ell$, $\ell \geq 2$  &
$\sfD_\ell$, $\ell \geq 4$  &
$\sfB\sfC_\ell$, $\ell \geq 1$
\\
\hline
type of $\NuCartan$ &
$\sfC_{\ell+1}$  &
$\sfA_{2\ell}^{(2)}$ &
$\sfC_\ell^{(1)}$ &
$\sfA_{2\ell-1}^{(2)}$ &
$\sfA_{2\ell}^{(2)}$
\\
with $\gre=\tom_1$ or $2\tom_1$ for $\sfB_1$ &&&&&
\\
\hline
\end{tabular}
\end{center}

\end{oss}

\smallskip

\subsection{First  properties of the extended Lie algebra}\label{ssec:triple}
From now on we assume $\Nug$ to be symmetrizable and we denote by
$\Nukappa$ the standard symmetric bilinear form on $\Nug$ as defined
in \cite{Kac}. In general it is not true that the
restriction of $\Nukappa$ to $\gog$ is a multiple of the Killing form
$\kappa$. So if we identify $\Nut^*$ with $\Nut$ using
$\Nukappa$ and we define $x\cech = \tfrac{2\,x}{\Nukappa(x,x)}$
for an element  $x$ of non zero length, then for $x\in \got$ this definition does not need to
agree with definition of $x\cech$ given in section 1.
However, if one has an ideal of $\gog$ which is simple for the action of
$G\rtimes\{\id,\grs\}$, then by the uniqueness of the $\grs$-invariant
bilinear form, the restriction of $\Nukappa$ to such an ideal must be a multiple
of the Killing form. So the two possible definitions of $x\cech$ coincide for elements
which belongs to such an ideal. In particular, they coincide for all elements in
$\Phi$ and $\tPhi$. For this reason we keep
the same symbol $\gra\cech$.

We list now some properties of the Lie Algebra $\Nug$.

\begin{prp}\label{prp:primeNug} \hfill

\begin{enumerate}[\indent i)]
\item $\Nug = \bigoplus _{i \in \mZ} \Nug_i$ and $\Nugrs (\Nug_i) =
  \Nug_{-i}$ for all $i \in \mZ$;
\item $\Nug_0 = \gog \oplus \got^*_\perp$ so $\gog$ is the semisimple part of a
Levi factor of
  $\Nug$ and any $\Nug_i$ is a $\gog$ module;
\item $\Nug_{-1} \isocan V_{\gre}$ as a $\gog$ module;
\item the subalgebra $\Nug_- \defi \bigoplus_{i>0} \gog_{-i}$ is
  generated by $\Nug_{-1}$;
\item the subalgebra $\Nug_+ \defi \bigoplus_{i>0} \gog_{i}$ is
  generated by $\Nug_{1}$;
\item for all $i\in \mZ$ we have $\dim \Nug_{i} < \infinito$;
\item for all $i\in \mZ$ we have $\Nug_{-i} \isocan \Nug_i^*$ as a
  $\gog$ module.
\end{enumerate}
\end{prp}

\begin{proof}
  $i)$ and $ii)$ are a direct consequence of the definition.

  To prove $iii)$ we show that $\provf_{\gra_0}$ generates $\Nug_{-1}$ as a $\gog$
  module and that it
  is a {highest} weight vector for the action of $\gog$ of weight
  $\gre$. The second claim is trivial since $[\prove_\gra,\provf_{\gra_0}]= 0$ for $\gra \in \Delta$
  and $\langle -\gra_0, \gra\cech\rangle = \langle \gre,
  \gra\cech\rangle$ by definition.
  Consider the subalgebra $\Nug_- = \bigoplus _{ i < 0}
  \Nug_i$. By \cite{Kac} it is contained in the
  subalgebra generated by the elements $\provf_i$. Hence $\Nug_{-1}$ it is
  generated by the elements of the form $[\provf_{\gra_{i_1}}\dots
  [\provf_{\gra_{i_m}}[\provf_{\gra_0}[\provf_{\gra_{j_1}}\dots[\provf_{\gra_{j_{n-1}}},\provf_{\gra_{j_n}}]\dots]$
  with $\gra_{i_1},\dots,\gra_{j_n}\in \Delta$. Since
  $x=[\provf_{\gra_{j_1}}\dots[\provf_{\gra_{j_{n-1}}},\provf_{\gra_{j_n}}]\dots] \in \gog$, we can rewrite
  the element above as
  $ -[\provf_{\gra_{i_1}}\dots[\provf_{\gra_{i_m}}[x,\provf_{\gra_0}]]\dots]$ which proves the claim.

  Similarly we observe that $\Nug_{-i-1}$ is the $\gog$ module spanned
  by $[\provf_{\gra_0},\Nug_{-i}]$.  Now if $x\in \gog$ and $y \in \Nug_{-i}$ we
  have $[x,[\provf_{\gra_0},y]]= [\provf_{\gra_0},[x,y]]+[[x,\provf_{\gra_0}],y] \in
  [\Nug_{-1},\Nug_{-i}]$. Hence $\Nug_{-i-1} = [\Nug_{-1}, \Nug_{-i}]$
  and this implies $iv)$. Point $v)$ is similar and $vi)$ follows by $iii)$ and $iv)$.

  Finally, to prove $vii)$ note that with respect to the non degenerate
  bilinear form $\Nukappa$ the subspace $\Nug_i$ is in duality with
  $\Nug_{-i}$.
\end{proof}

We introduce now a triple of elements in $\Nug$. By Lemma \ref{prp:primeNug}
we know $\Nug_1\isocan V_{\gre}^*$ and $\Nug_{-1}\isocan V_{\gre}$,
so we can choose spherical vectors $h_1\in \Nug_1$ and $h_{-1}\in \Nug_{-1}$ and
define $K =[h_1,h_{-1}]$.

\begin{lem}\label{lem:tripla}\hfill
\begin{enumerate}[i)]
\item $\kappa(h_1,h_{-1})\neq 0$ and $K=[h_1,h_{-1}]\neq 0$;
\item If $\Nug$ is not of affine type, then we can choose
$h_1$ and $h_{-1}$ in such a way that $h_{-1},K,h_1$ is an $sl(2)$ triple;
\item If $\Nug$ is of affine type then $ [K,h_1] = [K,h_{-1}] = 0 $.
\end{enumerate}
\end{lem}
\begin{proof}
 $H$ is reductive, there is only one line of elements fixed by
$H$ and $\kappa$ gives a $G$ equivariant isomorphism
between $\Nug_{-1}$ and $\Nug_1^*$, so we must have that
$\kappa(h_1,h_{-1}) \neq 0$.

The Lie bracket defines a surjective map
$\Nug_{-1}\otimes \Nug_1 \lra \Nug'_0 \defi \Nug_0 \cap \Nug'$
and $\Nug'_0 = \gog \oplus \mk C$ as a $\gog$ module.
The composition with the projection on the trivial factor is the only $G$ equivariant
map from $\Nug_{-1}\otimes \Nug_1$ to a trivial representation,
so it must be a non zero multiple of the map given by
$x_1\otimes x_{-1} \mapsto \kappa(x_1,x_{-1}) C$. In particular, $K \neq 0$
and, up to a nonzero scalar, we have $ K = C + x $ with $x \in \gog$.
Since $h_1$ and $h_{-1}$ are fixed by $H$, so is $x$ and hence either $x$ is
a non zero spherical vector or $x=0$. In the first case (see Remark \ref{oss:levi})
it is easy to prove that $\goh=Z_\gog(x)$. In particular, $x \in \goh$,
$[x,h_1]=0$ and $[x,h_{-1}]=0$ since $h_1$ and $h_{-1}$ are spherical.
So in both cases ($x =0$ or not) we have $[K,h_1]=[C,h_1]$ (and the same for $h_{-1}$).

Now $ii)$ and $iii)$ follows by the fact that if the diagram is affine, then
$C$ is central, and if it is not affine, then $C$ is a
nonzero scalar multiple of
$D$ and hence $K$ acts non-trivially on $\Nug_{1}$ and $\Nug_{-1}$.
\end{proof}

\subsection{The Weyl group of the \extended Lie algebra}
Notice that if $\gra \in \NuPhi^+$ and $\Nugrs(\gra) \neq \gra$, then
$\Nugrs(\gra) \in \NuPhi^-$. If $\Nug$ is finite dimensional, then
this implies that the maximal toral subalgebra $\Nut$ is maximally
split. In the infinite dimensional case we would like to consider this
property as the analogue for the toral subalgebra $\Nut$ to be maximally
split and we would like to prove for this situation the analogous
basic structural properties as in the finite dimensional case.

In \cite{Ric} the relation between the Weyl group
$W$ of the root system $\Phi$ and the Weyl group $\tW$ of the root system
$\tPhi$ is described.
Let $\gos \subset \got$ be the $(-1)$ eigenspace of the action of
$\grs$ on $\got$ and set
$W_1 =\{w\in W \st w(\gos)\subset \gos \}$ the subgroup of
$W$ preserving the span of spherical weights and
$W_2= \{w \in W \st w\ristretto_{\gos} =\id_{\gos}\}$
the subgroup of $W_1$ acting trivially on spherical weights. The restriction to $\gos$
gives an injective map $r:W_1 / W_2 \lra \Aut(\gos)$. The relation between $W$ and $\tW$
is given by the following Proposition.

\begin{prp}[Richardson \cite{Ric}, Proposition 4.7]\label{prp:WeWtilde}
$r$ defines an isomorphism between $W_1 / W_2$ and $\tW$.
\end{prp}

We generalize now this result to the extended situation. We prove first a weak form. Let $\Nugos$
be the $(-1)$-eigenspace of the action of $\Nugrs$ on $\Nut$ and for $i=0,\dots,\ell$ let $\ts_i$
the reflection of $\Nugos$ defined by the simple root $\tal_i$. Let also $\NuW \subset \Aut(\Nut)$
be the Weyl group of the root system $\NuPhi$ and $\NutW \subset \Aut(\Nugos)$ the group generated
by the reflections $\ts_i$ for $i=0,\dots,\ell$. As in the finite dimensional case define $\NuW_1
=\{w\in \NuW \st w(\Nugos)\subset \Nugos \mand w\ristretto_\Nugos \in \NutW  \}$ and $\NuW_2= \{w
\in \NuW \st w\ristretto_{\Nugos} =\id_{\Nugos}\}$. The restriction to $\Nugos$ gives an injective
map $\Nur:\NuW_1 / \NuW_2 \lra \Aut(\Nugos)$ and the analogue of Proposition \ref{prp:WeWtilde}
holds.

\begin{lem}\label{cor:NuWeWtilde}
$\Nur$ defines an isomorphism between $\NuW_1 / \NuW_2$ and $\NutW$.
\end{lem}

\begin{proof}
Note that $\{s_\gra \st \gra \in \Delta\}$
(resp. $\{\ts_i \st i=1,\dots,\ell\}$) generates a subgroup of $\NuW$ (resp. $\NutW$) isomorphic to
$W$ (resp. $\tW$) which acts trivially on $\Nut' \cap \got^*_\perp$
(resp. $\Nugos \cap \got^*_\perp$). So we have the following commutative diagram
$$
\xymatrix{
   W_1 \ar[r]^{r}   \ar@{}[d]|{\cap}  &  \tW   \ar@{}[d]|{\cap} \\
\NuW_1 \ar[r]^{\Nur}                  &  \NutW
}
$$
Hence it is clear that $\ts_1,\dots,\ts_\ell \in \Im \Nur$ by the finite case result
of Proposition \ref{prp:WeWtilde} and it remains to prove that $\ts_0\in \Im\Nur$.
But $\tal_0 = 2\gra_0$ and $\ts_0 = \Nur(s_{\gra_0})$.
\end{proof}

It is possible to describe
explicit covers of the generators of the Weyl group by
describing explicit elements $w_i \in \NuW_1$ such that $\Nur(w_i)
=\ts_i$.

For $i = 0, 1 ,\dots, \ell$ let
$\Si_i=\{\al \in \NugrD \st \tal=\tal_i\}\cup\De_0$.

\begin{prp}\label{prp:sollevamenti}\hfill

\begin{enumerate}[\indent i)]
\item Let $w_\grD$ be the longest element in $W$ (with respect to the simple roots
$\grD$), then $w_\grD(\grD_0) = -\grD_0$;
\item $w_\grD \comp \grs = \grs \comp w_\grD$;
\item Denote by $w_i$ the longest element of the Weyl group of $\Si_i$. Then
$w_i \in \NuW_1$ and $\Nur(w_i)=\ts_i$.
\end{enumerate}
\end{prp}

\begin{proof}
We prove first that $ii)$ implies $iii)$.
For $i=0$ it is trivial: $\Nugos$ is orthogonal to $\grD_0$ so
if $ w_{\grD_0}$ is the longest element of the Weyl group associated to
$\grD_0$, then
$ w_{\grD_0} \in \NuW_2$. Also notice that $\gra_0$ is not joined to $\grD_0$,
so we have $w_0 = s_{\gra_0} \comp w_{\grD_0}$
and $\Nur(w_0)=\ts_0$ follows from $\Nur(s_{\gra_0})=\ts_0$.

So we can reduce the proof to the finite dimensional case.
Let $\tom_h$ be a fundamental weight of $\tPhi$ orthogonal to $\tal_i$. Then $\tom_h$ is sum of
fundamental weights $\om_\gra$ orthogonal to any root in $\Si_i$. This shows that
$w_i(\tom_h)=\tom_h$. So it suffices to show that $w_i(\tal_i)=-\tal_i$.

Let $\got_i$ be the vector space spanned by $\Sigma_i$ and $\Phi_i$ the root system generated by
$\Sigma_i$. $\grs$ preserves $\Phi_i$, so by considering $\grs\ristretto_{\gog_i}$
we can assume that the rank of the involution $\grs$ is $1$.
In particular, $w_\grD = w_i$ in this case and it commutes with $\grs$, hence $w_i$
preserves $\gos_i \defi \gos\cap \got_i$. We have already seen that the orthogonal complement
to $\tal_i$ in $\gos_i$ is fixed by $w_i$ and hence, since $w_i$ is a real isometry,
$w_i(\tal_i)=\pm \tal_i$. Moreover, note that $\tal_i \in \mN[\Phi_i^+]$ so
$w_i(\tal_i) \in -\mN[\Phi_i^+]$ and hence $w_i(\tal_i) = -\tal_i$.

Now we prove $i)$ implies $ii)$. Notice first that to prove $ii)$
it is enough to examine the case of a simple involution.
Notice also that in the case of the flip: $\grs(x,y)=(y,x)$ the claim is trivial.
So we can assume that $\gog$ is simple. If $\gra \in W$, then
$\grs \comp s_\gra \comp \grs = s_{\grs(\gra)}$, so $\grs$ acts on $W$ by
conjugation. If $w_\grD(\grD_0) = - \grD_0$, then $w_\grD$ preserves $\Phi_0$.
Hence if we consider $w' = \grs \comp w_\grD \comp \grs$, then we have that it is an element
of the Weyl group that takes positive roots into negative roots so $w'=w_\grD$
and $w_\grD \comp \grs = \grs \comp w_\grD$.

Finally $i)$ is a special case of Lemma 15.5.8 in \cite{Springer}.
\end{proof}

Let now $\NuE=\grL\otimes_\mZ \mR + \mR\,\grg +\mR \grd \subset
\Nut^*$. Let $\NuPhi_{re}$ be the $\NuW$ orbit of $\NugrD$ and define
the subsets $A$ and $U$ of $E$ as
\begin{align*}
A & =\{x \in \NuE \st \langle x; \gra\cech\rangle \geq 0 \mforall \gra \in
\NugrD \} \\
U & =\{x\in \NuE \st \langle x; \gra\cech\rangle \geq 0 \mforall \gra \in
\NuPhi_{re} \text{ but a finite number}\}.
\end{align*}
Then $\NuW \,A = U$ and $A$ is a fundamental domain for the action of $\NuW$ on $U$. Define also
$\NutE = \grL^s \otimes_\mZ \mR + \mR \, \grg + \mR \, \grd$, let $\NutPhi_{re}$ be the $\NutW$
orbit of $\NutgrD$. Let $\tilde A$, $\tilde U$ be defined in the same way as $A$ and $U$, then
$\tilde U$ is stable by the action of $\NutW$ and $\tilde A$ is a fundamental domain for the action
of $\NutW$ on $\tilde U$.

\begin{lem}\hfill

\begin{enumerate}[\indent i)]
\item For all $x \in \NutPhi_{re}$ there exists $\gra \in
  \NuPhi_{re}$ such that $x = \tilde \gra = \gra-\Nugrs(\gra)$;
\item $A \cap \NutE =\tilde A$ and $U \cap \NutE =\tilde U$;
\end{enumerate}
\end{lem}

\begin{proof}
 $i)$ If $w\in \NuW_1$, then $w$ commutes with
 $\Nugrs$. Indeed $\Nut = \Nut_+ \oplus \Nugos$ where
 $\Nut_+$ is the subspace fixed by $\Nugrs$. By $\grs$ invariance
 of $\Nukappa$ this is an orthogonal decomposition of $\Nut$. So if
 $w\in \NuW$ preserves $\Nugos$ it also preserves $\Nut_+$ and by
 consequence commutes with $\Nugrs$.

 If $x \in \NutPhi_{re}$, then $x = w (\tbe)$ with
 $\grb\in \grD$ and, by Lemma \ref{cor:NuWeWtilde}, $w \in W_1$.
 So $x = w(\grb)-w(\Nugrs(\grb)) = w(\grb)-\Nugrs(w(\grb))= \tal$
with $\gra = w(\grb)\in \NuPhi_{re}$.

 $ii)$ The statement about $A$ is obvious since $\langle x, \tal\rangle
 = 2\langle x, \gra\rangle$ for all $x\in \Nugos$ and for all $\gra
 \in \NugrD$. Moreover, $U\cap \NutE \supset \NuW_1( A) \cap \NutE =
 \NuW_1(A\cap \NutE)=\NutW( \tilde A) =\tilde U$. Finally if $x \in U
 \cap \NutE$ by point $i)$ and $\langle x, \tal\rangle
 = 2\langle x, \gra\rangle$ we have also $x\in \tilde U$.
\end{proof}

In the next section we will need the following integral form of $\NutE$:
$\NuO= \Omega + \mZ \, \grg + \mZ \, \grd$.

\begin{cor}\label{lem:Waffine}\hfill

\begin{enumerate}[\indent i)]
\item $\NuW_1 = \{w \in \NuW \st w(\Nugos) = \Nugos \}$
\item If $\grl \in U \cap \NuO$ and $w\in \NuW$
  is such that $w(\grl) \in A$, then $w(\grl) \in \NuO$ and there
  exists $\tilde w \in \NutW$ such that $\tilde w (\grl)=w(\grl)$.
\end{enumerate}
\end{cor}

\begin{proof}
We prove $ii)$, the proof of $i)$ is similar.

Choose $\tilde w \in \NutW$ such that $\tilde w(\grl) \in
\tilde A \subset A$. By Corollary \ref{cor:NuWeWtilde}, $\tilde w$ is
the restriction to $\Nugos$ of an element of $\NuW$, so since $A$ is a
fundamental domain we have $\tilde w(\grl)=w(\grl)$.

Hence it is enough to prove that if $\grl \in \NuO$ and $\tilde w \in
\NutW$, then $\tilde w (\grl) \in \NuO$. This is clear if $\tilde w \in
\tW$ since $\tW$ preserves $\NuO$ and fixes $\grg$ and $\grd$. So it is
enough to consider the case $\tilde w = \tilde s_0 =s_0$. In this case the claim follows
from $\gra_0=\grd-\gre$ and $\langle\grl;\gra_0\cech \rangle \in \mZ$ if $\grl
\in \NuO$.
\end{proof}

If one tries to develop an analogue of the classical finite dimensional theory for this situation,
one of the first questions that one needs to address is to
clarify the relationship between $\NutPhi$  and the Cartan matrix
$\tilde{\mathbf A}$. More precisely we have the following conjecture.

\begin{con}\label{con:tildeallargato} Suppose that $\tPhi$ is not of type
  $\sfB\sfC$. Consider the realization of the Cartan matrix $\tilde{\mathbf A}$
  given by $(\Nugos,\NutPhi,\NutPhi\cech)$ and the root system $\Psi$ of its
  associated Kac-Moody algebra.  Then $\Psi = \NutPhi$.
\end{con}

When $\tPhi$ is of type $\sfB\sfC$ we could adjust the conjecture to
give it a reliable appearance but there seems to
be no general theory of nonreduced Kac-Moody root system.

Notice that the conjecture is
true in the case $\Nug$ is finite dimensional
or in the case $\Nugrs\ristretto_{\Nut} = -\id_{\Nut}$.
It is also easy to verify the conjecture in the case $\Nug$ is the
affinization of $\gog$
since we have an explicit description of the root system.

\section{The representation $Z$ and the Richardson variety $\Rich$}
\label{sez:schubert}

In this section we introduce a representation $Z$ of $\Nug$ and a Richardson variety $\Rich$ and we
prove the main technical results of the paper.

We keep the notation introduced in the previous section.
Moreover, from now on we fix a simple involution of $G$ and a subgroup of the
form $H_q$ (see section \ref{sez:1}) such that $\Omega_q$ is quadratic. We denote
$\Omega_q$ by $\Omega$ and $H_q$ by $H$. We denote by
$\gre_1,\dots,\gre_\ell$ the admissible basis of $\Omega$  as in
Convention \ref{conv:basee} and we choose $\gre = \gre_1$ in the
construction of $\Nug$ given in the previous section.
We set for convenience $\gre_0=0$ and
we define $\nugre_i = \gre_i - \langle \gre_i, \gra_0\cech\rangle \grg \in \Nut^*$
and for convenience $\nugre_0 = \nutom_0$.
In particular, the restricted root system
$\tPhi$ is of type $\sfA$, $\sfB$, $\sfC$ or $\sfBC$, and the Lie algebra
$\Nug$ is of finite or affine type (see Remark \ref{oss:nucartan}).

\subsection{The representation $Z$}
Let $Z$ be the integrable highest weight module of $\Nug$ with highest weight
$\nuom_0$ and let $z_0$ be an highest weight vector in $Z$.
We define a grading of $Z$ using the action of $D$ in the following way: let
$n_0=\langle\nuom_0,D\rangle$ and set
$$
Z_n = \{ z \in Z \st D\cdot z = (n+n_0) z \}
$$
This grading is compatible with the grading of $\Nug$ introduced
in the previous section and, by Proposition \ref{prp:primeNug}
$vi)$, each $Z_n$ is a $\gog$ module,
is finite dimensional and is zero for $n>0$.
We define the restricted dual $Z^*$ of $Z$ as $Z^* = \bigoplus_{n\geq0} (Z_{-n})^*$.
$Z^*$ is the integrable lowest weight module with lowest weight  $-\nuom_0$, and is graded
by the action of $D$ with $(Z^*)_n = (Z_{-n})^*$. We choose $z_0^*$ a
lowest weight vector
such that $\langle z_0, z_0^*\rangle =1$.

\begin{oss}
We have $n_0=0$ if $\Nug$ is of affine type and
$n_0=\tfrac{\ell+1}{2}$ if it is of finite type.
This can be easily computed by noticing that in finite type case
$\NutPhi$ is of type $\sfC_{\ell+1}$ and $D=\nuom_0\cech$.
\end{oss}

We need some information on the decomposition of $Z_n$ into $\gog$
modules.
We denote by $\leq$ the dominant order on $\Nut^*$ and we extend
the order $\leq_\grs$ to $\Nut^*$ by saying that $\mu \leq_\grs \grl$
if $\grl-\mu \in \mN[\NutgrD]$. Furthermore,
if $\grl \in \Nut^*$ is such that $\Nugrs(\grl)=-\grl$, then $\grl$ can
be written in the form $\grl=\sum_{i=0}^{\ell} a_i\, \gre_i + a\,\grg +
b\,\grd$ and we define
$$
\Nugr(\grl) \defi \sum_{i=0}^\ell i\, a_i - \langle D, \grl\rangle.
$$
Notice that we have
$$
\Nugr(\nuom_0) = n_0, \qquad
\Nugr(\tal_0)= \dots = \Nugr(\tal_{\ell - 1})= 0 \quad \mand \quad \Nugr(\tal_\ell)>0.
$$
More generally, if $\grl \in \Nut^*$, then we define
$\Nugr(\grl) = \tfrac{1}{2}\Nugr(\grl-\Nugrs(\grl))$. Recall that
$\NuO= \Omega + \mZ \, \grg + \mZ \, \grd$ and that
$\grD_0=\{\gra\in\NugrD\st \Nugrs(\gra)=\gra\}$.

\begin{prp}\label{prp:GmoduliZn}
Let $\grl \in \Nut^*$ be a weight of the $\Nug$ module $Z$. Then
$\Nugr(\grl) \leq \Nugr(\nuom_0)$ and moreover if $\grl \in \NuO$ and
$\Nugr(\grl) = \Nugr(\nuom_0)$, then $\grl \in \NutW_\ell(\nuom_0)$,
where $\NutW_\ell$ is the subgroup generated by $\ts_0,\ts_1, \dots,\ts_{\ell-1}$.
\end{prp}

\begin{proof}
The fact that $\Nugr(\grl) \leq \Nugr(\nuom_0)$ follows from $\grl
\leq \nuom_0$ and $\Nugr(\gra) = \Nugr(\tal) \geq 0$ for $\tal \in
\NugrD\senza \grD_0$ and $\Nugr(\gra)=0 $ if $\gra \in \grD_0$.

Assume now that $\grl \in \NuO$ and that $\Nugr(\grl) = \Nugr(\nuom_0)$. Let $\tw \in \NutW$ be
such that $\mu=\tw(\grl) \in \tilde A$ (see Corollary \ref{lem:Waffine}). By the description of the
weights of the integrable module $Z$ and Lemma~\ref{lem:Waffine} we have $\mu \leq \nuom_0$ and
$\mu \in \NuO$ and hence $\Nugr(\grl) = \Nugr(\mu) = \Nugr(\nuom_0)$.
Moreover, since $\Nugr(\tal_\ell)<0$ we can choose
$\tw$ in $\NutW_\ell$.

So it is enough to prove that $\nuom_0$ is the only element $\nu$ of
$\NuO \cap A$ such that $\nu \leq \nuom_0$ and
$\Nugr(\nu)=\Nugr(\nuom_0)$. Take $\nu$ with
these properties and consider $\tilde \nu = 2 \,\nu$. Then $\tilde \nu
\leq_\grs \tom_0$ and let $\tom_0 - \tilde \nu =\sum_{i=0}^\ell
b_i\,\tal_i$. Then from $\Nugr(\tilde \nu)=\Nugr(\nutom_0)$ and
$\Nugr(\tal_\ell)>0$ we deduce $b_\ell=0$. Now by
Proposition~\ref{prp:nucartan} and Remark~\ref{oss:nucartan},
the root system generated by $\tal_0,\dots,\tal_{\ell-1}$ is of
type $\sfC_\ell$ (numbered from $\ell-1$ to $0$). Further,
$\tilde \nu$ is a weight with respect to this root system, and
$\tilde \nu$ is less or equal to $\nutom_0$ with respect to the
dominant order of this root system (since $b_\ell=0$). A simple computation
for a root system of type $\sfC_\ell$ then shows that the elements
with these properties are given by the following list:
$$
\nutom_0 > \nugre_2 +\grd_0 > \nugre_4 + 2\grd_0 > \cdots
$$
where $\grd_0=0$ if $\Nug$ is of finite type and is equal to $\grd$ if
it is of affine type.  In particular, $\tilde \nu$ must be one of these weights and $\nu =
\tfrac{1}{2}\tilde \nu$ belongs to $\NuO$ only if $\tilde \nu = \nutom_0$
and $\nu=\nuom_0$.
\end{proof}

As we have already noticed in the proof of the proposition, the root system
generated by $\tal_0, \dots,\tal_{\ell-1}$ is always of type
$\sfC_\ell$ so we can easily compute the orbit $\NutW_\ell \nuom_0$.
In particular, we are interested in the weights in this orbit that are
dominant with respect to $\grD$ (or equivalently $\tDe$). We describe now these
weights. Recall that a root system of type $\sfC_\ell$ can be
realized in $\mR^\ell$, with standard basis $e_1,\dots,e_\ell$,
as the set $\{\pm e_i \pm e_j \st
i,j=1,\dots,\ell\}\senza\{0\}$ and $\gra_1^{\sfC}=e_1 - e_2,\dots,
\gra_{\ell-1}^{\sfC}=e_{\ell-1} - e_\ell, \gra_\ell^{\sfC}=2\,e_\ell$
is a simple basis. Then an element $x = \sum x_i
\, e_i$ is an integral weight if the coefficients $x_i$ are integers,
is a dominants weight w.r.t.
$\gra_1^{\sfC},\dots,\gra_{\ell-1}^{\sfC}$ if and only if $x_1\geq
\dots \geq x_\ell$ and the fundamental weight $\om_i^{\sfC}$ is the
element $\sum_{j\leq i}e_i$. The Weyl group is isomorphic to
$\sfS_\ell \ltimes (\mZ/2)^\ell$ where $\sfS_\ell$ acts by
permutations and $(\mZ/2)^\ell$ by changing the sign of the elements $e_i$.
Hence the elements of the Weyl group orbit of $\om_\ell^{\sfC}$ that
are dominant with respect to the first $\ell-1$ roots
are the elements $e_1+\dots+e_i-e_{i+1} \cdots - e_\ell$ for $i=0,\dots,\ell$.
In particular, there are $\ell+1$ of these elements.

\subsection{Some special Weyl group elements}\label{specialtau}
We describe the elements $\sutau_0,\dots, \sutau_\ell$ of $\NutW_\ell$ whose action on
$\nuom_0$ gives all the weights in $\NutW_\ell(\nuom_0)$ that are dominant
with respect to $\grD$. Let
$$
\sutau_0 = \id \quad\text{ and for } m=0,\dots,\ell-1 \qquad \sutau_{m+1} = \ts_0 \ts_1 \ts_{2}\cdots \ts_m \sutau_m,
$$
and define $\tau_m=w_\grD \sutau_i$. Set also $\sutau =\sutau_\ell$ and $\tau=\tau_\ell$.
Then

\begin{lem}\label{lem:primeX}\hfill

\begin{enumerate}[\indent $i$)]
\item For $i=0,\dots,\ell$ we have
$$
\sutau_i(\nuom_0) =
\begin{cases}
\nugre_i - \nuom_0           &\mif \Nug \text{ is of finite type};\\
\nugre_i - \nuom_0 - i \grd  &\mif \Nug \text{ is of affine type};
\end{cases}
$$
in particular {\color{mycolorP}$\sutau_i(\nuom_0)\ristretto _{\Nut'} = (\nugre_i
-\nuom_0)\ristretto _{\Nut'}$}; \item For $i=0,\dots,\ell$ we have $\langle\sutau_i(\nuom_0) , D
\rangle = n_0-i$; \item $\{\sutau_m(\nuom_0 ))\st m=0,\dots,\ell\} = \{\grl \in
  \NutW_\ell(\nuom_0) \st \grl \text{ dominant w.r.t. } \grD\}$.
\end{enumerate}
\end{lem}
\begin{proof}
To prove $i$) note that this is a computation which involves only objects related to
the Weyl group $\NutW$. So it is enough to notice that $\nugre_0 = \nutom_0 =2\,\nuom_0$
and that by Remark~\ref{oss:nucartan} (and also see Convention~\ref{conv:fundBC} for
the $\sfB\sfC_1$ case)
\begin{align*}
\tal_0 &= 2\, \nugre_0 - 2 \nugre_1 +  2\,\grd_0; \\
\tal_i &= 2\, \nugre_i - \nugre_{i-1} - \nugre_{i+1} \;\mfor 1\leq i \leq \ell -1;
\end{align*}
where $\grd_0 = 0$ if $\Nug$ is of finite type and $\grd_0=\grd$ if $\Nug$ is of affine type.

$ii)$ and $iii)$ now follows from $i)$ and the fact that by the
discussion above the set on the right side
in $iii)$ has $\ell+1$ elements.
\end{proof}

We now restate the results of this discussion in the form we will use
it in section \ref{sez:equazioni}. For $\grl \in \Omega$ and $\grl = \sum_i a_i \,
\gre_i$ we define $gr(\grl)=\sum_i i\,a_i$.

\begin{cor}\label{cor:GmoduliZn}
Let $\grl \in \Omega$ be such the $V^*_\grl$ appears as a $G$-module in
$(Z^*)_n$. Then $gr(\grl) \leq n$ and if $gr(\grl)=n$ then $n\leq
\ell$ and $\grl=\gre_n$. Moreover the multiplicity of $V^*_{\gre_n}$
in $(Z^*)_n$ is one.
\end{cor}
\begin{proof}
The first part of the Corollary is just a restatement of Proposition
\ref{prp:GmoduliZn}. The last statement follows from the fact that
each weight in the orbit $\NuW (\nuom_0)$ appears with multiplicity one.
\end{proof}

\subsection{The Schubert variety and the Richardson variety}
We denote by $\NuG$ the (minimal) Kac Moody group (see \cite{Kumar} pg. 228)
associated to the Lie algebra $\Nug$ and
by $\NuP$ the stabilizer of the line $\mk z_0$, so $\Gr=\NuG/\NuP$ its the associated
Grassmannian. On this Grassmannian we consider
the line bundle $\calL$ whose space of sections is the $\Nug$ module $Z^*$.

Let $\NuB$ be the Borel subgroup of $\NuG$ corresponding to the positive roots. Recall
that the $\NuB$ orbits in $\Gr$ are parametrized by $\NuW/\NuW_\NuP$
and that $\NuW_\NuP$, the Weyl group associated to $\NuP$, is equal to $W$.
For $w \in \NuW$ we denote by $[w]$ its class in $\NuW/ W$, we recall that
the set $\NuW/ W$ is partially ordered by  the inclusion relations
corresponding to the orbit closures of the $\NuB$-orbits.
In particular, the closure of $\NuB w \NuP/\NuP$
is given by all the orbits $\NuB w' \NuP/\NuP$ with $[w']\leq [w]$.

Consider the Schubert variety $\Schub_{\tau_m}\defi\overline{\NuB\tau_m\NuP/\NuP}$
and the module of sections $\grG(\Schub_{\tau_m}) = \grG(\Schub_{\tau_m}, \calL)$. This module is
a graded quotient of $Z^*$ and we denote by $\grG_n(\Schub_{\tau_m})$ its graded components.

\begin{lem}\label{lem:primeX2}For $m=1,\dots,\ell$ we have
$\Schub_{\tau_m}=\overline{\NuP{\tau_m}\NuP/\NuP}$, in
  particular, the Schubert varieties $\Schub_{\tau_m}$ are $G$ stable.
\end{lem}
\begin{proof}
Recall that the $\NuP$ orbits
in $\Gr$ are parametrized by $W\backslash \NuW/W$,
and if $w \in \NuW$ we denote by
$[w]_\NuP$ its class in $W\backslash \NuW/W$.
Since $\NuP w \NuP$ is the union of all classes
$\NuB w' \NuP$ with $[w']_\NuP = [w]_\NuP$,
our claim follows from the fact that $[w\tau_i] \leq [\tau_i]$ for all $w\in W$, or equivalently,
from $[w\sutau_i]\geq[\sutau_i]$ for all $w\in W$.
%
\end{proof}

Note that $\grG(\Schub_{\tau_m})$ is a $\NuP$ module, so it is also a $G$-module.
The following two theorems collect the essential properties of
$\grG(\Schub_{\tau_m})$ that we will need
for the constructions in the next section.
We describe first the structure of $\grG_n(\Schub_{\tau_m})$ as a $G$-module.

\begin{teo}\label{teo:sezioniX}Let $m \in \{1,\dots,\ell\}$, then
\begin{enumerate}[\indent $i$)]
\item for any $0\leq i \leq m$ we have $\Ga_i(\Schub_{\tau_m})\isocan V_{\gre_i}^*$ as a
$G$-module and
$\Ga_i(\Schub_{\tau_m})=0$ for any $i>m$;
\item $\grG(\Schub_{\tau_m}, \calL^{\otimes n}) \isocan \bigoplus_{0\leq i_1
  \leq \dots \leq  i_n \leq m}V^*_{\gre_{i_1}+\dots+\gre_{i_n}}$
as a $G$-module.
\end{enumerate}
\end{teo}

\begin{proof}
For $\rho \in \NuW$ define $S(\rho) = \{ \eta \in \NuW/W \st \eta \leq [\rho]\}$
and set
$$
S^+(\rho) = \{ \eta \in S(\rho) \st \eta (\nuom_0)\ristretto _ \got\text{\ is dominant
for the Lie algebra $\gog$}\}.
$$
For $\eta \in \NuW/W$ denote by $\sueta$ the minimal element in $W\eta$. Note
that if $\eta \in S(\rho)$, then $\eta \in S^+(\rho)$ if and only if $\sueta = \eta$,
in particular, $\widehat{[\tau_h]}=[\sutau_h]$. The first step in the proof is to show that
$S^+(\tau_m) =\{[\sutau_0],[\sutau_1],\dots, [\sutau_m]\}$.

Let $\eta \in S^+(\tau_m)$ and suppose $\eta\leq[\tau_h]$ for some
$0\leq h\leq m$.
We want to show that either $\eta=[\sutau_h]$ or $\eta\leq[\tau_{h-1}]$.
Once this is established, our claim follows by
induction on $h$ since $\eta\leq[\tau_m]$ by hypothesis.

Recall first that by Proposition \ref{prp:sollevamenti} $s_0=s_{\gra_0}$ does not appear in any
reduced expression for $\ts_1,\ts_2,\ldots,\ts_\ell$.
Hence there exists a reduced expression $s_0 s_{\grb_1}s_{\grb_2}\cdots s_{\grb_q}$ with $\grb_i
\in \NugrD$ for $\sutau_h$, and in turn there exists a reduced expression
$s_{\grg_1}s_{\grg_2}\cdots s_{\grg_p}s_0s_{\grb_1}s_{\grb_2}\cdots s_{\grb_q}$
{\color{mycolorP}for $\tau_h$} with $\grg_i \in \grD$ {\color{mycolorP} for all $1\leq i\leq p$}.

Also recall that $\eta$ is the minimal element in $W \eta$, so if $\eta\neq [e]$ and if $w$ is the
minimal element in $\NuW$ such that $[w]=\eta$ any reduced expression for $w$ must start with
$s_0$, let's say $s_0s_{\grd_1}s_{\grd_2}\cdots s_{\grd_r}$ is such an expression. By the
characterization of the Bruhat order in terms of subwords and since $\eta\leq[\tau_h]$, we can
choose the decomposition of $w$ such that $s_0s_{\grd_1}s_{\grd_2}\cdots s_{\grd_r}$ is a subword
of {\color{mycolorP} $s_{\grg_1}s_{\grg_2}\cdots s_{\grg_p}s_0s_{\grb_1}s_{\grb_2}\cdots
s_{\grb_q}$}. But {\color{mycolorP}$\grg_i \in \grD$ for all $1\leq i\leq p$}, hence
$s_0s_{\grd_1}\cdots s_{\grd_r}$ is a subword of $s_0s_{\grb_1}\cdots s_{\grb_q}=\sutau_h$;
{\color{mycolorP} this shows that $\eta\leq[\sutau_h]$ as elements of $\NuW/W$.}

{\color{mycolorP} {\color{mycolor}Next we show} that $[s_0\sutau_h]$ is the unique element in
$\NuW/W$ covered by $[\sutau_h]$ with respect to the {\color{mycolor}Bruhat} order. (If $a,b$ are
elements of a partially ordered set we say that $a$ \emph{covers} $b$ if $a>b$ and $a>c\geq b$
implies $c=b$.) {\color{mycolor} Recall that $\kappa' <\kappa$ for $\kappa,\kappa'\in \NuW/W$ if
and only if for the corresponding Demazure modules in $Z$ we have $Y_{\kappa'}\subseteq Y_\kappa$.
For $\sutau_h$ the Demazure module is generated by an extremal weight vector $v_{\sutau}$ of weight
$\sutau_h(\nuom_0) \ristretto_{\Nut'}=(\gre_h-\nuom_0)\ristretto_{\Nut'}$. It follows that $e_\al
v_{\sutau}=0$ for all root operators corresponding to a simple root $\al\not=\al_0$, and $e_{\al_0}
v_{\sutau}=v_{s_{\al_0}\sutau}$ is a generator for the Demazure module $Y_{s_{\al_0}\sutau}$. This
shows that a Demazure module properly contained in $Y_{\sutau}$ is also contained in
$Y_{s_{\al_0}\sutau}$, which proves the claim.} So we can now conclude that $\eta=[\sutau_h]$ or
$\eta\leq [s_0\sutau_h]\leq [\tau_{h-1}]$ and the claimed description of $S^+(\tau_m)$ is proved.}

Now we prove $i)$ using the LS--path branching rule \cite{L1}. Let $\B$ be the
LS--path model for the $\Nug$--module $Z$ and let $\B(\tau_m)$ be the path submodel for the
$\NuP$--module $\Ga(\Schub_{\tau_m})$ and recall that
$$
\Res_G^P\Ga(\Schub_{\tau_m},\calL^{\otimes n})\simeq\oplus_\pi V_{\pi(1)\ristretto_{\got}}^*
$$
where the sum runs over all LS--paths $\pi\in\mB(\tau_m)$ of degree $n$ such that
$\pi(x)\ristretto_{\got}$ belongs to the
dominant Weyl chamber of $\gog$ for all $0\leq x\leq 1$. Let us write such a path as
$\pi=\pi_1*\cdots*\pi_r$ with $\pi_i = \pi_{a_i\eta_h(\nuom_0)}$ for some elements
$\eta_1 <\dots < \eta_r$ in
$S(\tau_m)$ and some rational numbers $0 < a_1,\ldots,a_r$
such that $a_1+\dots+a_r = n$.
The requirement $\pi(x)\ristretto_\got$ dominant for all $x$
implies $\eta_r(\nuom_0) \ristretto _\got$ dominant or equivalently
$\eta_r \in S^+(\tau_m)$, so $\eta_r = [\sutau_h]$
for some $0\leq h\leq m$.

Now the requirement for $\pi$ to be a LS path implies that
$a_{r-1}\langle \sutau_h(\nuom_0), \gra_0\cech\rangle\in \mZ$.
But $\langle \sutau_h(\nuom_0), \gra_0\cech\rangle$ is equal to  $-1$ if $h>0$
and  to $1$ if $h=0$, so if $n=1$ this implies $a_r =1$ and $r=1$, $\pi = \pi_{\eta_h(\nuom_0)}$
and $\pi(1)\ristretto_\got =\gre_h$ which prove our claim
since $\langle\sutau_h(\nuom_0) , D \rangle = -h$ so $V_{\gre_h}^*$ is in degree $h$.

To simplify the presentation, we prove $ii)$ only in the case $n=2$, the proof for the
general case is completely analogous. In this case we can have $a_r=2$ and $r=1$,
$\pi = \pi_{2\eta_h(\nuom_0)}$ and $\pi(1)\ristretto_\got =2\gre_h$ or
$a_r=1$ and $r>1$. In this second case the requirement $\pi(x)\ristretto _\got$ dominant for all $x$
implies $(a_{r-1}\eta_{r-1}(\nuom_0)+ \sutau_h(\nuom_0) )\ristretto _\got$ is dominant.
Now note that if $a,b>0$ and $\eta<[\sutau_h]$ are such that
$(a \eta (\nuom_0)+ b \sutau_h(\nuom_0) )\ristretto _\got$ is dominant, then
$\eta \in S^+(\tau_m)$. Indeed, if $\gra \in \grD$, then
$\langle \sutau_h(\nuom_0), \gra\cech\rangle \neq 0$ implies
$\tal = \tal_h$. So it is enough to prove that
$\langle\eta(\nuom_0), \gra\cech \rangle \geq 0$ for all $\gra\in \grD$ such that $\tal = \tal_h$.
By construction we have $\sutau_h(\nuom_0) = \nuom_0 - \sum_{i\leq h-1}a_i\tal_i$ with
$a_i\in \mN$ and $\sutau_h(\nuom_0) = \nuom_0 - \sum_{\gra\in \NugrD\st\tal\neq \tal_h } b_\gra \gra$ with
$b_\gra \in \mN$. So if $\eta < [\sutau_h]$, then we must have
$\eta (\nuom_0) = \nuom_0 - \sum_{\gra\in \NugrD\st\tal\neq \tal_h } c_\gra \gra$, where $c_\gra \in \mN$.
In particular, $\langle\eta(\nuom_0), \gra\cech \rangle \geq 0$ for all $\gra\in \grD$ such that $\tal = \tal_h$.
\end{proof}

The previous theorem will be be more convenient for us in the
following form. For $m=1,\dots,\ell$ define $\Rich_m$ as the Richardson
subvariety of $\Schub_{\tau_m}$ defined by $z_0^*=0$ and set also
$\Rich = \Rich_\ell$.

\begin{cor} \label{cor:sezioniR}
For $m=1,\dots,\ell$ we have the following isomorphism of $G$-modules
$$\grG(\Rich_m, \calL^{\otimes n}) \isocan \bigoplus_{1\leq i_1
  \leq \dots \leq  i_n \leq m}V^*_{\gre_{i_1}+\dots+\gre_{i_n}}.$$
In particular $\grG_\Rich = \bigoplus _{n\geq 0} \grG(\Rich,\calL^n)$
  is isomorphic to $\mk[G/H]$ as a $G$-module.
\end{cor}

We need the following simple result in the proof of the next theorem:

\begin{lem} The module $V_{\gre_{i+1}}$ appears with multiplicity one
in the tensor product $V_{\gre_1}\otimes V_{\gre_i}$ for $i=1,2,\ldots,\ell-1$.
\end{lem}
\begin{proof}
Let us denote by $\B$ a path model for the $G$--module $V_{\gre_k}$ and denote by $\pi_{\gre_1}$ the
path $\Q\ni t\mapsto t\gre_1\in\La\otimes\Q$. We have the path tensor product formula (see
\cite{L1})
$$
V_{\gre_k}\otimes V_{\gre_1}\simeq\oplus V_{\eta(1)+\gre_1}
$$
where the sum runs on all paths $\eta\in\B$ such that the concatenation $\eta*\pi_{\gre_1}$ is
completely contained in the dominant Weyl chamber. So in order to obtain the module $V_{\gre_{k+1}}$
we must look for the paths in $\B$ ending in $\gre_{k+1}-\gre_1$.

Using the same description of restricted roots we have used in the proof of Lemma \ref{lem:primeX},
we have $s_{\tal_1}s_{\tal_2}\cdots s_{\tal_k}(\gre_k)=\gre_{k+1}-\gre_1$.

Since the restricted Weyl group is a quotient of a subgroup of the Weyl group of $G$ we have proved
that the weight $\gre_{k+1}-\gre_1$ is an extremal weight for the $G$--module $V_{\gre_k}$. This shows
that exactly one path in $\B$ ends in $\gre_{k+1}-\gre_1$ and finishes our proof.
\end{proof}

The non-vanishing of the following specific vector will be important for us in the next section.
Recall that by Proposition \ref{prp:primeNug} we have $\Nug_{1}\isocan V_{\gre_1}^*$,
so we can choose a spherical vector $h_{1} \in \Nug_{1}$.

\begin{teo}\label{teo:hnonsvanisce}
For any $0\leq i\leq\ell$ the element $h_1^i \cdot z^*_0 \ristretto_{\Schub_\tau}$ is a nonzero
section in $\grG_i(\Schub_\tau)$.
\end{teo}

\begin{proof}
Consider the enveloping algebra of $\Nug_+$: $\blU_+ = \blU(\Nug_+)$.
Notice that it is generated by $\Nug_1 \subset \Nug_+ \subset \blU_+$ and
that the map from $\blU_+$ to $Z^*$ given by $x \mapsto x\cdot z_0^*$ is surjective.
Moreover, $\blU_+$ and $Z^*$ are compatibly graded, hence for all $n>0$
we have a surjective morphism:
$$
\Nug_1^{\otimes n} \lra Z^*_n \quad \text{ given by } \quad
x_1\otimes\dots\otimes x_n \longmapsto x_1\cdot(x_2 \cdot (\dots x_n\cdot z_0^*)).
$$
Similarly we have a surjective map from $\Nug_1^{\otimes n}$ onto $\grG_n(\Schub_\tau)$
and by induction a surjective map
$$
a :  \Nug_1\otimes \grG_i(\Schub_\tau)\lra \grG_{i+1}(\Schub_\tau)
\quad \text{ given by } \quad x \otimes v \longmapsto x\cdot v.
$$
Now $\grG_i(\Schub_\tau)\isocan V_{\gre_i}^*$ and $\Nug_1 \isocan V_{\gre_1}^*$.
By the previous lemma, the multiplicity of $V_{\gre_{i+1}}^*$ in
$ V_{\gre_1}^*\otimes V_{\gre_i}^*$ is one. Since $a$ is $G$-equivariant,
the morphism $a$ must be equal to the projection $\pi^{\gre_1,\gre_i}_{\gre_{i+1}}$.
In particular, $a(h_{1} \otimes h_{\gre_i})\neq 0$ by Corollary \ref{cor:hallan}.
The image is $H$-invariant and must hence be a nonzero multiple of $h_{\gre_{i+1}}$,
which proves the claim by induction.
\end{proof}

\section{The equations of the symmetric variety}\label{sez:equazioni}

In this section we describe the relation between the symmetric space $G/H$
and the Grassmannian $\Gr$. The naive approach is the following:
let $h_{\meno 1} \in \Nug_{\meno 1}$ be fixed by $H$ as in \ref{ssec:triple}
and define $x = e^{h_{\meno 1}} (\mk z_0) \in \mP(Z)$. The point $x$ is
certainly fixed by $H$ (since both $h_{\meno 1}$ and $z_0$ are fixed by
$H$). So we can define an immersion $G/H \lra \Gr$ by $gH
\mapsto g x$, and we deduce the defining equations for $G/H$ from the
defining equations for $\Gr$.

Of course, this naive approach has a problem since the exponential map is
not defined for all elements in the Lie algebra in the affine case. Nevertheless,
the reader should keep this simple idea as a travel guide in mind.
To make the idea work despite the obvious mistake we have to go a
sometimes rather technical looking detour.

\subsection{The completion of $\calU^-$ and some notation for
Schubert varieties} In order to define $e^{h_{\meno 1}}$ we introduce a
completion of the negative unipotent subgroup of $\NuG$. Let
$\Nub^-$ be the Lie algebra of the Borel defined by the negative
roots and let $\NuB^-$ be the associated Borel subgroup. We define
$\Nun^-$ as the nilpotent radical of $\Nub^-$ and $\NuU^-$ as the
unipotent radical of $\NuB^-$. Also we denote by $\Nusun^-$ the
pro-Lie algebra $\prod_{\gra \in \NuPhi^-}\Nug_\gra$ and we define
$\NusuU^- \defi \exp(\Nusun^-)$ (see \cite{Kumar} pg. 221) and we
have an inclusion $\NuU^- \incluso \NusuU^-$.
In particular $e^{h_{\meno 1}}$ is an element of $\NusuU^-$.

The group $\NusuU^-$ does not act on $Z$ but for all finite
codimensional $\NuU^-$ submodules $J$ of $Z$  the action
of $\NuU^-$ on $Z/J$ extends uniquely to an action of $\NusuU^-$.
Moreover, if $J$ is $G$ stable, then the orbit $\NusuU^- z_0$ in $Z/J$ is
also stable by the action of $G$.

Let $\NuP^-$ be the parabolic subgroup opposite to $\NuP$. The
$\NuG$-orbit of the line $\mk \, z_0^*$ in $\mP(Z^*)$ is isomorphic to
$\NuG/\NuP^-$. For an element $\eta$ of the Weyl group let $\Schub_\eta$ be the
Schubert variety $\overline{\NuB\,\eta\,\NuP/\NuP}$ and denote by $\Schub\cech_\eta$ the Schubert
variety $\overline{\NuB^-\,\eta\,\NuP^-/\NuP^-}$. Let $Y_\eta \subset Z$ be the associated
Demazure module, i.e.,  $Y_\eta$ is the vector subspace of $Z$ generated by the cone
over $\Schub_\eta$. Similarly, let $Y\cech_\eta \subset Z^*$ be the associated Demazure
module. Denote by $J_\eta \subset Z$ (resp. $J\cech_\eta\subset Z^*$) the annihilator of
$Y\cech_\eta$ (resp. $Y_\eta$). Then $J_\eta$ is a $\NuU^-$ stable
complement of $Y_\eta$, and if $\Schub_\eta$ is $G$ stable, then $J_\eta$ is
also $G$ stable.

For an element $\eta$ of the Weyl group $\NuW $ we denote
by $\suA_\eta$ the orbit $\NusuU^- (\mk z_0) \subset Z/J_\eta$. If
$\Schub_{\eta'} \subset \Schub_{\eta}$, then we have an inclusion
$J_{\eta'} \supset J_{\eta}$ of the annihilators.
Denote by $p^\eta_{\eta'}$ the projection
$$
p^\eta_{\eta'}: \mP(Z/J_\eta)\senza \mP(J_{\eta'}/J_\eta)\longrightarrow\mP(Z/J_{\eta'})
$$
Note that $\suA_\eta \subset\mP(Z/J_\eta) \senza \mP(J_{\eta'}/J_\eta)$, so
$p^\eta_{\eta'}$ is well defined on $\suA_\eta$. Let
$A:=\NuU^- (\mk z_0)\subset \Gr$ be the open cell and set
$A_\eta=A\cap\Schub_\eta$.
The projection from $\mP(Z) \senza \mP(J_{\eta})$ to
$\mP(Z/J_{\eta})$ becomes an isomorphism when restricted to $A_\eta$, and
its image is contained in $\suA_\eta$.

\subsection{The immersion $\mi_\eta$}
Let $\Schub_\eta$ be the closure of a $\NuP$-orbit and set
$x_\eta\defi e^{h_{\meno 1}}(\mk z_0) \in \mP(Z/J_\eta)$. Consider the
$G$-equivariant map
\begin{equation}\label{eq:mi1}
\mi_\eta : G/H \lra \suA_\eta\subset \mP(Z/J_\eta); \quad
\mi_\eta(gH)= g x_\eta.
\end{equation}
Since $e^{h_{\meno 1}}z_0$ is fixed by $H$ the pull back
$\mi_{\eta}^*(\calO_{\mP(Z/J_\eta)})$ on $G/H$ is trivial and
we have an induced map $\mi^*_{\eta} : (Z/J_\eta)^* \isocan Y\cech_\eta \lra \mk[G/H]$. We
can normalize this map in such a way that $\mi^*_{\eta}((z^*_0))$ is
the constant function with value $1$ on $G/H$.
Note that if $\Schub_{\eta'} \subset \Schub_{\eta}$, then we have the following
commutative diagram
$$
\xymatrix{ G/H \ar[r]^{\mi_{\eta}} \ar[rd]_{\mi_{\eta'}}& \suA_{\eta}\ar@{}[r]|{\subset}
 \ar[d]^{p^\eta_{\eta'}} &
\mP(Z/J_\eta) \\
& \suA_{\eta'} \ar @{} [r]|{\subset} & \mP(Z/J_{\eta'})
 }
$$
If we normalize the pull back $(p^\eta_{\eta'})^*$ of the projections
$p^\eta_{\eta'}$ to map $z_0^*$ into $z_0^*$, then $(p^\eta_{\eta'})^*$
restricted to $Y\cech_{\eta'}$ is just given by the inclusion
$Y\cech_{\eta'} \subset Y\cech_{\eta}$. So if $f\in Y\cech_\eta$, then
$\mi^*_{\eta}(f)= \mi^*_{\eta'}(f)$ and we can define $\mi^*:
Z^*\lra \mC[G/H]$ as the limit of the maps $\mi^*_{\eta}$. Consider
the morphism of rings $S\mi^* : S(Z^*) \lra \mk[G/H]$ given by the
symmetric product of the map $\mi$. Recall that the ring
$\grG_{\Gr}\defi\Gamma_\calL(\Gr)=\bigoplus_{n\geq 0} \grG(\Gr, \calL^n)$
is a quotient of $S(Z^*)$, and let
$I\subset S(Z^*)$ be the ideal defining $\grG_{\Gr}$.

\begin{lem}\label{lem:Iphi}
$S\mi^*(I)=0$, so $S\mi^*$ determines a morphism of rings $\grf
:\grG_{\Gr} \lra \mk[G/H]$.
\end{lem}
\begin{proof}
Let $f \in I$. We can assume that $f$ is a homogeneous element
contained in the symmetric product of $Y\cech_\eta$ for an appropriate
$\eta$ so that $S\mi^*(f) = S\mi_\eta^*(f)$. We want to prove that
$f(\mi_\eta(x))=0$ for all $x \in G/H$. By \cite{Kumar} \S VII.3
there exists a Schubert variety $\Schub_\vartheta$ such that we have
$p^{\vartheta}_\eta(A_\vartheta) = \suA_\eta$. So let $y \in A_\vartheta$ be
such that $p^{\vartheta}_\eta (y)=\mi_\eta(x)$. Then
$f(y)=0$ since $y \in \Gr$. But notice that $f(y) = f(p^{\vartheta}_\eta (y))$
since $f$ in the symmetric product of $Y\cech_\eta$, so it is zero on $J_\eta$.
\end{proof}

\subsection{Standard monomial theory for $G/H$}
Now we use the morphism $\grf$ and the SMT for the ring $\Gamma_\Gr$
(see section \ref{ssec:SMTflag}) to construct a SMT for the ring
$\mk[G/H]$. Let $\mF = \mF_\calL$ be the basis of
$\Gamma(\Gr,\calL)=Z^*$ constructed in \cite{L:SMT} and denote
by $<$ the order on $\mF$. The construction can be fixed such that
$f_0 = z_0^*$ is the minimal element in $\mF$. Denote by $\mS\mM$
(respectively by $\mM$) the set of standard monomials (respectively the set of
monomials) in the elements of $\mF$. For $f\in \mF$ set
$g_f=\grf(f)$, we define similarly $g_m$ for $m\in \mM$.

For an element $\eta\in\NuW$ we define $\mF(\eta) =\{f \in \mF \st f\ristretto_{\calS_\eta}\neq
0\}$ and $\mF_0(\eta) =\mF(\eta)\senza\{f_0\}$. If $\eta$ is the special element $\eta=\tau$ (see
section~\ref{specialtau}), then we denote $\mF_0(\eta)$ just by $\mF_0$. Recall that $\mF_0=\{f\in
\mF \st f\ristretto_\Rich \neq 0\}$. Let $\mS\mM_0$ (respectively $\mM_0$) be the set of all
standard monomials (respectively all monomials) in the elements of $\mF_0$. Recall that, as in
\ref{ssec:SMTflag}, the set $\{m\ristretto_\Rich \st m \in \mS\mM_0\}$ is a $\mk$-basis of
$\grG_{\Rich}=\bigoplus_{n\geq 0}\grG(\Rich,\calL^n)$.

We are finally ready to apply all the various technical results of this and previous sections and
to conclude with our main theorem.

By Theorem \ref{teo:hnonsvanisce}, for all $f \in \mF_0$ the functions $g_f$ do not vanish
identically. Hence by Corollary \ref{cor:sezioniR} the set
$$
\mG_0=\{g_f\st f\in \mF_0\}
$$
is a $\mk$ basis of $\mV\defi V_{\gre_1}^*\oplus \dots\oplus V_{\gre_\ell}^*\subset \mk[G/H]$.
We introduce the following order on $\mG_0$ induced by the order on $\mF_0$: $g_f \prec g_{f'}$ iff $f < f'$.

\begin{teo}\label{teo:smt}
The set $\{g_m \st m \in \mS\mM _0\}$ is a basis of $\mk[G/H]$, hence
$(\mG_0,\prec)$ is a SMT for the ring $\mk[G/H]$.
\end{teo}

\begin{proof}
Let $E$ be the span in $\grG_{\Gr}$ of the monomials $g_m$ with
$m \in \mS\mM _0$, this set is $G$ stable. By
Corollary \ref{cor:sezioniR}, $E$ is isomorphic to $\bigoplus_{\grl\in \Omega} V_\grl^*$
as a $G$-module. Let $E_\grl$ be the $G$-submodule of $E$ isomorphic to $V^*_\grl$.
By Theorem \ref{teo:hnonsvanisce}, we know that $\grf(E_{\gre_i})\neq 0$, and hence
also $\grf(E_\grl) \neq 0$ because $\mk[G/H]$ is a domain (the product of the two highest weight
vectors in $E_\mu$ and $E_\nu$ is an highest weight vector in
$E_{\mu+\nu}$). So $\grf \ristretto_E$ is injective, and by the descriptions
of $\mk[G/H]$ and $\grG_\Rich$ as $G$-modules (Corollary~\ref{cor:sezioniR}),
it follows that the map is surjective.
\end{proof}

\begin{oss}\label{oss:smt}
In section \ref{ssec:SMTX} we gave a description for a SMT for the ring
$\grG_{\Xbar}$. In particular, by the description of $\mk[G/H]$ as the
quotient $\grG_{\Xbar} / (s_i=1)$ we obtain a set of generators of
$\mk[G/H]$. These generators coincide with the functions $g_f \in \mG_0$.
This follows from the fact that the $G$-modules we are considering
(the submodules $V_{\gre_i}$ of $Z^*_{i}$) are generated by extremal weight
vector of the modules $Z^*$, by the construction of the SMT in \cite{L:SMT} and
standard arguments.

Also it is not difficult to prove that the SMT of Theorem~\ref{teo:smt} is compatible with
$G$-modules in the following sense: there exists a filtration of $\mk[G/H]$ by $G$-modules $F_i$
with simple quotients such that for all $i$ the set $\{g_m\st m\in \mS\mM_0\}\cap F_i$ is a
$\mk$-basis of $F_i$.
\end{oss}

\subsection{Straightening relations for $\mk[G/H]$}
We describe now  straightening relations for the standard
monomial theory using the Pl\"ucker relations for the Grassmannian.
We denote by $<_t$ the total order on  $\mM$ and for
$f,f' \in \mF$ not comparable let $R_{f,f'} = f \, f' - P_{f,f'}\in I \cap S^2(Z^*)$
be the Pl\"ucker relation as in section \ref{ssec:SMTflag}.

Let $\mk[u]=\mk[u_f\mid f\in \mF_0]$ be the polynomial ring with generators
indexed by the elements of $\mF_0$. For a monomial $m=f_1\cdots f_s\in\mM_0$ let $u_m
=u_{f_1}\cdots u_{f_s}$ be the corresponding monomial in $\mk[u]$.
Denote by $\psi$ the morphism of rings from the polynomial algebra
$\mk[u]$ to $\mk[G/H]$ defined by $\psi(u_f)=g_f$ and let $Rel$ be the kernel
of this morphism.

We introduce on $\mk[u]$ a degree: for $f \in \mF_0 (\tau_i) \senza\mF_0(\tau_{i-1})$ let $u_f$ be of degree $i$
and we indicate by $gr(r)$ the degree of an element $r$ in $\mk[u]$. If $m,m' \in \mM_0$,
then we define $u_m\prec_t u_{m'}$ if $gr(u_m) < gr(u_{m'})$ or if  $gr(u_m) = gr(u_{m'})$ and $m <_t m'$.

This order has the properties explained in section \ref{ssec:SMT}. The compatibility of this order
with the order $\prec$ on $\mG_0$ follows from the compatibility of the order $<$ between
elements of $\mF$ with the dominant order of the
associated weights recalled in section \ref{ssec:SMTflag}.

Fix an element $\eta \geq \tau$ such that for all $f,f'\in \mF_0$
that are not comparable, the relation  $R_{f,f'}$ is in $S^2(Y\cech_\eta)$. Equivalently:
$P_{f,f'}$ is a polynomial in the functions in $\mF(\eta)$.
We define $\mF_1 = \mF(\eta) \senza \mF_0$.

For each $f \in \mF$ let $n_f = -\langle D, weight(f)\rangle +n_0$.   
Recall that $weight(f)$ is the weight of $f$ with respect to $\Nut$, so
$f \in (Z^*)_{n_f}$. Note that the set $\{f \st f \notin \mF_0, \mand n_f= n\}$ is a
$G$-stable complement for $V_{\gre_n}^*$
in $(Z^*)_n$ for $n=1,\dots,\ell$ and is equal to $(Z^*)_n$
otherwise. Hence if $f \in \mF_1$, then by Corollary \ref{cor:GmoduliZn}
we have
$$
g_f \in \bigoplus_{\grl \in \Omega \mand gr(\grl)< n_f} V_{\grl}^*.
$$
In particular, for each $f \in \mF_1$ we can choose an element
$F_f(u) \in \mk[u]$ such that $gr(F_f) < n_f$ and  such that
$\psi(F_f) = F_f((g_{f'})_{f'\in \mF_0}) =g_f$. We set also
$F(u)=(F_f(u))_{f\in \mF_1}$.

\begin{oss}\label{oss:polF}
The computation of the polynomials $F$ depends only on the expansion of $e^{h_{-1}} z^*_0$
and on the representation theory of $G$ and not anymore on the geometry of $G/H$.
Indeed, once $e^{h_{-1}}z_0^*$ is computed, we can
determine the map $\grf$, hence the decomposition of the functions $g_f$ in the irreducible factors in $\mk[G/H]$
(we have explicit bases of the irreducible modules given for example by the basis in \cite{L:SMT}). Now given an element in
$V_\grl^* \subset \mk[G/H]$, we have $\grl=\sum n_i \gre_i$ and $V_\grl^*$ appears with multiplicity one in
the tensor product $TP=V_{\gre_1}^{\otimes n_1}\otimes \cdots \otimes V_{\gre_1}^{\otimes n_1}$. In particular,
the $G$-equivariant projection $\pi$ from $TP$ to $V_\grl^*$ is unique up to scalar. Now consider product the map from $TP$ 
to the ring $\mk[G/H]$ followed by the projecting onto $V_\grl^*$. This is also a $G$-equivariant
non zero map, so it has to coincide with $\pi$ up to a non zero scalar. By fixing highest weight vectors, this scalar
can be normalized to be $1$. So the functions $F_f$ are determined by the decomposition of the tensor product $TP$.
\end{oss}

We will use the set $\mG_1 = \{g_f \st f \in \mF_1\}$  as a set of  auxiliary
variables, so for each $f \in \mF_1$  we introduce a new
variable $v_f$ and we set $v=(v_f)_{f \in \mF_1}$.

For non-comparable elements  $f,f' \in \mF_0$
we have the polynomials $R_{f,f'}$ and $P_{f,f'}$ in
the symmetric algebra $\sfS(\mF_0\cup \mF_1)$.
Let $R_{f,f'}(u,v)$ and $P_{f,f'}(u,v)$ be the polynomials
obtained by substituting an element $h\in \mF_0\cup \mF_1$
by $u_h$ if $h \in \mF_0$ and $v_h$ if $h \in \mF_1$,
so $R_{f,f'}(u,v)= u_f\,u_{f'} - P_{f,f'}(u,v)$.
Note that $P_{f,f'}(u,v)$ is a homogeneous polynomial of degree
two which is the sum of monomials of the form
$u_{f_1}\,u_{f_2}$ or $u_{f_1}\,v_{f_2}$ or
or $v_{f_1}\,v_{f_2}$, where $f_1\, f_2 <_t f\,f'$
and $n_{f_1}+n_{f_2} = n_{f}+n_{f'}$ (by the fact that the
relations are $\Nut$ homogeneous).

Let now $\psi_1$ is the morphism of rings from the polynomial ring
$\mk[u,v]$ to $\mk[G/H]$ defined by $\psi_1(u_f)=g_f$ if $f
\in \mF_0$ and $\psi_1(v_f)=g_f$ if $f \in \mF_1$, and let
$Rel_1$ be the kernel of this map. By Lemma \ref{lem:Iphi} and by the definition above,
we have the following equations in $\mk[u,v]$:
\begin{alignat}{3}
v_f              & = F_f(u) \qquad &(\text{mod } Rel_1)  &\qquad&&\mforall f \in \mF_1; \label{eq:vF}\\
R_{f,f'} (u,v)& =0                &(\text{mod } Rel_1)  &&&\mforall
f, f' \in \mF_0 \text{ that are not comparable}. \label{eq:R}
\end{alignat}
Now we can substitute  equations \eqref{eq:vF} in equations
\eqref{eq:R} and define
\begin{align*}
\hat P _{f,f'} (u) &= P_{f,f'} (u, F(u)) \\
\hat R _{f,f'} (u) &= R_{f,f'} (u, F(u)) =
u_{f}\,u_{f'} - \hat P _{f,f'} (u)
\end{align*}
for all $f, f' \in \mF_0$ that are not comparable.  The new polynomials
$\hat R _{f,f'} (u)$ obtained in this way are obviously elements of $Rel\subset \mk[u]$.
More precisely, the following theorem states that these polynomials form a set of
straightening relations.

\begin{teo}\label{teo:equazioni}
The relations $\hat R_{f,f'}(u)$ for $f,f' \in \mF_0$ that are not comparable
are a set of straightening relations for the order $\prec_t$ introduced above.
In particular they generate the ideal $Rel=\ker \psi$ in $\mk[u]$.
\end{teo}

\begin{proof}
We have to prove for all $f,f' \in \mF_0$ that are not comparable:
the polynomial $\hat P _{f,f'} (u)$ is a sum of monomials $u_m \prec
 u_{f}\, u_{f'}$.

Let $u_{f_1}v_{f_2}$ be a monomial which appears in $P_{f,f'} (u,v)$.
Then $gr(u_{f_1}F_{f_2}(u)) = gr(u_{f_1}) + gr(F_{f_2}(u)) < n_{f_1} +n_{f_2} =
n_{f}+n_{f'}$, by the discussion above, so all the monomials
which appears in $u_{f_1}F_{f_2}(u)$ are $\prec$ of
$u_{f}u_{f'}$. Similarly we can treat the monomials
$v_{f_1}v_{f_2}$. Finally the monomials
$u_{f_1}u_{f_2}$ which appear in $P _{f,f'} (u,v)$ are
such that $f_1 \,f_2 <_t f\,f'$ so  $u_{f_1}u_{f_2}\prec
u_{f}u_{f'}$.
This proves that the relations $\hat R_{f,f'}(u)$ for $f,f' \in \mF_0$
are a set of straightening relations. The second part of the statement follows now by
Theorem \ref{teo:smt} and Lemma~\ref{lem:SMT}.
\end{proof}

Despite the fact that the computation of the polynomials $F_f$ depends only on the expansion of
$e^{h_{-1}}$ and the representation theory of $G$, it seems complicated to get explicit
formulas and check basic properties for these polynomials. For example, by Corollary~\ref{cor:quadratic}
we know that the relations in the generators $\mG_0$ are quadratic. However, a priori the
relations $\hat R_{f,f'}$ can be of higher degree. From this point of view it is natural to ask whether it is
possible to fix $\eta\ge\tau$ such that the functions $F_f$ can be chosen to be linear in the generators,
in this case it would be clear that the relations $\hat R_{f,f'}$ are quadratic.

A more precise way to state this is the following: let $\eta$ be
minimal such that $p^{\eta}_\tau(\Schub_\eta)$ contains $x_\tau$.
Is it true that any spherical module in $\grG(\Schub_\eta,\calL)$ is
one of the modules $V_{\gre_i}^*$? In the last
section we show that in the case $\Nug$ is of finite type this
question has an affirmative answer.

\begin{oss}
We have seen above that the coordinate ring of $G/H$ and $\grG_\Rich$ have similar properties.
Indeed, we can perform a two steps flat and $G$-equivariant deformation of $\mk[G/H]$ to $\grG_\Rich$.
Let $\grG_\tau = \bigoplus _{n\geq 0}\grG(\calS_\tau,\calL^n)$ and define $A$ to be the quotient of
$\grG_\tau$ modulo the ideal generated by $(f_0-1)$. It is clear that the ring $A$ can be deformed to
$\grG_\Rich$ in a flat and $G$ equivariant way. We  exhibit now a deformation of $\mk[G/H]$ to $A$. To
this order we need first to change the choice of our generators $\mF_1$.

Let $\mF'_1$ be a set of elements such that:
 \begin{enumerate}[\noindent i)]
   \item $f_0\in \mF'_1$;
   \item $\mF'_1$ is a basis of the vector space generated by $\mF_1$;
   \item the elements of $\mF'_1$ are $\Nut$ homogeneous and compatible with $G$-modules; in
particular for each $f$ in $\mF'_1$ there exists an irreducible submodule $M$ of $\grG(\Gr,\calL)$
such that $f \in M$, and let $M\isocan V_{\grl_f}^*$. If $\grl_f = \sum a_i \gre_i$ and $f \in Z_n$
then we define also $\tilde{n}_f=n-\sum i \,a_i$ and notice that this number is bigger than $0$ if
$f\neq f_0$.
 \end{enumerate}
  Notice that conditions $i)$ and $ii)$ are compatibles since the vector space spanned by $\mF_1$ is
$G$ stable and $\Nut$ homogeneous.
  With this choice of generators for each $f \in \mF_1'$ then $\grf(f)$ is in the image of the product
$$m: S^{a_1}(V^*_{\gre_1})\otimes \cdots \otimes S^{a_\ell}(V^*_{\gre_\ell})\lra \mk[G/H]$$
where $\grl_f = \sum a_i \gre_i$. In particular there exists an element $F'_f \in
S^{a_1}(V^*_{\gre_1})\otimes \cdots \otimes S^{a_\ell}(V^*_{\gre_\ell})$ such that
$m(F'_f)=\grf(f)$. We consider $F'_f$ as a multihomogeneous  polynomial in the variables $f\in
\mF_0$.

  Finally notice that the old basis $\mF_1$ can be written in terms of the basis $\mF'_1$.
So we can write the relations $R_{f,f'}$ with respect to this new basis by expressing the elements
in $\mF_1$ as a linear combinations of elements of $\mF_1'$. We call these relations $R'_{f,f'}$.

   Now consider $u$ a set of variables as in the previous discussion and a set of new variables
$v'=(v'_{f'})_{f'\in \mF'_1}$. Consider now in the polynomial ring $\mk[u,v,t]$ the ideal generated
by $R'_{f,f'}(u,v'))$ for $f,f'\in \mF_0$ not comparable and by the elements $v_f - t^{\tilde{n}_f}
F'_f(u)$; let $B$ be the quotient of $\mk[u,v,t]$ under this ideal and finally for $a \in \mC^*$
let $B_a=B/(t-a)$.

   Now notice that there is a $\mk^*$-action on $B$ defined for all $z\in \mk^*$ by $z\cdot u_f = z^n u_f $ if $f \in
\mF_0\cap Z_n$ and by $z\cdot v_f= z^n v_f$ if $f \in \mF'_1\cap Z_n$. Finally notice that $B_0
\isocan A$ and that $B_1 \isocan \mk[G/H]$. In particular $B$ gives the claimed flat deformation
from $\mk[G/H]$ to $A$.

\end{oss}


\section{The finite case}\label{sez:finito}
In the case $\Nug$ is of finite type (or equivalently by Proposition \ref{prp:finitoaffine}:
when $\Phi$ is of type $\sfA_\ell$) part of the proof and construction described in the previous
paragraphs can be simplified and also some other additional properties hold. In this section we
describe some of these special properties.

\begin{prp}\hfill

\begin{enumerate}[\indent i)]
\item $\Schub_\tau$ is a codimension one Schubert variety in $\Gr$;
\item $\grG_i(\Gr) = \grG_i(\Schub_\tau)$ for $i=0,\dots,\ell-1$;
\item $\grG_\ell(\Gr) = \mk$ and $\grG_i(\Gr) = 0$ for $i>\ell$;
\item $h_{1}^i{z^*_0} \neq 0$ for all $i=0,\dots,\ell$.
\end{enumerate}
\end{prp}

\begin{proof}
To prove $i)$ it is enough to show that $[s_0 \tau] = [w_{\NugrD} ]$ in $\NuW/W$ or
equivalently, since $s_0 = \ts_0$ that $\tau (\om_0) = \om_0 -\gre_1$.
This is a computation essentially in the restricted
root system that in this case we know to be of type $\sfA_\ell$. We have
$\tau = w_\grD \sutau$ so by Lemma \ref{lem:primeX} we  have
$\tau(\om_0) = w_\grD(\gre_\ell-\om_0) = w_\grD (\gre_\ell)- \om_0$.
Now $w_\grD (\gre_\ell)= - \gre_1= \gre_\ell- (\tal_1 +\dots +\tal_\ell)$
so $\tau (\om_0) = \om_0 -\gre_1$.

$ii)$ and $iii)$ follows immediately and $iv)$ follows using Lemma \ref{lem:tripla}.
\end{proof}

\begin{oss}
In the case the restricted root system is of type $\sfB$, $\sfC$ or $\sfB\sfC$
a similar computation gives $\tau(\om_0)= - \gre_\ell+3\om_0$.
\end{oss}

The theory developed in the previous section becomes
particularly simple in this case and we restate parts of Theorem \ref{teo:smt}
and Theorem \ref{teo:equazioni} in the following more explicit way.

\begin{teo}\label{teo:finito}
\begin{enumerate}[\noindent $i$)]
\item $\mF=\mF_0 \cup \{f_0,f_1\}$ where $f_1$ is an
 highest weight vector in $Z^*$;
\item we can normalize $f_{0}$ and $f_{1}$ in such a way that $F_{f_0}
  = F_{f_1} = 1$;
\item The map $\grf$ induces the following isomorphism:
$$ \mk[G/H]\isocan \frac{\grG_{\Gr}} { (f_{0}=f_{1}=1) }. $$
\end{enumerate}
\end{teo}

It should also be pointed out that for some of the involutions in
which the restricted root system is of type $\sfA$ the
results described where already obtained as special cases by other
authors: in particular in the case of the group $SL(n)$ this was
obtained by De Concini, Eisenbud and Procesi in \cite{DEP} and for
symmetric quadrics by Strickland \cite{Str} and Musili \cite{Musili,Musili2}.
There are other three families of involutions in which the restricted root system is of type
$\sfA$: the case of antisymmetric quadrics, a family for the group
$SO(n)$ in which the restricted root system is of type $\sfA_1$ and
that for this reason is particularly simple, and the involution of
$\sfE_6$ with fixed point subalgebra of type $\sfF_4$. In the last part
of this section we want to make as explicit as possible
the case of this exceptional involution.
In this case we have that $\NuG$ is of type $\sfE_7$.
In the picture below we have numbered the nodes of
$\sfE_7$ following the notations of the previous sections and we have colored the
nodes according to the Satake diagram of the corresponding involution.
$$
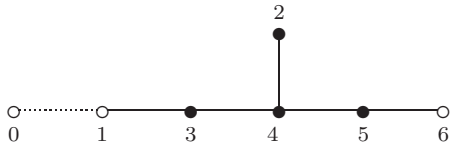
\xymatrix{
  & & & & *{\overset{\mathstrut 2}{\bullet}} \\
   & *{\underset{\mathstrut 0}{\circ}}
           \ar @{.} @<+1.3ex> []+<0.4ex,0ex>;[r]-<0.4ex,0ex>
    & *{\underset{\mathstrut 1}{\circ}}
         \ar @{-} @<+1.3ex> []+<0.4ex,0ex>;[r]-<0.4ex,0ex>
    & *{\underset{\mathstrut 3}{\bullet}} \ar @{-} @<+1.3ex> []+0;[r]-<0.4ex,0ex>
    & *{\underset{\mathstrut 4\;\;}{\bullet}} \ar @{-} @<+1.3ex> []+0;[r]-<0.4ex,0ex>
         \ar @{-}[]+<0ex,1.3ex>;[u]-<0ex,1.2ex>
    & *{\underset{\mathstrut 5}{\bullet}} \ar @{-} @<+1.3ex> []+0;[r]-<0.4ex,0ex>
    & *{\underset{\mathstrut 6}{\circ}}
 }
$$
The module $Z$ is of dimension $56$ and it is a minuscule module
so we can identify an element of the basis $\mF$ by giving its weight.
Also  $S^2 Z \isocan Z_{2\om_0} \oplus \Nug$,
so the Pl\"ucker relations are generated as $\NuG$ modules by the
following single relation:
$$
x_0\,y_0 - x_1\,y_1 + x_2\,y_2 - x_3\,y_3 + x_4\,y_4 -x_5\,y_5=0
$$
where $x_0 = z_0^*$ and
$x_1=\provf_0 (x_0)$,
$x_2=\provf_1 (x_1)$,
$x_3=\provf_3 (x_2)$,
$x_4=\provf_4 (x_3)$,
$x_5=\provf_5 (x_4)$,
$y_5=\provf_2 (x_4)$,
$y_4=\provf_5 (y_5)$,
$y_3=\provf_4 (y_4)$,
$y_2=\provf_3 (y_3)$,
$y_1=\provf_1 (y_2)$ and
$y_0=\provf_0 (y_1)$ where $\provf_i = \provf_{\gra_i}$ are the
Chevalley generators.

\section{Appendix: SMT from above and from below}
We have actually proved the existence of two bases for the coordinate
ring of the symmetric space. One basis is given by the standard monomials $\{g_m \st m \in \mS\mM _0\}$
obtained by restricting the standard monomial theory on the (affine) Grassmannian
to the symmetric space (Theorem~\ref{teo:smt}). 
The other basis comes from below in the following sense: it is obtained via lifting and
pull back from the SMT for the multicone over the closed orbit in the wonderful
compactification (Proposition~\ref{liftsmt}). To be more precise,
in the last case we have a description of $\mk[X_q]$ as the quotient
$$
  \frac{\grG_{\Omega_q}}{(s_i-1 \st i = 1,\dots,\ell)} \isocan\mk[X_q].
$$
By Theorem~\ref{SMTonWonderful}, $\grG_{\Omega_q}$ has as basis monomials
of the form $s^{\mu} m^X $ where $s^{\mu}$ is a product of the $s_i$ and the $m^X$
are appropriate lifts of the standard monomials on the closed orbit $Y$ in the
wonderful compactification (see sections~\ref{ssec:SMTX} and \ref{afirstdescription}).
So the images $\overline{m^X}$ of the $m^X$ also form a basis for the coordinate ring.
We would like to compare these two bases and the two different indexing systems.
\begin{teo}\label{equalbasis}
The possible choices relevant for the construction of the two bases can be arranged
such that the two bases coincide.
\end{teo}
Before we come to the proof, note that
this comparison is also interesting from the combinatorial point of view.
The definition of a standard monomial on a Grassmannian is rather straightforward,
see also~\ref{ssec:SMTflag}. The set of generators of the ring is indexed by certain
LS-paths of shape $\nuom_0$. For details see \cite{L1}, we recall here
only the properties needed in the following. An LS-path of shape $\nuom_0$ is
a pair of sequences $\pi=(\underline{x},\underline{a})$, where
$\underline{x}=(x_1,\ldots,x_r)$ is a strictly increasing sequence
(in the Bruhat order) of elements in $\NuW/\NuW_{\nuom_0}$
(here $\NuW$ is the Weyl group of $\Nug$ and  $\NuW_{\nuom_0}$
is the stabilizer of $\nuom_0$), and
$\underline{a}=(1>a_1>\ldots>a_{r-1}>0)$ is a strictly decreasing sequence of
rational numbers (satisfying certain properties, see  \cite{L1}).

Let  $\eta=(\underline{\kappa},\underline{b})$ be a second LS-path of
shape where $\underline{\kappa}=(\kappa_1,\ldots,\kappa_s)$. We say
\begin{equation}\label{grassmannstandard}
\pi\le \eta\quad \text{if and only if}\quad x_r \le \kappa_1.
\end{equation}
Note that $\pi\ge \kappa$ and $\kappa\ge \pi$ implies $r=s=1$ and hence
$\pi=\eta=(x)$. By definition, a product
\begin{equation}\label{grassmannmonomialstandard}
f_{\pi_1}\cdots f_{\pi_s}\quad \text{is standard if and only if}\quad
\pi_1\le \pi_2\le \ldots\le \pi_s.
\end{equation}
As mentioned in Remark~\ref{oss:smt},  in the multicone picture the definition of
a standard monomial is much more involved. The generators are again indexed
by certain LS-paths, but of a different type.
Let $\ep_1,\ldots,\ep_n$ be the generators of the admissible lattice.
The generators of type $\ep_i$ (see~\ref{ssec:SMTflag}) are indexed by
LS-paths of type $\ep_i$, i.e., pairs of sequences $\pi=(\underline{x},\underline{a})$,
where $\underline{x}=(x_1,\ldots,x_r)$ is a strictly increasing sequence
(in the Bruhat order) of elements in the quotient $W/W_{\ep_i}$, and $\underline{a}$
is a strictly decreasing sequence of rational numbers (satisfying certain conditions, see
\cite{L1}). By a {\it defining sequence} for $\pi$ we mean a weakly increasing sequence
$\underline{\tilde{x}}=(\tilde{x}_1,\ldots,\tilde{x}_r)$
of elements in $W$ such that $\tilde{x}_j\equiv x_j \bmod W_{\ep_i}$.
Given LS-paths $\pi_{1,1},\ldots,\pi_{1,a_1}$,$\ldots$, $\pi_{n,1},\ldots,\pi_{n,a_n}$,
where $\pi_{i,j}=(\underline{x}^{i,j},\underline{a}^{i,j})$ is an LS-path of type $\ep_i$, the monomial
$$
\underbrace{
f_{\pi_{1,1}}\cdots f_{\pi_{1,a_1}}}_{\text{type}\,\ep_1}\cdot
\underbrace{
f_{\pi_{2,1}}\cdots f_{\pi_{2,a_2}}}_{\text{type}\,\ep_2}
\cdots
\underbrace{
f_{\pi_{\ell,1}}\cdots f_{\pi_{\ell,a_\ell}}}_{\text{type}\,\ep_\ell}
$$
is called standard if there exist defining sequences
$\underline{\tilde{x}}^{i,j}$
for the $\pi_{i,j}=(\underline{x}^{i,j},\underline{a}^{i,j})$
such that the defining sequences give rise to a weakly increasing sequence of Weyl group elements:
\begin{equation}\label{multiconestandard}
\underbrace{{\tilde{x}}^{1,1}_1\le {\tilde{x}}^{1,1}_2 \le \ldots \le
\tilde{x}^{1,1}_r}_{\underline{{\tilde{x}}}^{1,1}}\le
\underbrace{\tilde{x}^{1,2}_1\le \cdots\le \tilde{x}^{1,2}_p}_{\underline{{\tilde{x}}}^{1,2}}\le\cdots\le
\underbrace{\tilde{x}^{n,a_n}_1\le  \ldots\le
\tilde{x}^{n,a_n}_s}_{\underline{\tilde{x}}^{n,a_n}}
\end{equation}
This definition of a standard monomial is far away from the definition given
in section~\ref{ssec:SMT}. It depends on the choice of the enumeration
$\ep_1,\ldots,\ep_n$ of the basis of the admissible lattice, and there is no obvious
canonical choice.

In the case where the admissible lattice is the weight lattice, there exist
special "nice enumerations`` for certain groups (see \cite{LMS}, for a Young
diagram like version see \cite{L3}). In these cases the definition above simplifies
dramatically and becomes similar to the one above for the Grassmannian.
The bijection below together with the comparison theorem above gives a beautiful
geometric interpretation of this combinatorial fact and provides yet another connection
between Young tableau like indexing systems and combinatorics of the affine Weyl group.

\begin{proof}
We have already pointed out in Remark~\ref{oss:smt} that the possible
choices for the set of generators can be arranged for
both constructions such that the generators actually coincide.
It remains to prove that the notion of a {\it standard monomial}
coincides for both constructions.

Let us recall a few facts and definitions
related to LS-paths. By Lemma~\ref{lem:primeX}, we can enumerate the basis of the lattice
$\ep_1,\ldots,\ep_\ell$ such that there exist elements in $\NuW$
(the enumeration is different from the one in the lemma above)
\begin{equation}\label{enumerationchoice}
\sutau_1>\sutau_2>\ldots > \sutau_\ell\quad \text{and}\quad \sutau_h(\nuom_0)\vert_{\got}=\ep_h,
\ h=1,\ldots,\ell;
\end{equation}
and the $\sutau_j$ are of minimal length with this property.
Consider first an LS-path $\pi=(\underline{x},\underline{a})$ of type $\ep_i$,
where $\underline{x}=(x_1,\ldots,x_r)$. By abuse of notation we
write also $x_j\in W$ for a minimal representative. By the definition of an LS-path and
by~(\ref{enumerationchoice}) it follows that
$$
\nupi=({^e\!{x}},\underline{a}),\quad\text{where}\quad {^e\!{x}}=(x_1 \sutau_i,\ldots, x_r\sutau_i)
$$
is an LS-path of type $\nuom_0$. So the map $\pi\mapsto \nupi$ defines an injective
(and also surjective) map between the union $\bigcup_{i=1}^\ell\{\text{LS-paths of type}\,\ep_i\}$
and the set of LS-paths standard on the Richardson variety $\Rich$, i.e., the associated sections
do not vanish identically on $\Rich$.

It remains to check that the notion of a standard monomial in both pictures
is the same. To not get drowned in indices, we consider only a product of two
elements. Let $\pi=(\underline{x},\underline{a})$ be of type $\ep_i$
and $\eta=(\underline{y},\underline{b})$ of type $\ep_j$ such that $i>j$ and
$f_\pi f_\eta$ is standard. By definition, this implies that we can find
defining sequences  $(\tilde{x}_1,\ldots,\tilde{x}_r)$ for $\underline{x}=(x_1,\ldots,x_r)$
and  $(\tilde{y}_1,\ldots,\tilde{y}_s)$ for $\underline{y}=(y_1,\ldots,y_s)$ such that in $W$ we have
$$
\tilde{x}_1\le\ldots\le\tilde{x}_r\le \tilde{y}_1\le\ldots\le\tilde{y}_s\quad
\text{and hence in $\NuW$:}\ \tilde{x}_1\sutau_i\le\ldots\le\tilde{x}_r\sutau_i\le
\tilde{y}_1\sutau_j\le\ldots\le\tilde{y}_s\sutau_j.
$$
Recall that an element $\tilde{x}_k$ is of the form $x_kw_k$, where $w_k$ is an element
in the stabilizer $W_{\ep_i}$. Similarly, $\tilde{y}_m$ is of the form $y_m w_m$, where $w_m\in W_{\ep_j}$.
So the linearly ordered sequence in $\NuW$ above gives rise to a linearly ordered sequence
$$
 {x}_1\sutau_i\le\ldots\le {x}_r\sutau_i\le  {y}_1\sutau_j\le\ldots\le {y}_s\sutau_j
\quad\text{in}\quad\NuW/\NuW_{\nuom_0}
$$
Now by (\ref{grassmannmonomialstandard}) this implies:
$f_\nupi f_\nueta$ is a standard monomial.

This argument extends to arbitrary standard monomials on the multicone. Summarizing,
$\mk[X_q]$ has as a basis the standard monomials from below, i.e.,
the classes $\overline{m^X}$, where the $m^X$ are appropriate lifts of the standard
monomials (with respect to the enumeration of the basis of the admissible lattice chosen above)
on the closed orbit $Y$ in the wonderful compactification. The map defined on the set of standard
monomials
$$
f_{\pi_1}\cdots f_{\pi_s}\mapsto f_{\nupi_1}\cdots f_{\nupi_s}
$$
induces a bijection between the standard monomials from below and the standard monomials
from above, i.e., the set $\{g_m \st m \in \mS\mM _0\}$ (Theorem~\ref{teo:smt}).
\end{proof}

\section{Appendix B: The Satake diagrams}

In this Appendix we list the Satake diagrams of all involutions. We add a node to a Satake diagram
as described in the previous sections and we obtain in this way the extended Dynkin diagram; this
special node is joined to the rest of the diagram with a dotted line (or lines). Beside each
diagram we indicate the Lie algebra $\lih$ of the set of fixed points, the type of the restricted
root system and the isogeny type of the group, in particular `SC' means `simply connected' and
`ADJ' means `adjoint'.

\vskip 0.5cm
\newpage

%
%

\begin{center}
\begin{tabular}{|c|c|c|}
\hline ${\scriptstyle\goh}$ & ${\scriptstyle\textrm{extended Satake diagram}}$ & ${\scriptstyle\textrm{roots and isogeny type}}$\\

\hline
& & \\

${\scriptstyle\gos\goo(\ell+1)}$ &  $\xymatrix{
      *{\overset{\mathstrut 0}{\circ}}
           \ar@{:>} @<-0.95ex> []+<0.42ex,0ex>;[r]-<0.42ex,0ex>
    & *{\overset{\mathstrut 1}{\circ}}
           \ar @{-} @<-0.95ex> []+<0.42ex,0ex>;[r]-<0.42ex,0ex>
    & *{\overset{\mathstrut 2}{\circ}}
           \ar @{-} @<-0.95ex> []+<0.42ex,0ex>;[r]-<1.50ex,0ex>
    & *{\overset{\mathstrut\mbox{}}{\cdots}}
           \ar @{-} @<-0.95ex> []+<1.50ex,0ex>;[r]-<0.42ex,0ex>
    & *{\overset{\mathstrut \ell-1}{\circ}}
           \ar @{-} @<-0.95ex> []+<0.42ex,0ex>;[r]-<0.42ex,0ex>
    & *{\overset{\mathstrut \ell}{\circ}}\\
    &
    &
    &
    & *{\overset{\mathstrut 0}{\circ}}
           \ar @{.} @<-0.40ex> []+<0.32ex,0ex>;[r]-<1.30ex,0ex>
           \ar @{:>} @<-0.95ex> []+<0.32ex,0ex>;[r]-<0.42ex,0ex>
           \ar @{.} @<-1.50ex> []+<0.32ex,0ex>;[r]-<1.30ex,0ex>
    & *{\overset{\mathstrut 1}{\circ}}
} $

&

$ \xymatrix{ {\scriptstyle\ell\geq 1\quad\sfA_\ell \quad
\textrm{SC}}\\
{\scriptstyle\ell=1\quad\sfB_1 \quad \textrm{ADJ}} }$\\
& & \\
\hline

& & \\
$ {\scriptstyle\gos\gop(2\ell+2)}$ & $\xymatrix{
    &  *{\overset{\mathstrut 0}{\circ}}
           \ar @{.} @<-0.1ex>[]-<0ex,1.50ex>;[d]-<0ex,0.42ex> \\
      *{\overset{\mathstrut \mbox{}}{\bullet}}
           \ar @{-} @<-0.95ex> []+<0.42ex,0ex>;[r]-<0.42ex,0ex>
    & *{\overset{\mathstrut \;\;\;2}{\circ}}
           \ar @{-} @<-0.95ex> []+<0.42ex,0ex>;[r]-<0.42ex,0ex>
    & *{\overset{\mathstrut \mbox{}}{\bullet}}
           \ar @{-} @<-0.95ex> []+<0.42ex,0ex>;[r]-<1.50ex,0ex>
    & *{\overset{\mathstrut\mbox{}}{\cdots}}
           \ar @{-} @<-0.95ex> []+<1.50ex,0ex>;[r]-<0.42ex,0ex>
    & *{\overset{\mathstrut 2\ell}{\circ}}
           \ar @{-} @<-0.95ex> []+<0.42ex,0ex>;[r]-<0.42ex,0ex>
    & *{\overset{\mathstrut\mbox{}}{\bullet}}\\
    &  *{\overset{\mathstrut 0}{\circ}}
           \ar @{:>} @<-0.1ex>[]-<0ex,1.50ex>;[d]-<0ex,0.42ex> \\
      *{\overset{\mathstrut \mbox{}}{\bullet}}
           \ar @{-} @<-0.95ex> []+<0.42ex,0ex>;[r]-<0.42ex,0ex>
    & *{\overset{\mathstrut \;\;\;\;\;2}{\circ}}
           \ar @{-} @<-0.95ex> []+<0.42ex,0ex>;[r]-<0.42ex,0ex>
    & *{\overset{\mathstrut \mbox{}}{\bullet}}
    &
    &
    &
} $ & $\xymatrix{\\
    {\scriptstyle\ell\geq 1\quad\sfA_\ell \quad \textrm{SC}}\\
    \\
    {\scriptstyle\ell=1\quad\sfB_1 \quad \textrm{ADJ}}
}$\\
& & \\
\hline

& & \\ ${\scriptstyle\gos\gol(\ell)\oplus\gos\gol(\myRank+1-\ell)}$ & $\xymatrix{
    & *{\overset{\mathstrut 1}{\circ}}
           \ar @{-} @<-0.95ex> []+<0.42ex,0ex>;[r]-<0.42ex,0ex>
    & *{\overset{\mathstrut 2}{\circ}}
           \ar @{-} @<-0.95ex> []+<0.42ex,0ex>;[r]-<1.50ex,0ex>
    & *{\overset{\mathstrut\mbox{}}{\cdots}}
           \ar @{-} @<-0.95ex> []+<1.50ex,0ex>;[r]-<0.42ex,0ex>
    & *{\overset{\mathstrut \ell}{\circ}}
           \ar @{-} @<-0.95ex> []+<0.42ex,0ex>;[r]-<0.42ex,0ex>
    & *{\overset{\mathstrut \mbox{}}{\bullet}}
           \ar@{-}[]-<0ex,1.40ex>;[d]+<0ex,0.50ex>\\
    & \ar@{<->}@/_0.8ex/[]+<0ex,3.50ex>;[dd]-<0ex,3.50ex>
    & \ar@{<->}@/_0.6ex/[]+<0ex,3.50ex>;[dd]-<0ex,3.50ex>
    &
    & \ar@{<->}@/_0.5ex/[]+<0ex,3.50ex>;[dd]-<0ex,3.50ex>
    & \bullet \ar@{-}[]-<0ex,0.40ex>;[d]+<0ex,0.70ex>\\
      {\scriptstyle 0}\;\;\circ\;\;
           \ar@{.} @<-0.30ex> []+<0.50ex,0.50ex>;[uur]
           \ar@{.} @<0.10ex> []+<0.70ex,-0.50ex>;[ddr]+<-0.50ex,1.20ex>
    &
    &
    &
    &
    & \vdots \ar@{-}[]-<0ex,1.80ex>;[d]+<0ex,0.50ex>\\
    &
    &
    &
    &
    & \bullet \ar@{-}[]-<0ex,0.40ex>;[d]+<0ex,1.50ex>\\
    & *{\underset{\mathstrut\myRank}{\circ}}
           \ar @{-} @<1.15ex> []+<0.42ex,0ex>;[r]-<0.42ex,0ex>
    & *{\underset{\mathstrut\myRank-1}{\circ}}
           \ar @{-} @<1.15ex> []+<0.42ex,0ex>;[r]-<1.50ex,0ex>
    & *{\underset{\mathstrut\mbox{}}{\cdots}}
           \ar @{-} @<1.15ex> []+<1.50ex,0ex>;[r]-<0.42ex,0ex>
    & *{\underset{\mathstrut\myRank+1-\ell}{\circ}}
           \ar @{-} @<1.15ex> []+<0.42ex,0ex>;[r]-<0.42ex,0ex>
    & *{\underset{\mathstrut \mbox{}}{\bullet}}
}$ & ${\scriptstyle\sfBC_\ell \quad \textrm{SC=ADJ}} $\\
& & \\
\hline

& & \\
${\scriptstyle\gos\gol(\ell)\oplus\gos\gol(\ell)\oplus\mathbb{C}}$ & $\xymatrix{
    & *{\overset{\mathstrut 1}{\circ}}
           \ar @{-} @<-0.95ex> []+<0.42ex,0ex>;[r]-<0.42ex,0ex>
           \ar@{<->}@/_0.8ex/[]+<0ex,-2.0ex>;[dd]-<0ex,-2.0ex>
    & *{\overset{\mathstrut 2}{\circ}}
           \ar @{-} @<-0.95ex> []+<0.42ex,0ex>;[r]-<1.50ex,0ex>
           \ar@{<->}@/_0.8ex/[]+<0ex,-2.0ex>;[dd]-<0ex,-2.0ex>
    & *{\overset{\mathstrut\mbox{}}{\cdots}}
           \ar @{-} @<-0.95ex> []+<1.50ex,0ex>;[r]-<0.42ex,0ex>
    & *{\overset{\mathstrut \ell-1}{\circ}}
           \ar@{<->}@/_0.8ex/[]+<0ex,-2.0ex>;[dd]-<0ex,-2.0ex>\\
      {\scriptstyle0}\;\;\circ\;\;
           \ar@{.} @<-0.50ex> []+<0.40ex,0.70ex>;[ur]+<-0.70ex,-0.70ex>
           \ar@{.} @<0.10ex> []+<0.70ex,-0.50ex>;[dr]+<-0.50ex,1.20ex>
    &
    &
    &
    &
    & \circ\;\;{\scriptstyle\ell}
           \ar@{-}[]+<-1.3ex,0.3ex>;[ul]+<0.2ex,-1.3ex>
           \ar@{-}[]+<-1.3ex,-0.3ex>;[dl]+<0.3ex,1.5ex>\\
    & *{\underset{\mathstrut 2\ell-1}{\circ}}
           \ar @{-} @<1.15ex> []+<0.42ex,0ex>;[r]-<0.42ex,0ex>
    & *{\underset{\mathstrut 2\ell-2}{\circ}}
           \ar @{-} @<1.15ex> []+<0.42ex,0ex>;[r]-<1.50ex,0ex>
    & *{\underset{\mathstrut\mbox{}}{\cdots}}
           \ar @{-} @<1.15ex> []+<1.50ex,0ex>;[r]-<0.42ex,0ex>
    & *{\underset{\mathstrut\ell+1}{\circ}}
} $ & ${\scriptstyle\sfC_\ell \quad \textrm{SC}}$\\
& & \\
\hline
\end{tabular}
\end{center}

\newpage

\begin{center}
\begin{tabular}{|c|c|c|}
\hline ${\scriptstyle\goh}$ & ${\scriptstyle\textrm{extended Satake diagram}}$ & ${\scriptstyle\textrm{roots and isogeny type}}$\\

\hline
& & \\
${\scriptstyle\gos\goo(\ell)\oplus\gos\goo(2\myRank+1-\ell)}$ & $\xymatrix{
      *{\overset{\mathstrut 0}{\circ}}
           \ar @{:>} @<-0.95ex> []+<0.42ex,0ex>;[r]-<0.42ex,0ex>
    & *{\overset{\mathstrut 1}{\circ}}
           \ar @{-} @<-0.95ex> []+<0.42ex,0ex>;[r]-<0.42ex,0ex>
    & *{\overset{\mathstrut 2}{\circ}}
           \ar @{-} @<-0.95ex> []+<0.42ex,0ex>;[r]-<1.50ex,0ex>
    & *{\overset{\mathstrut\mbox{}}{\cdots}}
           \ar @{-} @<-0.95ex> []+<1.50ex,0ex>;[r]-<0.42ex,0ex>
    & *{\overset{\mathstrut \ell}{\circ}}
           \ar @{-} @<-0.95ex> []+<0.42ex,0ex>;[r]-<0.42ex,0ex>
    & *{\overset{\mathstrut \mbox{}}{\bullet}}
           \ar @{-} @<-0.95ex> []+<0.42ex,0ex>;[r]-<1.50ex,0ex>
    & *{\overset{\mathstrut\mbox{}}{\cdots}}
           \ar @{-} @<-0.95ex> []+<1.50ex,0ex>;[r]-<0.42ex,0ex>
    & *{\overset{\mathstrut\mbox{}}{\bullet}}
           \ar @{=>} @<-0.95ex> []+<0.42ex,0ex>;[r]-<0.42ex,0ex>
    & *{\overset{\mathstrut\mbox{}}{\bullet}}\\
    &
    &
    &
    &  *{\overset{\mathstrut 0}{\circ}}
           \ar @{.} @<-0.1ex>[]-<0ex,1.50ex>;[d]-<0ex,0.42ex>\\
    &
    &
    & *{\overset{\mathstrut 1}{\circ}}
           \ar @{-} @<-0.95ex> []+<0.42ex,0ex>;[r]-<0.42ex,0ex>
    & *{\overset{\mathstrut \;\;\;2}{\circ}}
           \ar @{-} @<-0.95ex> []+<0.42ex,0ex>;[r]-<0.42ex,0ex>
    & *{\overset{\mathstrut \mbox{}}{\bullet}}
           \ar @{-} @<-0.95ex> []+<0.42ex,0ex>;[r]-<1.50ex,0ex>
    & *{\overset{\mathstrut\mbox{}}{\cdots}}
           \ar @{-} @<-0.95ex> []+<1.50ex,0ex>;[r]-<0.42ex,0ex>
    & *{\overset{\mathstrut\mbox{}}{\bullet}}
           \ar @{=>} @<-0.95ex> []+<0.42ex,0ex>;[r]-<0.42ex,0ex>
    & *{\overset{\mathstrut\mbox{}}{\bullet}}\\
    &
    &
    & *{\overset{\mathstrut 0}{\circ}}
           \ar @{.} @<-0.95ex> []+<0.42ex,0ex>;[r]-<0.42ex,0ex>
    & *{\overset{\mathstrut 1}{\circ}}
           \ar @{-} @<-0.95ex> []+<0.42ex,0ex>;[r]-<0.42ex,0ex>
    & *{\overset{\mathstrut \mbox{}}{\bullet}}
           \ar @{-} @<-0.95ex> []+<0.42ex,0ex>;[r]-<1.50ex,0ex>
    & *{\overset{\mathstrut\mbox{}}{\cdots}}
           \ar @{-} @<-0.95ex> []+<1.50ex,0ex>;[r]-<0.42ex,0ex>
    & *{\overset{\mathstrut\mbox{}}{\bullet}}
           \ar @{=>} @<-0.95ex> []+<0.42ex,0ex>;[r]-<0.42ex,0ex>
    & *{\overset{\mathstrut\mbox{}}{\bullet}}\\
} $ & $ \xymatrix{ {\scriptstyle\ell\geq1\quad\sfB_\ell\quad\textrm{ADJ}}\\
                   \\
                   {\scriptstyle\ell=2\quad\sfC_2\quad\textrm{SC}}\\
                   {\scriptstyle\ell=1\quad\sfA_1\quad\textrm{SC}}}$\\
& & \\
\hline

& &\\
${\scriptstyle\gog\gol(\ell)}$ & $\xymatrix{
      *{\overset{\mathstrut 0}{\circ}}
           \ar @{:>} @<-0.95ex> []+<0.42ex,0ex>;[r]-<0.42ex,0ex>
    & *{\overset{\mathstrut 1}{\circ}}
           \ar @{-} @<-0.95ex> []+<0.42ex,0ex>;[r]-<0.42ex,0ex>
    & *{\overset{\mathstrut 2}{\circ}}
           \ar @{-} @<-0.95ex> []+<0.42ex,0ex>;[r]-<1.50ex,0ex>
    & *{\overset{\mathstrut\mbox{}}{\cdots}}
           \ar @{-} @<-0.95ex> []+<1.50ex,0ex>;[r]-<0.42ex,0ex>
    & *{\overset{\mathstrut \ell-1}{\circ}}
           \ar @{<=} @<-0.95ex> []+<0.42ex,0ex>;[r]-<0.42ex,0ex>
    & *{\overset{\mathstrut \ell}{\circ}}
    &\\
    &
    &
    &
    & *{\overset{\mathstrut 1}{\circ}}
           \ar @{<=} @<-0.95ex> []+<0.42ex,0ex>;[r]-<0.42ex,0ex>
    & *{\overset{\mathstrut 2}{\circ}}
           \ar @{<:} @<-0.95ex> []+<0.42ex,0ex>;[r]-<0.42ex,0ex>
    & *{\overset{\mathstrut 0}{\circ}}\\
} $ & $\xymatrix{ {\scriptstyle\ell\geq 2\quad\sfC_\ell\quad\textrm{SC}}\\
                  {\scriptstyle\ell=2\quad\sfB_2\quad\textrm{ADJ}}}$\\
& & \\
\hline

& & \\
${\scriptstyle\gos\gop(2\ell)\oplus\gos\gop(2\myRank-2\ell)}$ & $\xymatrix{
    &  *{\overset{\mathstrut 0}{\circ}}
           \ar @{.} @<-0.1ex>[]-<0ex,1.50ex>;[d]-<0ex,0.42ex> \\
      *{\overset{\mathstrut \mbox{}}{\bullet}}
           \ar @{-} @<-0.95ex> []+<0.42ex,0ex>;[r]-<0.42ex,0ex>
    & *{\overset{\mathstrut \;\;\;2}{\circ}}
           \ar @{-} @<-0.95ex> []+<0.42ex,0ex>;[r]-<0.42ex,0ex>
    & *{\overset{\mathstrut \mbox{}}{\bullet}}
           \ar @{-} @<-0.95ex> []+<0.42ex,0ex>;[r]-<1.50ex,0ex>
    & *{\overset{\mathstrut\mbox{}}{\cdots}}
           \ar @{-} @<-0.95ex> []+<1.50ex,0ex>;[r]-<0.42ex,0ex>
    & *{\overset{\mathstrut 2\ell}{\circ}}
           \ar @{-} @<-0.95ex> []+<0.42ex,0ex>;[r]-<0.42ex,0ex>
    & *{\overset{\mathstrut \mbox{}}{\bullet}}
           \ar @{-} @<-0.95ex> []+<0.42ex,0ex>;[r]-<1.50ex,0ex>
    & *{\overset{\mathstrut\mbox{}}{\cdots}}
           \ar @{-} @<-0.95ex> []+<1.50ex,0ex>;[r]-<0.42ex,0ex>
    & *{\overset{\mathstrut\mbox{}}{\bullet}}
           \ar @{<=} @<-0.95ex> []+<0.42ex,0ex>;[r]-<0.42ex,0ex>
    & *{\overset{\mathstrut\mbox{}}{\bullet}}\\
}$ & ${\scriptstyle\sfBC_\ell\quad\textrm{SC=ADJ}}$\\
& & \\
\hline

& & \\
${\scriptstyle\gos\gop(2\ell)\oplus\gos\gop(2\ell)}$ & $\xymatrix{
    &  *{\overset{\mathstrut 0}{\circ}}
           \ar @{.} @<-0.1ex>[]-<0ex,1.50ex>;[d]-<0ex,0.42ex> \\
      *{\overset{\mathstrut \mbox{}}{\bullet}}
           \ar @{-} @<-0.95ex> []+<0.42ex,0ex>;[r]-<0.42ex,0ex>
    & *{\overset{\mathstrut \;\;\;2}{\circ}}
           \ar @{-} @<-0.95ex> []+<0.42ex,0ex>;[r]-<0.42ex,0ex>
    & *{\overset{\mathstrut \mbox{}}{\bullet}}
           \ar @{-} @<-0.95ex> []+<0.42ex,0ex>;[r]-<1.50ex,0ex>
    & *{\overset{\mathstrut\mbox{}}{\cdots}}
           \ar @{-} @<-0.95ex> []+<1.50ex,0ex>;[r]-<0.42ex,0ex>
    & *{\overset{\mathstrut 2\ell-2}{\circ}}
           \ar @{-} @<-0.95ex> []+<0.42ex,0ex>;[r]-<0.42ex,0ex>
    & *{\overset{\mathstrut\mbox{}}{\bullet}}
           \ar @{<=} @<-0.95ex> []+<0.42ex,0ex>;[r]-<0.42ex,0ex>
    & *{\overset{\mathstrut 2\ell}{\circ}}
    &\\
    &
    &
    &
    &
    & *{\overset{\mathstrut\mbox{}}{\bullet}}
           \ar @{<=} @<-0.95ex> []+<0.42ex,0ex>;[r]-<0.42ex,0ex>
    & *{\overset{\mathstrut 2}{\circ}}
           \ar @{<:} @<-0.95ex> []+<0.42ex,0ex>;[r]-<0.42ex,0ex>
    & *{\overset{\mathstrut 0}{\circ}}\\
    &
    &
    &
      *{\overset{\mathstrut \mbox{}}{\bullet}}
           \ar @{-} @<-0.95ex> []+<0.42ex,0ex>;[r]-<0.42ex,0ex>
    & *{\overset{\mathstrut 2}{\circ}}
           \ar @{-} @<-0.95ex> []+<0.42ex,0ex>;[r]-<0.42ex,0ex>
    & *{\overset{\mathstrut\mbox{}}{\bullet}}
           \ar @{<=} @<-0.95ex> []+<0.42ex,0ex>;[r]-<0.42ex,0ex>
    & *{\overset{\mathstrut 4}{\circ}}
           \ar @{.} @<-0.95ex> []+<0.42ex,0ex>;[r]-<0.42ex,0ex>
    & *{\overset{\mathstrut 0}{\circ}}\\
} $ & $\xymatrix{ \\
                  {\scriptstyle\ell\geq1\quad\sfC_\ell\quad\textrm{SC}}\\
                  {\scriptstyle\ell=1\quad\sfB_1\quad\textrm{ADJ}}\\
                  {\scriptstyle\ell=2\quad\sfB_2\quad\textrm{ADJ}}}$\\
& & \\
\hline
\end{tabular}
\end{center}

\newpage

\begin{center}
\begin{tabular}{|c|c|c|}
\hline ${\scriptstyle\goh}$ & ${\scriptstyle\textrm{extended Satake diagram}}$ & ${\scriptstyle\textrm{roots and isogeny type}}$\\

\hline
& & \\
${\scriptstyle\gos\goo(\ell)\oplus\gos\goo(2\myRank-\ell)}$ & $\xymatrix{
    &
    &
    &
    &
    &
    &
    &
    & *{\overset{\mathstrut\mbox{}}{\bullet}}\\
      *{\overset{\mathstrut 0}{\circ}}
           \ar @{:>} @<-0.95ex> []+<0.42ex,0ex>;[r]-<0.42ex,0ex>
    & *{\overset{\mathstrut 1}{\circ}}
           \ar @{-} @<-0.95ex> []+<0.42ex,0ex>;[r]-<0.42ex,0ex>
    & *{\overset{\mathstrut 2}{\circ}}
           \ar @{-} @<-0.95ex> []+<0.42ex,0ex>;[r]-<1.50ex,0ex>
    & *{\overset{\mathstrut\mbox{}}{\cdots}}
           \ar @{-} @<-0.95ex> []+<1.50ex,0ex>;[r]-<0.42ex,0ex>
    & *{\overset{\mathstrut \ell}{\circ}}
           \ar @{-} @<-0.95ex> []+<0.42ex,0ex>;[r]-<0.42ex,0ex>
    & *{\overset{\mathstrut \mbox{}}{\bullet}}
           \ar @{-} @<-0.95ex> []+<0.42ex,0ex>;[r]-<1.50ex,0ex>
    & *{\overset{\mathstrut\mbox{}}{\cdots}}
           \ar @{-} @<-0.95ex> []+<1.50ex,0ex>;[r]-<0.42ex,0ex>
    & *{\overset{\mathstrut\mbox{}}{\bullet}}
           \ar @{-} []+<0.2ex,-0.7ex>;[ru]+<0ex,-1.0ex>
           \ar @{-} []+<0.2ex,-1.2ex>;[rd]+<0ex,-1.0ex>
    &\\
    &
    &
    &
    &
    &
    &
    &
    & *{\overset{\mathstrut\mbox{}}{\bullet}}\\
    &
    &
    &
    &
    &
    &
    &
    & *{\overset{\mathstrut\mbox{}}{\bullet}}\\
    &
    &
    & *{\overset{\mathstrut 0}{\circ}}
           \ar @{:>} @<-0.95ex> []+<0.42ex,0ex>;[r]-<0.42ex,0ex>
    & *{\overset{\mathstrut 1}{\circ}}
           \ar @{-} @<-0.95ex> []+<0.42ex,0ex>;[r]-<0.42ex,0ex>
    & *{\overset{\mathstrut \mbox{}}{\bullet}}
           \ar @{-} @<-0.95ex> []+<0.42ex,0ex>;[r]-<1.50ex,0ex>
    & *{\overset{\mathstrut\mbox{}}{\cdots}}
           \ar @{-} @<-0.95ex> []+<1.50ex,0ex>;[r]-<0.42ex,0ex>
    & *{\overset{\mathstrut\mbox{}}{\bullet}}
           \ar @{-} []+<0.2ex,-0.7ex>;[ru]+<0ex,-1.0ex>
           \ar @{-} []+<0.2ex,-1.2ex>;[rd]+<0ex,-1.0ex>
    &\\
    &
    &
    &
    &
    &
    &
    &
    & *{\overset{\mathstrut\mbox{}}{\bullet}}\\
    &
    &
    &
    & *{\overset{\mathstrut 0}{\circ}}
           \ar @{.} ;[d]+<0ex,-0.5ex>
    &
    &
    &
    & *{\overset{\mathstrut\mbox{}}{\bullet}}\\
    &
    &
    & *{\overset{\mathstrut 1}{\circ}}
           \ar @{-} @<-0.95ex> []+<0.42ex,0ex>;[r]-<0.42ex,0ex>
    & *{\overset{\mathstrut \;\;\;2}{\circ}}
           \ar @{-} @<-0.95ex> []+<0.42ex,0ex>;[r]-<0.42ex,0ex>
    & *{\overset{\mathstrut \mbox{}}{\bullet}}
           \ar @{-} @<-0.95ex> []+<0.42ex,0ex>;[r]-<1.50ex,0ex>
    & *{\overset{\mathstrut\mbox{}}{\cdots}}
           \ar @{-} @<-0.95ex> []+<1.50ex,0ex>;[r]-<0.42ex,0ex>
    & *{\overset{\mathstrut\mbox{}}{\bullet}}
           \ar @{-} []+<0.2ex,-0.7ex>;[ru]+<0ex,-1.0ex>
           \ar @{-} []+<0.2ex,-1.2ex>;[rd]+<0ex,-1.0ex>
    &\\
    &
    &
    &
    &
    &
    &
    &
    & *{\overset{\mathstrut\mbox{}}{\bullet}}\\
}$ & $\xymatrix{\\
                {\scriptstyle\ell\geq 1\quad\sfB_\ell\quad\textrm{ADJ}}\\
                \\
                \\
                {\scriptstyle\ell=1\quad\sfA_1\quad\textrm{SC}}\\
                \\
                \\
                {\scriptstyle\ell=2\quad\sfC_2\quad\textrm{SC}}}$ \\
& & \\
\hline

& & \\
${\scriptstyle\gos\goo(\ell)\oplus\gos\goo(\ell+2)}$ & $\xymatrix{
    &
    &
    &
    &
    & *{\overset{\mathstrut\mbox{}}{\circ\;\;{\scriptstyle\ell}}}
           \ar@{<->}@<-0.5ex> @/^0.8ex/[]+<0ex,-2.0ex>;[dd]-<0ex,0ex>\\
      *{\overset{\mathstrut 0}{\circ}}
           \ar @{:>} @<-0.95ex> []+<0.42ex,0ex>;[r]-<0.42ex,0ex>
    & *{\overset{\mathstrut 1}{\circ}}
           \ar @{-} @<-0.95ex> []+<0.42ex,0ex>;[r]-<0.42ex,0ex>
    & *{\overset{\mathstrut 2}{\circ}}
           \ar @{-} @<-0.95ex> []+<0.42ex,0ex>;[r]-<1.50ex,0ex>
    & *{\overset{\mathstrut\mbox{}}{\cdots}}
           \ar @{-} @<-0.95ex> []+<1.50ex,0ex>;[r]-<0.42ex,0ex>
    & *{\overset{\mathstrut \ell-2\;\;}{\circ}}
           \ar @{-} []+<0.3ex,-0.7ex>;[ru]+<-1.35ex,-1.1ex>
           \ar @{-} []+<0.3ex,-1.2ex>;[rd]+<-1.4ex,-0.6ex>
    &\\
    &
    &
    &
    &
    & *{\overset{\mathstrut\mbox{}}{\;\;\;\circ\;\;{\scriptstyle\ell+1}}}
} $ & ${\scriptstyle\sfB_\ell\quad\textrm{ADJ}}$\\
& & \\
\hline
\end{tabular}
\end{center}

\newpage

\begin{center}
\begin{tabular}{|c|c|c|}
\hline ${\scriptstyle\goh}$ & ${\scriptstyle\textrm{extended Satake diagram}}$ & ${\scriptstyle\textrm{roots and isogeny type}}$\\

\hline
& & \\
${\scriptstyle\gog\gol(2\ell)}$ & $\xymatrix{
    & *{\overset{\mathstrut 0}{\circ}}
           \ar @{.} ;[d]+<0ex,-0.5ex>
    &
    &
    &
    & *{\overset{\mathstrut\mbox{}}{\bullet\;\;\;}}\\
      *{\overset{\mathstrut\mbox{}}{\bullet}}
           \ar @{-} @<-0.95ex> []+<0.42ex,0ex>;[r]-<0.42ex,0ex>
    & *{\overset{\mathstrut\;\;\;2}{\circ}}
           \ar @{-} @<-0.95ex> []+<0.42ex,0ex>;[r]-<0.42ex,0ex>
    & *{\overset{\mathstrut\mbox{}}{\bullet}}
           \ar @{-} @<-0.95ex> []+<0.42ex,0ex>;[r]-<1.50ex,0ex>
    & *{\overset{\mathstrut\mbox{}}{\cdots}}
           \ar @{-} @<-0.95ex> []+<1.50ex,0ex>;[r]-<0.42ex,0ex>
    & *{\overset{\mathstrut 2\ell-2\;\;\;}{\circ}}
           \ar @{-} []+<0.3ex,-0.7ex>;[ru]+<-1.35ex,-1.1ex>
           \ar @{-} []+<0.3ex,-1.2ex>;[rd]+<-1.4ex,-0.6ex>
    &
    &\\
    &
    &
    &
    &
    & *{\overset{\mathstrut\mbox{}}{\;\circ\;\;{\scriptstyle 2\ell}}}\\
    &
    &
    &
    &
    & *{\overset{\mathstrut\mbox{}}{\bullet\;\;\;\;\,}}\\
    &
    &
    & *{\overset{\mathstrut\mbox{}}{\bullet}}
           \ar @{-} @<-0.95ex> []+<0.42ex,0ex>;[r]-<0.42ex,0ex>
    & *{\overset{\mathstrut 2}{\circ}}
           \ar @{-} []+<0.3ex,-0.7ex>;[ru]+<-1.7ex,-1.2ex>
           \ar @{-} []+<0.3ex,-1.2ex>;[rd]+<-1.7ex,-0.6ex>
    &
    &\\
    &
    &
    &
    &
    & *{\overset{\mathstrut 4\;\;\;\;}{\circ\;\;\;\;\,}}
           \ar @{<:} @<-0.95ex> []+<-1.00ex,0ex>;[r]-<0.42ex,0ex>
    & *{\overset{\mathstrut 0}{\circ}}
} $ & $\xymatrix{\\
                 {\scriptstyle\ell\geq 2\quad\sfC_\ell\quad\textrm{SC}}\\
                 \\
                 \\
                 {\scriptstyle\ell=2\quad\sfB_2\quad\textrm{ADJ}}}$ \\
& & \\
\hline

& & \\
${\scriptstyle\gog\gol(2\ell+1)}$ & $\xymatrix{
    & *{\overset{\mathstrut 0}{\circ}}
           \ar @{.} ;[d]+<0ex,-0.5ex>
    &
    &
    &
    &
    & *{\overset{\mathstrut\mbox{}}{\circ\;\;{\scriptstyle 2\ell}}}
           \ar@{<->}@<-0.7ex> @/^0.8ex/[]+<0ex,-2.0ex>;[dd]-<0ex,0ex>\\
      *{\overset{\mathstrut\mbox{}}{\bullet}}
           \ar @{-} @<-0.95ex> []+<0.42ex,0ex>;[r]-<0.42ex,0ex>
    & *{\overset{\mathstrut\;\;\;2}{\circ}}
           \ar @{-} @<-0.95ex> []+<0.42ex,0ex>;[r]-<0.42ex,0ex>
    & *{\overset{\mathstrut\mbox{}}{\bullet}}
           \ar @{-} @<-0.95ex> []+<0.42ex,0ex>;[r]-<1.50ex,0ex>
    & *{\overset{\mathstrut\mbox{}}{\cdots}}
           \ar @{-} @<-0.95ex> []+<1.50ex,0ex>;[r]-<0.42ex,0ex>
    & *{\overset{\mathstrut 2\ell-2}{\circ}}
           \ar @{-} @<-0.95ex> []+<0.42ex,0ex>;[r]-<0.42ex,0ex>
    & *{\overset{\mathstrut\mbox{}}{\bullet}}
           \ar @{-} []+<0.3ex,-0.7ex>;[ru]+<-1.7ex,-1.2ex>
           \ar @{-} []+<0.3ex,-1.2ex>;[rd]+<-1.7ex,-0.6ex>
    &\\
    &
    &
    &
    &
    &
    & *{\overset{\mathstrut\mbox{}}{\;\;\;\circ\;\;{\scriptstyle 2\ell+1}}}
} $ & ${\scriptstyle\sfBC_\ell\quad\textrm{SC=ADJ}}$\\
& & \\
\hline

& & \\
${\scriptstyle\gos\goo(10)\oplus\mathbb{C}}$ & $\xymatrix{
    &
    & *{\overset{\mathstrut 0}{\circ}}
           \ar @{.} ;[d]+<0ex,-0.5ex>\\
    &
    & *{\overset{\mathstrut\mbox{}}{\,\,\,\,\circ\;{\scriptstyle 2}}}
           \ar @{-} []+<0ex,-1.4ex> ;[d]+<0ex,-0.5ex>
    &
    &\\
      *{\overset{\mathstrut 1}{\circ}}
           \ar @{-} @<-0.95ex> []+<0.42ex,0ex>;[r]-<0.42ex,0ex>
    & *{\overset{\mathstrut 3}{\bullet}}
           \ar @{-} @<-0.95ex> []+<0.42ex,0ex>;[r]-<0.42ex,0ex>
    & *{\overset{\mathstrut\;\;\;4}{\bullet}}
           \ar @{-} @<-0.95ex> []+<0.42ex,0ex>;[r]-<0.42ex,0ex>
    & *{\overset{\mathstrut 5}{\bullet}}
           \ar @{-} @<-0.95ex> []+<0.42ex,0ex>;[r]-<0.42ex,0ex>
    & *{\overset{\mathstrut 6}{\circ}}\\
      \mbox{}
           \ar @{<->} @/_3.0ex/ @<5.0ex>;[rrrr]
    &
    &
    &
    & \mbox{}
}$ & ${\scriptstyle\sfBC_2\quad\textrm{SC=ADJ}}$\\
& & \\
\hline

& & \\
${\scriptstyle\sfF_4}$ & $\xymatrix{
    &
    &
    & *{\overset{\mathstrut\mbox{}}{\,\,\,\,\circ\;{\scriptstyle 2}}}
           \ar @{-} []+<0ex,-1.4ex> ;[d]+<0ex,-0.5ex>\\
      *{\overset{\mathstrut 0}{\circ}}
           \ar @{.} @<-0.95ex> []+<0.42ex,0ex>;[r]-<0.42ex,0ex>
    & *{\overset{\mathstrut 1}{\circ}}
           \ar @{-} @<-0.95ex> []+<0.42ex,0ex>;[r]-<0.42ex,0ex>
    & *{\overset{\mathstrut 3}{\bullet}}
           \ar @{-} @<-0.95ex> []+<0.42ex,0ex>;[r]-<0.42ex,0ex>
    & *{\overset{\mathstrut\;\;\;4}{\bullet}}
           \ar @{-} @<-0.95ex> []+<0.42ex,0ex>;[r]-<0.42ex,0ex>
    & *{\overset{\mathstrut 5}{\bullet}}
           \ar @{-} @<-0.95ex> []+<0.42ex,0ex>;[r]-<0.42ex,0ex>
    & *{\overset{\mathstrut 6}{\circ}}\\
}$ & ${\scriptstyle\sfA_2\quad\textrm{SC}}$\\
& & \\
\hline
\end{tabular}
\end{center}

\newpage

\begin{center}
\begin{tabular}{|c|c|c|}
\hline ${\scriptstyle\goh}$ & ${\scriptstyle\textrm{extended Satake diagram}}$ & ${\scriptstyle\textrm{roots and isogeny type}}$\\

\hline

& & \\
${\scriptstyle\sfE_6\oplus\mathbb{C}}$ & $\xymatrix{
    &
    &
    & *{\overset{\mathstrut\mbox{}}{\,\,\,\,\circ\;{\scriptstyle 2}}}
           \ar @{-} []+<0ex,-1.4ex> ;[d]+<0ex,-0.5ex>\\
      *{\overset{\mathstrut 0}{\circ}}
           \ar @{.} @<-0.95ex> []+<0.42ex,0ex>;[r]-<0.42ex,0ex>
    & *{\overset{\mathstrut 1}{\circ}}
           \ar @{-} @<-0.95ex> []+<0.42ex,0ex>;[r]-<0.42ex,0ex>
    & *{\overset{\mathstrut 3}{\bullet}}
           \ar @{-} @<-0.95ex> []+<0.42ex,0ex>;[r]-<0.42ex,0ex>
    & *{\overset{\mathstrut\;\;\;4}{\bullet}}
           \ar @{-} @<-0.95ex> []+<0.42ex,0ex>;[r]-<0.42ex,0ex>
    & *{\overset{\mathstrut 5}{\bullet}}
           \ar @{-} @<-0.95ex> []+<0.42ex,0ex>;[r]-<0.42ex,0ex>
    & *{\overset{\mathstrut 6}{\circ}}
           \ar @{-} @<-0.95ex> []+<0.42ex,0ex>;[r]-<0.42ex,0ex>
    & *{\overset{\mathstrut 7}{\circ}}\\
}$ & ${\scriptstyle\sfC_3\quad\textrm{SC}}$\\
& & \\
\hline

& & \\
${\scriptstyle\gos\goo(9)}$ & $\xymatrix@1{
      *{\overset{\mathstrut 1}{\bullet}}
           \ar @{-} @<-0.95ex> []+<0.42ex,0ex>;[r]-<0.42ex,0ex>
    & *{\overset{\mathstrut 2}{\bullet}}
           \ar @{=>} @<-0.95ex> []+<0.32ex,0ex>;[r]-<0.42ex,0ex>
    & *{\overset{\mathstrut 3}{\bullet}}
           \ar @{-} @<-0.95ex> []+<0.42ex,0ex>;[r]-<0.42ex,0ex>
    & *{\overset{\mathstrut 4}{\circ}}
           \ar @{.} @<-0.95ex> []+<0.42ex,0ex>;[r]-<0.42ex,0ex>
    & *{\overset{\mathstrut 0}{\circ}}
} $ & ${\scriptstyle\sfBC_1\quad\textrm{SC=ADJ}}$\\
& & \\
\hline

& & \\
${\scriptstyle\gos\gol(\ell+1)}$ & $\xymatrix{
    & *{\overset{\mathstrut 1}{\circ}}
           \ar @{-} @<-0.95ex> []+<0.42ex,0ex>;[r]-<0.42ex,0ex>
           \ar@{<->}@/^0.8ex/[]+<0ex,-2.0ex>;[dd]-<0ex,-2.0ex>
    & *{\overset{\mathstrut 2}{\circ}}
           \ar @{-} @<-0.95ex> []+<0.42ex,0ex>;[r]-<1.50ex,0ex>
           \ar@{<->}@/^0.8ex/[]+<0ex,-2.0ex>;[dd]-<0ex,-2.0ex>
    & *{\overset{\mathstrut\mbox{}}{\cdots}}
           \ar @{-} @<-0.95ex> []+<1.50ex,0ex>;[r]-<0.42ex,0ex>
    & *{\overset{\mathstrut \ell}{\circ}}
           \ar@{<->}@/^0.8ex/[]+<0ex,-2.0ex>;[dd]-<0ex,-2.0ex>\\
      {\scriptstyle0}\;\;\circ\;\;
           \ar@{.} @<-0.50ex> []+<0.40ex,0.70ex>;[ur]+<-0.70ex,-0.70ex>
           \ar@{.} @<0.10ex> []+<0.70ex,-0.50ex>;[dr]+<-0.50ex,1.20ex>
    &
    &
    &
    &\\
    & *{\underset{\mathstrut\mbox{}}{\circ}}
           \ar @{-} @<1.15ex> []+<0.42ex,0ex>;[r]-<0.42ex,0ex>
    & *{\underset{\mathstrut\mbox{}}{\circ}}
           \ar @{-} @<1.15ex> []+<0.42ex,0ex>;[r]-<1.50ex,0ex>
    & *{\underset{\mathstrut\mbox{}}{\cdots}}
           \ar @{-} @<1.15ex> []+<1.50ex,0ex>;[r]-<0.42ex,0ex>
    & *{\underset{\mathstrut\mbox{}}{\circ}}\\
    &
    &
    &
    & *{\overset{\mathstrut 1}{\circ}}
           \ar@{<->}@/^0.8ex/[]+<0ex,-2.0ex>;[dd]-<0ex,-2.0ex>\\
    &
    &
    & {\scriptstyle0}\;\;\circ\;\;
           \ar@{:>} @<-0.50ex> []+<0.40ex,0.70ex>;[ur]+<-0.70ex,-0.70ex>
           \ar@{:>} @<0.10ex> []+<0.55ex,-0.35ex>;[dr]+<-0.50ex,1.20ex>
    &\\
    &
    &
    &
    & *{\underset{\mathstrut\mbox{}}{\circ}}
}$ & $\xymatrix{\\
                {\scriptstyle\ell\geq1\quad\sfA_\ell \quad \textrm{SC}}\\
                \\
                \\
                {\scriptstyle\ell=1\quad\sfB_1 \quad \textrm{ADJ}}
                }$\\
& & \\
\hline

& & \\
${\scriptstyle\gos\goo(2\ell+1)}$ & $\xymatrix{
    & *{\overset{\mathstrut 1}{\circ}}
           \ar @{-} @<-0.95ex> []+<0.42ex,0ex>;[r]-<0.42ex,0ex>
           \ar@{<->}@/^0.8ex/[]+<0ex,-2.0ex>;[dd]-<0ex,-2.0ex>
    & *{\overset{\mathstrut 2}{\circ}}
           \ar @{-} @<-0.95ex> []+<0.42ex,0ex>;[r]-<1.50ex,0ex>
           \ar@{<->}@/^0.8ex/[]+<0ex,-2.0ex>;[dd]-<0ex,-2.0ex>
    & *{\overset{\mathstrut\mbox{}}{\cdots}}
           \ar @{-} @<-0.95ex> []+<1.50ex,0ex>;[r]-<0.42ex,0ex>
    & *{\overset{\mathstrut \ell-1}{\circ}}
           \ar @{=>} @<-0.95ex> []+<0.32ex,0ex>;[r]-<0.42ex,0ex>
           \ar@{<->}@/^0.8ex/[]+<0ex,-2.0ex>;[dd]-<0ex,-2.0ex>
    & *{\overset{\mathstrut \ell}{\circ}}
           \ar@{<->}@/^0.8ex/[]+<0ex,-2.0ex>;[dd]-<0ex,-2.0ex>\\
      {\scriptstyle0}\;\;\circ\;\;
           \ar@{.} @<-0.50ex> []+<0.40ex,0.70ex>;[ur]+<-0.70ex,-0.70ex>
           \ar@{.} @<0.10ex> []+<0.70ex,-0.50ex>;[dr]+<-0.50ex,1.20ex>
    &
    &
    &
    &
    &\\
    & *{\underset{\mathstrut\mbox{}}{\circ}}
           \ar @{-} @<1.15ex> []+<0.42ex,0ex>;[r]-<0.42ex,0ex>
    & *{\underset{\mathstrut\mbox{}}{\circ}}
           \ar @{-} @<1.15ex> []+<0.42ex,0ex>;[r]-<1.50ex,0ex>
    & *{\underset{\mathstrut\mbox{}}{\cdots}}
           \ar @{-} @<1.15ex> []+<1.50ex,0ex>;[r]-<0.42ex,0ex>
    & *{\underset{\mathstrut\mbox{}}{\circ}}
           \ar @{=>} @<1.15ex> []+<0.32ex,0ex>;[r]-<0.42ex,0ex>
    & *{\underset{\mathstrut\mbox{}}{\circ}}
} $ & ${\scriptstyle\sfB_\ell \quad \textrm{ADJ}}$\\
& & \\
\hline

& &\\
${\scriptstyle\gos\gop(2\ell)}$ & $\xymatrix{
    & *{\overset{\mathstrut 1}{\circ}}
           \ar @{-} @<-0.95ex> []+<0.42ex,0ex>;[r]-<0.42ex,0ex>
           \ar@{<->}@/^0.8ex/[]+<0ex,-2.0ex>;[dd]-<0ex,-2.0ex>
    & *{\overset{\mathstrut 2}{\circ}}
           \ar @{-} @<-0.95ex> []+<0.42ex,0ex>;[r]-<1.50ex,0ex>
           \ar@{<->}@/^0.8ex/[]+<0ex,-2.0ex>;[dd]-<0ex,-2.0ex>
    & *{\overset{\mathstrut\mbox{}}{\cdots}}
           \ar @{-} @<-0.95ex> []+<1.50ex,0ex>;[r]-<0.42ex,0ex>
    & *{\overset{\mathstrut \ell-1}{\circ}}
           \ar @{<=} @<-0.95ex> []+<0.42ex,0ex>;[r]-<0.32ex,0ex>
           \ar@{<->}@/^0.8ex/[]+<0ex,-2.0ex>;[dd]-<0ex,-2.0ex>
    & *{\overset{\mathstrut \ell}{\circ}}
           \ar@{<->}@/^0.8ex/[]+<0ex,-2.0ex>;[dd]-<0ex,-2.0ex>\\
      {\scriptstyle0}\;\;\circ\;\;
           \ar@{.} @<-0.50ex> []+<0.40ex,0.70ex>;[ur]+<-0.70ex,-0.70ex>
           \ar@{.} @<0.10ex> []+<0.70ex,-0.50ex>;[dr]+<-0.50ex,1.20ex>
    &
    &
    &
    &
    &\\
    & *{\underset{\mathstrut\mbox{}}{\circ}}
           \ar @{-} @<1.15ex> []+<0.42ex,0ex>;[r]-<0.42ex,0ex>
    & *{\underset{\mathstrut\mbox{}}{\circ}}
           \ar @{-} @<1.15ex> []+<0.42ex,0ex>;[r]-<1.50ex,0ex>
    & *{\underset{\mathstrut\mbox{}}{\cdots}}
           \ar @{-} @<1.15ex> []+<1.50ex,0ex>;[r]-<0.42ex,0ex>
    & *{\underset{\mathstrut\mbox{}}{\circ}}
           \ar @{<=} @<1.15ex> []+<0.42ex,0ex>;[r]-<0.32ex,0ex>
    & *{\underset{\mathstrut\mbox{}}{\circ}}\\
}$ & ${\scriptstyle\sfC_\ell \quad \textrm{SC}}$\\
& & \\
\hline
\end{tabular}
\end{center}


\bibliographystyle{amsplain}

\def\dbar{\leavevmode\hbox to 0pt{\hskip.2ex \accent"16\hss}d}
\providecommand{\bysame}{\leavevmode\hbox to3em{\hrulefill}\thinspace}
\providecommand{\MR}{\relax\ifhmode\unskip\space\fi MR }
\providecommand{\MRhref}[2]{%
  \href{http://www.ams.org/mathscinet-getitem?mr=#1}{#2}
}
\providecommand{\href}[2]{#2}

\end{document}